\documentclass[smallextended,final]{svjour3}
\usepackage{amssymb}
\usepackage[left=2.80cm, right=2.80cm, top=1.60cm, bottom=1.60cm]{geometry}
\usepackage{amsmath}        
\usepackage{amssymb}        
\usepackage{amsfonts}       
\usepackage{tikz}           
\usepackage{newtxtext}      
\usepackage[colorlinks,linkcolor=blue,anchorcolor=blue,citecolor=blue]{hyperref}
\usepackage[numbers,sort&compress]{natbib}  
\smartqed
\newcommand{\sCurrentSectionNum}{0}
\newcommand{\sSectionNum}[2]{\renewcommand{\sCurrentSectionNum}{\csname the#1\endcsname}}

\newcommand{\SectionWithLabel}[1]{\section{#1} \label{#1} \sSectionNum{section}{#1}}
\newcommand{\SubSectionWithLabel}[1]{\subsection{#1} \label{#1} \sSectionNum{subsection}{#1}}

\newcounter{sDefinitionCount}
\newcounter{sRemarkCountt}
\newcounter{sDiscussionCount}
\newcounter{sPropositionCount}

\makeatletter
\@addtoreset{sDefinitionCount}{section}
\@addtoreset{sDefinitionCount}{subsection}
\@addtoreset{sDefinitionCount}{subsubsection}
\@addtoreset{sDefinitionCount}{paragraph}
\@addtoreset{sDefinitionCount}{subparagraph}
\@addtoreset{sRemarkCountt}{section}
\@addtoreset{sRemarkCountt}{subsection}
\@addtoreset{sRemarkCountt}{subsubsection}
\@addtoreset{sRemarkCountt}{paragraph}
\@addtoreset{sRemarkCountt}{subparagraph}
\@addtoreset{sDiscussionCount}{section}
\@addtoreset{sDiscussionCount}{subsection}
\@addtoreset{sDiscussionCount}{subsubsection}
\@addtoreset{sDiscussionCount}{paragraph}
\@addtoreset{sDiscussionCount}{subparagraph}
\@addtoreset{sPropositionCount}{section}
\@addtoreset{sPropositionCount}{subsection}
\@addtoreset{sPropositionCount}{subsubsection}
\@addtoreset{sPropositionCount}{paragraph}
\@addtoreset{sPropositionCount}{subparagraph}
\makeatother
\newcommand{\sLimit} {\mathop{\lim}\limits}
\newcommand{\sProof} {{\noindent\bf Proof.\ }}

\newcommand{\FangKuai}[1] {{\setlength{\parskip}{0.7em} \noindent\bf #1}} 
\newcommand{\sDefinition} {\refstepcounter{sDefinitionCount} \renewcommand{\thesDefinitionCount}{\sCurrentSectionNum.\arabic{sDefinitionCount}} \FangKuai{Definition\;\thesDefinitionCount.\ }}
\newcommand{\sRemark} {\refstepcounter{sRemarkCountt} \renewcommand{\thesRemarkCountt}{\sCurrentSectionNum.\arabic{sRemarkCountt}} \FangKuai{Remark\;\thesRemarkCountt.\ }}
\newcommand{\sDiscussion} {\refstepcounter{sDiscussionCount} \renewcommand{\thesDiscussionCount}{\sCurrentSectionNum.\arabic{sDiscussionCount}} \FangKuai{Discussion\;\thesDiscussionCount.\ }}
\newcommand{\sProposition} {\refstepcounter{sPropositionCount} \renewcommand{\thesPropositionCount}{\sCurrentSectionNum.\arabic{sPropositionCount}} \FangKuai{Proposition\;\thesPropositionCount.\ }}

\newcommand{\DenoteBy} {\triangleq}

\newcommand{\Jiao} {\cap}

\newcommand{\IfAndOnlyIf} {\Leftrightarrow}
\newcommand{\DengJiaYu} {\cong}
\newcommand{\TongTaiYu} {\simeq}

\newcommand{\Summation} {\sum\limits}

\newcommand{\LeftAngle} {\left\langle}
\newcommand{\RightAngle} {\right\rangle}
\newcommand{\LeftGathered} {\left\{\begin{gathered}}
\newcommand{\RightGathered} {\end{gathered}\right.}
\newcommand{\LeftList} {\left.\begin{gathered}}
\newcommand{\RightList} {\end{gathered}\right.}

\title{\boldmath Affine connection representation of gauge fields}
\author{Zhao-Hui Man}

\institute{Zhao-Hui Man, Yinchuan, Ningxia, China. \\ 
           \email{shetcslion@163.com} 
           \at * This article has been published (URL: \url{https://doi.org/10.1155/2022/4573981}) in the journal \emph{Advances in High Energy Physics}, Volume 2022, Article ID 4573981, May 17, 2022, which is identical to this preprint in content.
}

\begin{document} \maketitle 

\begin{abstract}
There are two ways to unify gravitational field and gauge field. One is to represent gravitational field as principal bundle connection, and the other is to represent gauge field as affine connection. Poincar\'{e} gauge theory and metric-affine gauge theory adopt the first approach. This paper adopts the second. In this approach:

(i) Gauge field and gravitational field can both be represented by affine connection; they can be described by a unified spatial frame.

(ii) Time can be regarded as the total metric with respect to all dimensions of internal coordinate space and external coordinate space. On-shell can be regarded as gradient direction. Quantum theory can be regarded as a geometric theory of distribution of gradient directions. Hence, gauge theory, gravitational theory, and quantum theory all reflect intrinsic geometric properties of manifold.

(iii) Coupling constants, chiral asymmetry, PMNS mixing and CKM mixing arise spontaneously as geometric properties in affine connection representation, so they are not necessary to be regarded as direct postulates in the Lagrangian anymore.

(iv) The unification theory of gauge fields that are represented by affine connection can avoid the problem that a proton decays into a lepton in theories such as $SU(5)$. 

(v) There exists a geometric interpretation to the color confinement of quarks.

In the affine connection representation, we can get better interpretations to the above physical properties, therefore, to represent gauge fields by affine connection is probably a necessary step towards the ultimate theory of physics. 
\end{abstract}

\keywords{affine connection representation of gauge fields \and geometric meaning of coupling constant \and time metric \and reference-system \and distribution of gradient directions}
\subclass{Primary 51P05, 70A05, 53Z05 \and Secondary 53C05, 58A05}


\setcounter{tocdepth}{3} \setcounter{secnumdepth}{3}
\tableofcontents
\setcounter{section}{0}

\SectionWithLabel{Introduction}

\SubSectionWithLabel{Background and purpose}

We know that in gauge theory, the field strength and the gauge-covariant derivative 
$$F^a_{\mu\nu}=\partial_\mu A^a_\nu -\partial_\nu A^a_\mu + gf^{abc}A^b_\mu A^c_\nu, \ \ \ \ \ \ \ \ \ \  D_\mu=\partial_\mu-igT^aA_\mu^a$$
both contain a coupling constant $g$, which measures the strength of interaction. A problem is that why is there a coupling constant $g$?

If we represent gauge fields by affine connection, we can obtain a nice interpretation. For example, if we use $\Gamma_{MNP}$ to represent gauge potentials, it is not hard to find some specific conditions to turn the curvature tensor $R_{NPQ}^M$ to 
\begin{equation} \label{EqCurvaturePre}
\begin{aligned}
R_{MNPQ} &= \partial_P \Gamma_{MNQ} - \partial_Q \Gamma_{MNP} + \Gamma_{MHP}\Gamma_{NQ}^H  -\Gamma_{NP}^H\Gamma_{MHQ} \\
&=\partial_P \Gamma_{MNQ} - \partial_Q \Gamma_{MNP} + G^{RH}(\Gamma_{MHP}\Gamma_{RNQ} - \Gamma_{RNP}\Gamma_{MHQ}).\\
\end{aligned}
\end{equation}
Thus, $R_{MNPQ}$ can be used to represent field strength. In addition, for any $\rho_M$, we see that
\begin{equation} \label{EqChargePre}
\rho_{M;P}=\partial_P \rho_M - \Gamma^H_{MP}\rho_H = \partial_P \rho_M - G^{RH}\Gamma_{RMP}\rho_H.
\end{equation}
Eq.(\ref{EqCurvaturePre}) and Eq.(\ref{EqChargePre}) mean that the coupling constant $g$ may have a geometric meaning, which originates from $G^{RH}$.

This implies that only when affine connection is adopted to represent gauge field can some physical properties be better interpreted. On the other hand, in the general relativity theory, gravitational field is also described by affine connection, so it is convenient to describe gravitational field and gauge field uniformly by affine connection. Therefore, it is necessary to study the affine connection representation of gauge fields. This is the basic motivation of this paper.

There are the following two ways to unify gravitational field and gauge field.

One way is to represent gravitational field as principal bundle connection. We can take the transformation group $Gravi(3,1)$ of gravitational field as the structure group of principal bundle to establish a gauge theory of gravitational field, the local transformation group of which is in the form of $Gravi(3,1) \otimes Gauge(n)$, e.g. Poincar{\'{e}} gauge theory \cite{Kibble1961,Sciama1962,Harnad1976,Heyde1976,Cho1976,Hehl1976_0,Macdowell1977,Chkareuli2017,Chkareuli2017_2,Obukhov2018,ObukhovHehl2020} and metric-affine gauge theory \cite{Hehl1976_1,Hehl1976_2,Hehl1976_3,Hehl1976_4,Hehl1977,Lord1978,Hehl1995,Sobreiro2010,Sobreiro2011,CiriloLombardo2010,CiriloLombardo2011_1,CiriloLombardo2011_2}. This way can be interpreted intuitively as
$$\fbox{\shortstack[c]{\mbox{gravitation} \ \mbox{theory}}}
 \xrightarrow{\mbox{\ \ be incorporated into\ \ }}
\fbox{\shortstack[c]{\mbox{the framework of} \ \mbox{gauge theory.}}}$$

The other way is to represent gauge field as affine connection. This is the approach adopted by this paper. Gravitational field and gauge field can both be described by affine connection. Besides, we will also establish an affine connection representation of elementary particles. This way can be interpreted intuitively as
$$\fbox{\shortstack[c]{\mbox{gauge} \ \mbox{theory}}}
 \xrightarrow{\mbox{\ \ be incorporated into\ \ }}
\fbox{\shortstack[c]{\mbox{the framework of} \ \mbox{gravitation theory.}}}
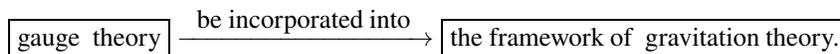
$$

\SubSectionWithLabel{Ideas and methods}

We divide the problem of establishing affine connection representation of gauge fields into three parts as follows. 

(I) Which affine connection is suitable for describing not only gravitational field, but also gauge field and elementary particle field?

(II) How to describe the evolution of these fields in affine connection representation?

(III) What are the concrete forms of electromagnetic, weak, and strong interaction fields in affine connection representation?

{\bf For the problem (I)}. On a Riemannian manifold $(M,G)$, the metric tensor can be expressed as $G_{MN}=\delta_{AB}B^A_M B^B_N$ and $G^{MN}=\delta^{AB}C^M_A C^N_B$, where $B^A_M$ and $C^M_A$ are semi-metrics, or to say frame fields. It is evident that semi-metric is more fundamental than metric, so we hope $B^A_M$ or $C^M_A$ is regarded as a unified frame field of gravitational field and gauge field, and the frame transformation of $B^A_M$ or $C^M_A$ is regarded as gauge transformation. Hence, we need a more general manifold $(M,B^A_M)$ rather than the Riemannian manifold $(M,G)$. 

Next, we put metric and semi-metric together to construct a new connection, which is not only an affine connection, but also a connection on a fibre bundle. In this way, gravitational field and various gauge fields can be unified on a manifold $(M,B^A_M)$ that is defined by semi-metric. 

In addition, we notice that in the theories based on principal bundle connection representation: 

(1) Several complex-valued functions which satisfy the Dirac equation, are sometimes used to refer to a charged lepton field $l$, and sometimes a neutrino field $\nu$. It is not clear how to distinguish these field functions $l$ and $\nu$ by inherent geometric constructions. 

(2) Gauge potentials are abstract; they have no inherent geometric constructions. In other words, the Levi-Civita connection $\Gamma^\mu_{\nu\rho}$ of gravity is constructed by the metric $g_{\mu\nu}$, however it is not explicit what geometric quantity the connection $A^a_\mu$ of gauge field is constructed by.

By contrast, in the affine connection representation of this paper, we are able to use the semi-metrics $B^A_M$ and $C_A^M$ of internal coordinate space to endow particle fields $l$ and $\nu$ and gauge fields $A^a_\mu$ with geometric constructions. Thus, they are not only irreducible representations of group, but also possessed of concrete geometric entities. 

{\bf For the problem (II)}. There is a fundamental difficulty that time is effected by gravitational field, but not effected by gauge field. This leads to an essential difference between the description of evolution of gravitational field and that of gauge field. In this case, it seems difficult to obtain a unified theory of evolution in affine connection representation. Nevertheless, we find that, we can define time as the total metric with respect to all dimensions of internal coordinate space and external coordinate space, and define evolution as one-parameter group of diffeomorphism, to overcome the above difficulty. 

Now that gauge field and gravitational field are both represented as affine connection, then the properties that are related to gauge field, such as charge, current, mass, energy, momentum, and action, must have corresponding affine representations. Thus, Yang-Mills equation, energy-momentum equation, and Dirac equation are turned into geometric properties in gradient direction, in other words, on-shell evolution is characterized by gradient direction. Correspondingly, quantum theory can be interpreted as a geometric theory of distribution of gradient directions. 

{\bf For the problem (III)}. The basic idea is that on a $\mathfrak{D}$-dimensional manifold, the components $B_m^a$ and $C_a^m$ of semi-metrics $B_M^A$ and $C_A^M$ with $m,a\in \{4,5,\cdots,\mathfrak{D}\}$ are regarded as the frame field of electromagnetic, weak, and strong interactions. The other components of $B_M^A$ and $C_A^M$ are regarded as the frame field of gravitation. 

We take the affine connection as 
\begin{equation}\label{AffineConnectionDefinitionA}\begin{aligned}
\Gamma^M_{NP} \DenoteBy \frac{1}{2}\left(\left[^M_{NP}\right] + \left\{^M_{NP}\right\}\right) 
&= \frac{1}{2}\left[C_A^M \left(D_P B_N^A\right) + \left\{^M_{NP}\right\}\right] 
= \frac{1}{2}\left[C_A^M \left(D_C B_N^A\right)b^C_P + \left\{^M_{NP}\right\}\right] \\
&= \frac{1}{2}\left[C_A^M \left(\frac{\partial B_N^A}{\partial \zeta^C} + \left(^A_{BC}\right) B_N^B \right)b^C_P
+ \frac{1}{2}G^{MQ}\left(\frac{\partial G_{NQ}}{\partial x^P}+\frac{\partial G_{PQ}}{\partial x^N}-\frac{\partial G_{NP}}{\partial x^Q}\right)\right] \\
&= \frac{1}{2}\left[\left(C_A^M \frac{\partial B_N^A}{\partial x^P} + C_A^M \left(^A_{BP}\right) B_N^B \right)
+ \frac{1}{2}G^{MQ}\left(\frac{\partial G_{NQ}}{\partial x^P}+\frac{\partial G_{PQ}}{\partial x^N}-\frac{\partial G_{NP}}{\partial x^Q}\right)\right], \\
\end{aligned}\end{equation}
where $b^C_P \DenoteBy \frac{\partial \zeta^C}{\partial x^P}$ is a local coordinate transformation, $\left\{^M_{NP}\right\}$ is Christoffel symbol, $G_{MN}=\delta_{AB}B^A_M B^B_N$, 
$$\left[^M_{NP}\right] \DenoteBy C_A^M \left(D_P B_N^A\right) = C_A^M \frac{\partial B_N^A}{\partial x^P} + C_A^M \left(^A_{BP}\right) B_N^B $$
is said to be a gauge connection, and $\Gamma^M_{NP}$ is said to be a holonomic connection. $\left(^A_{BP}\right) \DenoteBy \left(^A_{BC}\right)b^C_P$.  
$$\left(^A_{BC}\right) \DenoteBy \frac{1}{2} C_{A'}^A \left(\frac{\partial B_B^{A'}}{\partial \zeta^C} + \frac{\partial B_C^{A'}}{\partial \zeta^B}\right)$$
is said to be a torsion-free simple connection. Thus, 

\newpage

$$\Gamma_{MNP} = \frac{1}{2}\left(\left[MNP\right] + \left\{MNP\right\}\right) = \frac{1}{2}\left[\delta_{AD}B_M^D \left(\frac{\partial B_N^A}{\partial x^P} + \left(^A_{BP}\right) B^B_N \right)
+ \frac{1}{2}\left(\frac{\partial G_{NM}}{\partial x^P}+\frac{\partial G_{PM}}{\partial x^N}-\frac{\partial G_{NP}}{\partial x^M}\right)\right].$$
For the sake of simplicity, we firstly consider the affine connection representation of gauge fields without gravitation. That is to say, let 
$$s,i,j=1,2,3; \ \ \ \ \ \ a,m,n,l,q=4,5,\cdots,\mathfrak{D}; \ \ \ \ \ \ A,B,M,N,P=1,2,\cdots,\mathfrak{D};$$
and consider a $\mathfrak{D}$-dimensional manifold $\left(M,B^A_M\right)$ that satisfies the following conditions: 

(i) $B^s_i=\delta^s_i, \ \ B^a_i=0, \ \ B^s_m=0$;

(ii) $G_{ij}=\delta_{ij}, \ \ G_{mn}=const, \ \ G_{mi}=0$;

(iii) When $m\neq n$, $G_{mn}=0$.

\noindent Thus, $\left\{MNP\right\}=0$, $\left[MNP\right]\neq 0$ in general. The components $\Gamma_{mnP}$ of $\Gamma_{MNP}=\frac{1}{2}\left[MNP\right]$ with $m,n\in \{4,5,\cdots,\mathfrak{D}\}$ describe gauge potentials of electromagnetic, weak, and strong interactions. We also use the affine connection $\Gamma^M_{NP}$ to construct elementary particle fields $\rho_{MN}$. The components $\rho_{mn}$ of $\rho_{MN}$ with $m,n\in \{4,5,\cdots,\mathfrak{D}\}$ describe field functions of leptons and quarks. 

The components $G^{mn}$ of $G^{MN}$ with $m,n\in \{4,5,\cdots,\mathfrak{D}\}$ describe coupling constants of particle fields $\rho_{mn}$ and gauge potentials $\Gamma_{mnP}$. The other components of $G^{MN}$ are the metrics of gravitational field. The other components of $\rho_{MN}$ and $\Gamma_{MNP}$ provide possible candidates for dark matters and their interactions.

\SubSectionWithLabel{Content and organization}

In this paper, we are going to show how to construct the affine connection representation of gauge fields. Sections are organized as follows.

Corresponding to the problem (I), in section \ref{Mathematical preparations}, we make some necessary mathematical preparations, and discuss the coordinate transformation and frame transformation of the above connection. Meanwhile, in order to make the languages that are used to describe gauge field and gravitational field unified and harmonized, we generalize the notion of reference system and give it a strict mathematical definition. The reference system in conventional sense is just only defined on a local coordinate neighborhood, and it has only $(1+3)$ dimensions. But in this paper, we define the concept of reference-system over the entire manifold. It is possessed of more dimensions but diffenrent from Kaluza-Klein theory \cite{Kaluza1921,Klein1926_1,Klein1926_2} and string theories \cite{Goddard1973,Polchinski1996arXiv,Taylor1998arXiv,Bachas1999arXiv,Giveon1998,Wess1974,Hanany1997,Douglas2004,Haber1985,Gunion1986_1,Gunion1986_2,Gunion1988,Gunion1993arXiv}. Thus, both of gravitational field and gauge field are regarded as special cases of such a concept of reference-system.

Corresponding to the problem (II), in section \ref{The evolution in affine connection representation of gauge fields}, we establish the general theory of evolution in affine connection representation of gauge fields, and in section \ref{Affine connection representation of gauge fields in classical spacetime}, we discuss the application of this general theory of evolution to $(1+3)$-dimensional classical spacetime. 

Corresponding to the problem (III), in sections \ref{Affine connection representation of the gauge field of weak-electromagnetic interaction} to \ref{Affine connection representation of the unified gauge field}, we show concrete forms of affine connection representations of electromagnetic, weak, and strong interaction fields.

Some important topics are organized as follows.

(1) Time is regarded as the total metric with respect to all spatial dimensions including external coordinate space and internal coordinate space, see Definition \ref{DefinitionOfTimeMetric} and Remark \ref{DiscussionOfMetricOfClassicalSpacetime} for detail. The CPT inversion is interpreted as the composition of full inversion of coordinates and full inversion of metrics, see section \ref{Inversion transformation in affine connection representation} for detail. The conventional $(1+3)$-dimensional Minkowski coordinate $x^\mu$ originates from the general $\mathfrak{D}$-dimensional coordinate $x^M$. The construction method of extra dimensions is different from those of Kaluza-Klein theory and string theory, see section \ref{Classical spacetime reference-system} for detail.

(2) On-shell evolution is characterized by gradient direction field. See sections \ref{On-shell evolution as a gradient},  \ref{On-shell evolution of potential field and affine connection representation of Yang-Mills equation},  \ref{On-shell evolution of general charge and affine connection representation of mass, energy, momentum, and action} and \ref{Affine connection representation of classical spacetime evolution} for detail. Quantum theory is regarded as a geometric theory of distribution of gradient directions. We show two dual descriptions of gradient direction. They just exactly correspond to the Schr{\"{o}}dinger picture and the Heisenberg picture. In these points of view, the gravitational theory and quantum theory become coordinated. They have a unified description of evolution, and the definition of Feynman propagator is simplified to a stricter form. See sections \ref{Two dual descriptions of gradient direction field} and  \ref{Quantum evolution as a distribution of gradient directions} for detail.

(3) Yang-Mills equation originates from a geometric property of gradient direction. We show the affine connection representation of Yang-Mills equation. See section \ref{On-shell evolution of potential field and affine connection representation of Yang-Mills equation} and section \ref{From classical spacetime back to full-dimensional space} for detail.

(4) Energy-momentum equation originates from a geometric property of gradient direction. We show the affine connection representation of mass, energy, momentum and action, see section \ref{On-shell evolution of general charge and affine connection representation of mass, energy, momentum, and action}, Definition \ref{DefinitionOfMinkowskiEnergyMomentum} and Discussion \ref{DiscussionOfMinkowskiEnergyMomentumEquation} for detail. Furthermore, we also show the affine connection representation of Dirac equation, see section \ref{Affine connection representation of Dirac equation} for detail.

(5) Why do not neutrinos participate in the electromagnetic interactions? And why do not right-handed neutrinos participate in the weak interactions with $W$ bosons? In the theory of this paper, they are natural and geometric results of affine connection representation of gauge fields; therefore, they are not necessary to be regarded as postulates anymore. See Proposition \ref{PropositionOfLeptonEvolutionForm} and Proposition \ref{PropositionOfMNSMixing} for detail.

(6) In section \ref{Affine connection representation of the unified gauge field}, we give new interpretations to PMNS mixing of leptons, CKM mixing of quarks, and color confinement. That is to say, in affine connection representation of gauge fields, these physical properties can be interpreted as geometric properties on manifold.

\SectionWithLabel{Mathematical preparations}

\SubSectionWithLabel{Geometric manifold}

In order to make the languages that are used to describe gauge field and gravitational field unified and harmonized, we adopt the following definition. 

\sDefinition\label{DefinitionOfGeneralReferenceSystem} Let $M$ be a $\mathfrak{D}$-dimensional connected smooth real manifold. $\forall p\in M$, take a coordinate chart $(U_p,\varphi_{Up})$ on a neighborhood $U_p$ of $p$. They constitute a coordinate covering 
$$\varphi \DenoteBy \{(U_p,\varphi_{Up})\}_{p\in M},$$ 
which is said to be a {\bf point-by-point covering}. For the sake of simplicity, $U_p$ can be denoted by $U$, and $\varphi_{Up}$ by $\varphi_U$.

Let $\varphi$ and $\psi$ be two point-by-point coverings. For the two coordinate frames $\varphi_U$ and $\psi_U$ on the neighborhood $U$ of point $p$, if
$$f_p\DenoteBy\varphi_U\circ\psi_U^{-1}: \psi_U(U)\to\varphi_U(U),\ \ \xi^A \mapsto x^M$$
is a smooth homeomorphism, $f_p$ is called a {\bf local reference-system}. 

If every $p\in M$ is endowed with a local reference-system $f(p)$, and we require the semi-metrics $B_M^A$ and $C_A^M$ in Eq.(\ref{DefinitionFormulaSemiMetrics}) to be smooth real functions on $M$, then 
\begin{equation}f:M\to REF,\ \ p\mapsto f(p)\end{equation}
is said to be a {\bf reference-system} on $M$, and $(M,f)$ is said to be a {\bf geometric manifold}.

\SubSectionWithLabel{Metric and semi-metric}

In the absence of a special declaration, the indices take values as $A,B,C,D,E=1,2,\cdots ,\mathfrak{D}$ and $M,N,P,Q,R=1,2,\cdots ,\mathfrak{D}$. The derivative functions 
\begin{equation}
b_M^A\DenoteBy\frac{{\partial\xi^A}}{{\partial x^M}},\ \ \ \ \ \ \ c_A^M\DenoteBy\frac{{\partial x^M}}{{\partial\xi^A}}
\end{equation}
of $f(p)$ on $U_p$ define the semi-metrics (or to say frame field) $B_M^A$ and $C_A^M$ of $f$ on the manifold $M$, that are
\begin{equation}\label{DefinitionFormulaSemiMetrics}
B_M^A :M\to\mathbb{R},\ \ p\mapsto \ B_M^A(p)\DenoteBy (b_{f(p)})_M^A(p),\ \ \ \ \ \ \ \ C_A^M :M\to\mathbb{R},\ \ p\mapsto \ C_A^M(p)\DenoteBy (c_{f(p)})_A^M(p).
\end{equation}
Let $\delta_{AB}=\delta^{AB}=\delta_B^A = Kronecker(A,B)$ and $\varepsilon_{MN}=\varepsilon^{MN}=\varepsilon_N^M = Kronecker(M,N)$. The metric tensors of $f$ are 
\begin{equation}G_{MN}=\delta_{AB}B_M^A B_N^B,\ \ \ \ \ \ \ \ H_{AB}=\varepsilon_{MN}C_A^M C_B^N.\end{equation}
Similarly, it can also be defined that ${\bar b}_A^M \DenoteBy\frac{\partial\xi_A}{\partial x_M}, \ \ {\bar c}_M^A \DenoteBy\frac{\partial x_M}{\partial\xi_A}$ and corresponding ${\bar B}_A^M, \ \ {\bar C}_M^A$.

\SubSectionWithLabel{Gauge transformation in affine connection representation}

$\forall p\in M$, $f(p)\DenoteBy\rho_U\circ\psi_U^{-1}$ induces local reference-system transformations
$$\begin{aligned}
L_{f(p)}: \ \ &k(p)\DenoteBy\psi_U\circ\varphi_U^{-1}\ \ \mapsto \ \ \rho_U\circ\varphi_U^{-1} = f(p)\circ k(p),\hfill\\
R_{f(p)}: \ \ &h(p)\DenoteBy\varphi_U\circ\rho_U^{-1}\ \ \mapsto \ \ \varphi_U\circ\psi_U^{-1} = h(p)\circ f(p),\hfill\\
\end{aligned}$$
and reference-system transformations on the manifold $M$
\begin{equation}L_f:p\mapsto L_{f(p)}, \ \ \ \ R_f:p\mapsto R_{f(p)}.\end{equation}

\newpage

\noindent We also speak of $L_f$ and $R_f$ as {\bf (affine) gauge transformations}. 

(i) $L_f$ and $R_f$ are identical transformations if and only if $[B_M^A]$ of $f$ is an identity matrix.

(ii) $L_f$ and $R_f$ are flat transformations if and only if $\forall p_1,p_2\in M$, $B_M^A(p_1) = B_M^A(p_2)$.

(iii) $L_f$ and $R_f$ are orthogonal transformations if and only if $\delta_{AB}B_M^A B_N^B=\varepsilon_{MN}$.

\noindent The totality of all reference-system transformations on $M$ is denoted by $GL(M)$, which is a subgroup of $\displaystyle \bigotimes_{p\in M} GL(\mathfrak{D},\mathbb{R})_p$, where $\bigotimes$ represents external direct product.

\SubSectionWithLabel{Coordinate transformation of holonomic connection and frame transformation of gauge connection}

Suppose there are reference-systems $g$ and $\mathfrak{g}$ on the manifold $M$, denote $\mathcal{G}\DenoteBy \mathfrak{g} \circ g$, and $\forall p\in M$, on the neighborhood $U$ of $p$, $g(p)$ and $\mathfrak{g}(p)$ satisfy
$$(U,x^M)\xleftarrow{g(p)}(U,\zeta^A)\xleftarrow{\mathfrak{g}(p)}(U,\beta^{A'}).$$
On the geometric manifold $(M,\mathfrak{g})$ we define {\bf torsion-free simple connection} $D$ and its coefficients $(^A_{BC})_\mathfrak{g}$ by
\begin{equation} \label{DefinitionTorsionFreeSimpleConnection}
D \frac{\partial}{\partial \zeta^B} \DenoteBy (\omega_\mathfrak{g})^A_B \otimes \frac{\partial}{\partial \zeta^A} = (^A_{BC})_\mathfrak{g} d\zeta^C \otimes \frac{\partial}{\partial \zeta^A}
= \frac{1}{2} (C_\mathfrak{g})_{A'}^A \left(\frac{\partial (B_\mathfrak{g})_B^{A'}}{\partial \zeta^C} + \frac{\partial (B_\mathfrak{g})_C^{A'}}{\partial \zeta^B}\right) d\zeta^C \otimes \frac{\partial}{\partial \zeta^A}.
\end{equation}
Then, we can compute the absolute derivative of the frame field $\frac{\partial}{\partial x^N}$
$$\begin{aligned}
D\frac{\partial}{\partial x^N} &= D\left((B_g)_N^B \frac{\partial}{\partial \zeta^B}\right) 
= d(B_g)_N^B \otimes \frac{\partial}{\partial \zeta^B} + (B_g)_N^B D \frac{\partial}{\partial \zeta^B} \\
&= \frac{\partial (B_g)_N^B}{\partial \zeta^C} d\zeta^C \otimes \frac{\partial}{\partial \zeta^B} + (B_g)_N^B (^A_{BC})_\mathfrak{g} d\zeta^C \otimes \frac{\partial}{\partial \zeta^A} 
= \left( \frac{\partial (B_g)_N^A}{\partial \zeta^C} + (B_g)_N^B (^A_{BC})_\mathfrak{g} \right) d\zeta^C \otimes \frac{\partial}{\partial \zeta^A}, \\
\end{aligned}$$
Thus, it is obtained that 
$$D_C (B_g)_N^A = \frac{\partial (B_g)_N^A}{\partial \zeta^C} + (B_g)_N^B (^A_{BC})_\mathfrak{g} \ .$$
Denote $D_P \DenoteBy (b_{g(p)})^C_P D_C$; thus we can define on $(M,\mathcal{G})$ the required {\bf gauge connection}, which is 
\begin{equation} \label{DefinitionGaugeConnection}
\left[^M_{NP}\right]_\mathcal{G} \DenoteBy (C_g)_A^M D_P (B_g)_N^A = (C_g)_A^M \frac{\partial (B_g)_N^A}{\partial x^P} + (C_g)_A^M \left(^A_{BP}\right)_\mathfrak{g} (B_g)_N^B \ .
\end{equation}
It is important that $\left[^M_{NP}\right]_\mathcal{G}$ is not only an affine connection on $(M,\mathcal{G})$, but also a connection on frame bundle. 

\noindent {\bf (I). $\left[^M_{NP}\right]_\mathcal{G}$ as an affine connection.} Under the coordinate transformation
$\displaystyle L_{k(p)}: (U,x^M) \to (U,x^{M'})$, $b^M_{M'}\DenoteBy\frac{\partial x^M}{\partial x^{M'}}$, $c_M^{M'}\DenoteBy\frac{\partial x^{M'}}{\partial x^M}$,
$(B_g)^A_M \mapsto (B_g)^A_{M'} = b^M_{M'} (B_g)^A_M,\ \ (C_g)^M_A \mapsto (C_g)^{M'}_A = c^{M'}_M (C_g)^M_A$. Consequently, the gauge connection $\left[^M_{NP}\right]_\mathcal{G}$ is transformed according to
\begin{equation} \label{SimpleConnectionCoordinateTransformation}
L_{k(p)}: \ \ \left[^M_{NP}\right]_\mathcal{G} \ \ \mapsto \ \ \left[^{M'}_{N'P'}\right]_\mathcal{G} =c_M^{M'} \left[^M_{NP}\right]_\mathcal{G} b_{N'}^N b_{P'}^P +c_M^{M'}\frac{{\partial b_{N'}^M}}{{\partial x^{P'}}}.
\end{equation}
Due to Eq.(\ref{SimpleConnectionCoordinateTransformation}), under the coordinate transformation, the {\bf holonomic connection}
\begin{equation}\label{DefinitionFormulaOfHolonomicConnection}\begin{aligned}
(\Gamma_\mathcal{G})^M_{NP} &\DenoteBy \frac{1}{2}\left(\left[^M_{NP}\right]_\mathcal{G} + \left\{^M_{NP}\right\}_\mathcal{G}\right) \\
&= \frac{1}{2}\left[\left((C_g)_A^M \frac{\partial (B_g)_N^A}{\partial x^P} + (C_g)_A^M \left(^A_{BP}\right)_\mathfrak{g} (B_g)_N^B \right)
+ \frac{1}{2}(G_\mathcal{G})^{MQ}\left(\frac{\partial (G_\mathcal{G})_{NQ}}{\partial x^P}+\frac{\partial (G_\mathcal{G})_{PQ}}{\partial x^N}-\frac{\partial (G_\mathcal{G})_{NP}}{\partial x^Q}\right)\right] \\
\end{aligned}\end{equation}
is transformed according to
\begin{equation} \label{HolonomicConnectionCoordinateTransformation}
L_{k(p)}: \ \ (\Gamma_\mathcal{G})_{NP}^M \ \ \mapsto \ \ (\Gamma_\mathcal{G})_{N'P'}^{M'} =c_M^{M'}(\Gamma_\mathcal{G})_{NP}^M b_{N'}^N b_{P'}^P +c_M^{M'}\frac{{\partial b_{N'}^M}}{{\partial x^{P'}}}.
\end{equation}

\noindent {\bf (II). $\left[^M_{NP}\right]_\mathcal{G}$ as a connection on frame bundle.} Under the frame transformation 
$\displaystyle L_k: (M,\mathcal{G}) \mapsto (M,\mathcal{G}'), \ \ \frac{\partial}{\partial x^M} \mapsto \frac{\partial}{\partial x^{M'}} = (B_k)^M_{M'} \frac{\partial}{\partial x^M} \ , $
$(B_g)^A_M \mapsto (B_{g'})^A_{M'} = (B_k)^M_{M'} (B_g)^A_M,\ \ (C_g)^M_A \mapsto (C_{g'})^{M'}_A = (C_k)^{M'}_M (C_g)^M_A$. Consequently, the gauge connection $\left[^M_{NP}\right]_\mathcal{G}$ is tranformed according to
\begin{equation} \label{GaugeConnectionFrameTransformation}
L_k: \ \ \left[^M_{NP}\right]_\mathcal{G} \ \ \mapsto \ \ \left[^{M'}_{N'P'}\right]_{\mathcal{G}'} = \left[^{M'}_{N'P}\right]_{\mathcal{G}'} b^P_{P'} = \left( (C_k)_M^{M'} \left[^M_{NP}\right]_\mathcal{G} (B_k)_{N'}^N + (C_k)_M^{M'}\frac{{\partial (B_k)_{N'}^M}}{{\partial x^{P}}} \right) b^P_{P'}
\end{equation}
Eq.(\ref{SimpleConnectionCoordinateTransformation}) and Eq.(\ref{GaugeConnectionFrameTransformation}) show that $\left[^M_{NP}\right]_\mathcal{G}$ is not only an affine connection, but also a connection on frame bundle. 

\newpage

Apply Eq.(\ref{SimpleConnectionCoordinateTransformation})$\sim$(\ref{GaugeConnectionFrameTransformation}) to the curvature tensors
$$\begin{aligned}
&\left[^M_{NPQ}\right] \DenoteBy\frac{\partial\left[^M_{NQ}\right]}{\partial x^P}-\frac{\partial\left[^M_{NP}\right]}{\partial x^Q}+\left[^M_{HP}\right]\left[^H_{NQ}\right] - \left[^H_{NP}\right]\left[^M_{HQ}\right],\\
&\left\{^M_{NPQ}\right\} \DenoteBy\frac{\partial\left\{^M_{NQ}\right\}}{\partial x^P}-\frac{\partial\left\{^M_{NP}\right\}}{\partial x^Q}+\left\{^M_{HP}\right\}\left\{^H_{NQ}\right\} - \left\{^H_{NP}\right\}\left[^M_{HQ}\right\},\\
&R^M_{NPQ} \DenoteBy\frac{\partial\Gamma^M_{NQ}}{\partial x^P}-\frac{\partial\Gamma^M_{NP}}{\partial x^Q}+\Gamma^M_{HP}\Gamma^H_{NQ} - \Gamma^H_{NP}\Gamma^M_{HQ},\\
\end{aligned}$$
then it is obtained that
\begin{equation} \label{HolonomicCurvatureFrameTransformation}
\begin{aligned}
&L_k: \ \ \left[^M_{NPQ}\right]_\mathcal{G} \ \ \mapsto \ \ \left[^{M'}_{N'P'Q'}\right]_{\mathcal{G}'}=\left[^{M'}_{N'PQ}\right]_{\mathcal{G}'} b^P_{P'} b^Q_{Q'} = \left((C_k)_M^{M'} \left[^M_{NPQ}\right]_\mathcal{G} (B_k)_{N'}^N\right) b^P_{P'} b^Q_{Q'} \ , \\
&L_{k(p)}: \ \ \left[^M_{NPQ}\right]_\mathcal{G} \ \ \mapsto \ \ \left[^{M'}_{N'P'Q'}\right]_\mathcal{G} = c_M^{M'} \left[^M_{NPQ}\right]_\mathcal{G} b_{N'}^N b^P_{P'} b^Q_{Q'} \ , \\
&L_{k(p)}: \ \ \left\{^M_{NPQ}\right\}_\mathcal{G} \ \ \mapsto \ \ \left\{^{M'}_{N'P'Q'}\right\}_\mathcal{G} = c_M^{M'} \left\{^M_{NPQ}\right\}_\mathcal{G} b_{N'}^N b^P_{P'} b^Q_{Q'} \ , \\
&L_{k(p)}: \ \ (R_\mathcal{G})^M_{NPQ} \ \ \mapsto \ \ (R_\mathcal{G})^{M'}_{N'P'Q'} = c_M^{M'} (R_\mathcal{G})^M_{NPQ} b_{N'}^N b^P_{P'} b^Q_{Q'} \ . \\
\end{aligned}
\end{equation}
We see from Eq.(\ref{HolonomicCurvatureFrameTransformation}) that the $\left[^M_{NPQ}\right]_\mathcal{G}$ without gravitation is both a curvature tensor of affine connection and a curvature tensor on frame bundle, and that the $(R_\mathcal{G})^M_{NPQ}$ with gravitation is a curvature tensor of affine connection, but not a curvature tensor on frame bundle. In other words, under the gauge transformation $L_k$, $\left[^M_{NPQ}\right]_\mathcal{G}$ and $\left[^{M'}_{N'PQ}\right]_{\mathcal{G}'}$ represent the same physical state, while $(R_\mathcal{G})^M_{NPQ}$ and $(R_{\mathcal{G}'})^{M'}_{N'PQ}$ represent different physical states. This shows that the gravitational field in $(R_\mathcal{G})^M_{NPQ}$ makes the gauge frames $B^A_M$ and $C^M_A$ have physical effects.

\SectionWithLabel{The evolution in affine connection representation of gauge fields}

Now that we have the required affine connection, next we have to solve the problem that how to describe the evolution in affine connection representation. 

In the existing theories, time is effected by gravitational field, but not effected by gauge field. This leads to an essential difference between the description of evolution of gravitational field and that of gauge field. In this case, it is difficult to obtain a unified theory of evolution in affine connection representation. We adopt the following way to overcome this difficulty. 

\SubSectionWithLabel{The relation between time and space}

\sDefinition\label{DefinitionOfTimeMetric}\label{DefinitionOfInternalSpaceAndExternalSpace} Suppose $M=P\times N$ and $r\DenoteBy\dim P = 3$. Let 
$$A,B,M,N=1,\cdots,\mathfrak{D}\ ;\ \ \ \ \ \ s,i=1,\cdots,r\ ;\ \ \ \ \ \ a,m=r+1,\cdots,\mathfrak{D}.$$
On a geometric manifold $(M,f)$, the $d\xi^0$ and $dx^0$ which are defined by
\begin{equation}
\LeftList
 (d\xi^0)^2 \DenoteBy \Summation_{A=1}^\mathfrak{D}{(d\xi^A)^2} = \delta_{AB}d\xi^A d\xi^B =G_{MN}dx^M dx^N,  \hfill\\
 (dx^0)^2 \DenoteBy \Summation_{M=1}^\mathfrak{D}{(dx^M)^2} = \varepsilon_{MN}dx^M dx^N =H_{AB}d\xi^A d\xi^B  \hfill\\ 
\RightList
\end{equation}
are said to be {\bf total space metrics} or {\bf time metrics}. We also suppose 
$$\LeftList
 (d\xi^{(P)})^2 \DenoteBy\Summation_{s=1}^r{(d\xi^s)^2},\ \ \ \ (d\xi^{(N)})^2 \DenoteBy\Summation_{a=r + 1}^\mathfrak{D}{(d\xi^a)^2},\hfill\\
 (dx^{(P)})^2 \DenoteBy\Summation_{i=1}^r{(dx^i)^2},\ \ \ \ (dx^{(N)})^2 \DenoteBy\Summation_{m=r + 1}^\mathfrak{D}{(dx^m)^2}.\hfill\\
\RightList$$
$d\xi^{(N)}$ and $dx^{(N)}$ are regarded as proper-time metrics. For convenience, $P$ is said to be {\bf external space} and $N$ is said to be {\bf internal space}. 

\newpage

\sRemark The above definition implies a new viewpoint about time and space. The relation between time and space in this way is different from the Minkowski coordinates $x^\mu\ (\mu=0,1,2,3)$. Time and space are not the components on an equal footing anymore, but have a relation of total to component. It can be seen later that time reflects the total evolution in the full-dimensional space, while a specific spatial dimension reflects just a partial evolution in a specific direction.

\SubSectionWithLabel{Evolution path as a submanifold}

\sDefinition\label{DefinitionOfMotionAndInteraction}\label{DefinitionOfEvolutionPath} Let there be reference-systems $f$, $g$, $\mathfrak{f}$, $\mathfrak{g}$ on a manifold $M$, such that $\forall p\in M$, on the neighborhood $U$ of $p$, 
\begin{equation}\label{DefinitionFormulaOfEvolution}(U,\alpha^{A'})\xrightarrow{\mathfrak{f}(p)}(U,\xi^A)\xrightarrow{f(p)}(U,x^M)\xleftarrow{g(p)}(U,\zeta^A)\xleftarrow{\mathfrak{g}(p)}(U,\beta^{A'}).\end{equation}
Denote $\mathcal{F}\DenoteBy \mathfrak{f} \circ f$ and $\mathcal{G}\DenoteBy \mathfrak{g} \circ g$, then we say $\mathcal{F}$ and $\mathcal{G}$ {\bf motion relatively and interact mutually}, and also we say that $\mathcal{F}$ evolves in $\mathcal{G}$, or $\mathcal{F}$ evolves on the geometric manifold $(M,\mathcal{G})$. Meanwhile, $\mathcal{G}$ evolves in $\mathcal{F}$, or we say $\mathcal{G}$ evolves on $(M,\mathcal{F})$.

From Eq.(\ref{DefinitionGaugeConnection}) we know that in $\mathcal{F}$ and $\mathcal{G}$, gauge fields originate from $\mathfrak{f}$ and $\mathfrak{g}$, and gravitational fields $(G_\mathcal{F})_{MN}$ and $(G_\mathcal{G})_{MN}$ are effected by $\mathfrak{f}$ and $\mathfrak{g}$, respectively. We are going to describe their evolutions step by step in the following sections.

Let there be a one-parameter group of diffeomorphisms 
$$\varphi_X :M\times\mathbb{R}\to M$$ 
acting on $M$, such that $\varphi_X (p,0)=p$. Thus, $\varphi_X$ determines a smooth tangent vector field $X$ on $M$. If $X$ is nonzero everywhere, we say $\varphi_X$ is {\bf a set of evolution paths}, and $X$ is an {\bf evolution direction field}. Let $T\subseteq\mathbb{R}$ be an interval; then, the regular imbedding
\begin{equation}L_p \DenoteBy \varphi_{X,p} :T\to M,\ t \mapsto \varphi_X(p,t)\end{equation}
is said to be an {\bf evolution path} through $p$. The tangent vector $\frac{d}{dt}\DenoteBy [L_p]=X(p)$ is called an {\bf evolution direction} at $p$. For the sake of simplicity, we also denote $L_p\DenoteBy L_p(T) \subset M$; then, 
\begin{equation}\pi: L_p \to M,\ q \mapsto q\end{equation}
is also a regular imbedding. If it is not necessary to emphasize the point $p$, $L_p$ is denoted by $L$ concisely.

In order to describe physical evolution, next we are going to strictly describe the mathematical properties of the reference-systems $f$ and $g$ which are sent onto the evolution path $L$. 

\sDefinition\label{DefinitionOfPathMapOfCoordinate}\label{DefinitionOfEvolutionOfSlackTight} Let the time metrics of $(U,\xi^A)$, $(U,x^M)$, and $(U,\zeta^A)$ be $d\xi^0$, $dx^0$, and $d\zeta^0$, respectively. On $U_L \DenoteBy U\Jiao L_p$ we have parameter equations
\begin{equation}\label{CoordinateEvolution}
\LeftList
    \xi^A =\xi^A(x^0),\ \ \ \ x^M = x^M(\xi^0),\ \ \ \ \zeta^A =\zeta^A(x^0),\hfill\\
    \xi^0 =\xi^0(x^0),\ \ \ \ \ \ x^0 =x^0(\xi^0),\ \ \ \ \ \ \zeta^0 =\zeta^0(x^0).\hfill\\  
\RightList
\end{equation}
Take $f$ for example, according to Eq.(\ref{CoordinateEvolution}), on $U_L$ we define
$$\begin{aligned}
&b_0^A \DenoteBy\frac{d\xi^A}{dx^0}, \ \ \ \ \ \ \ b_0^0 \DenoteBy\frac{d\xi^0}{dx^0}, \ \ \ \ \ \ 
\varepsilon_0^M \DenoteBy \frac{dx^M}{dx^0}=b_0^0 c_0^M =b_0^A c_A^M, \\
&c_0^M \DenoteBy\frac{dx^M}{d\xi^0}, \ \ \ \ \ c_0^0 \DenoteBy\frac{dx^0}{d\xi^0}, \ \ \ \ \ \ \delta_0^A \DenoteBy \frac{d\xi^A}{d\xi^0}=c_0^0 b_0^A =c_0^M b_M^A. \\
\end{aligned}$$
Define $d\xi_0 \DenoteBy\frac{dx^0}{d\xi^0}dx^0$ and $dx_0 \DenoteBy\frac{d\xi^0}{dx^0}d\xi^0$, which induce $\frac{d}{d\xi_0}$ and $\frac{d}{dx_0}$, such that $\LeftAngle{\frac{d}{d\xi_0},d\xi_0}\RightAngle =1$, $\LeftAngle{\frac{d}{dx_0},dx_0}\RightAngle =1$. So we can also define
$$\begin{aligned}
&\bar b_A^0 \DenoteBy\frac{d\xi_A}{dx_0},\ \ \ \ \ \ \ \bar b_0^0 \DenoteBy\frac{d\xi_0}{dx_0},\ \ \ \ \ \  \bar{\varepsilon}_M^0 \DenoteBy \frac{dx_M}{dx_0}=\bar b_0^0 \bar c_M^0 =\bar b_A^0 \bar c_M^A, \\
&\bar c_M^0 \DenoteBy\frac{dx_M}{d\xi_0},\ \ \ \ \ \bar c_0^0 \DenoteBy\frac{dx_0}{d\xi_0}, \ \ \ \ \ \ \bar{\delta}_A^0 \DenoteBy \frac{d{\xi}_A}{d\bar{\xi}_0}=\bar c_0^0 \bar b_A^0 =\bar c_M^0 \bar b_A^M. \\
\end{aligned}$$
They determine the following smooth functions on the entire $L$, similar to section \ref{Metric and semi-metric}, that
$$\LeftList
  B_0^A :L \to\mathbb{R},\ \ p\mapsto B_0^A(p)\DenoteBy (b_{f(p)})_0^A(p),\ \ \ \ \ \ \ \ C_0^M :L \to\mathbb{R},\ \ p\mapsto C_0^M(p)\DenoteBy (c_{f(p)})_0^M(p),\hfill\\ 
 \bar B_A^0 :L \to\mathbb{R},\ \ p\mapsto\bar B_A^0(p)\DenoteBy(\bar b_{f(p)})_{A}^0(p),\ \ \ \ \ \ \ \ \bar C_M^0 \ :L \to\mathbb{R},\ \ p\mapsto\bar C_M^0(p)\DenoteBy(\bar c_{f(p)})_{M}^0(p),\hfill\\ 
  B_0^0 :\ L \to\mathbb{R},\ \ p\mapsto B_0^0(p)\DenoteBy (b_{f(p)})_0^0(p),\ \ \ \ \ \ \ \ \ C_0^0 \ \ :L \to\mathbb{R},\ \ p\mapsto C_0^0(p)\DenoteBy (c_{f(p)})_0^0(p),\hfill\\ 
 \bar B_0^0 :\ L \to\mathbb{R},\ \ p\mapsto\bar B_0^0(p)\DenoteBy(\bar b_{f(p)})_{0}^0(p),\ \ \ \ \ \ \ \ \ \bar C_0^0 \ \ :L \to\mathbb{R},\ \ p\mapsto\bar C_0^0(p)\DenoteBy(\bar c_{f(p)})_{0}^0(p).\hfill\\ 
\RightList$$
For convenience, we still use the notations $\varepsilon$ and $\delta$ and have the following smooth functions.
$$\begin{aligned}
&\varepsilon_0^M \DenoteBy B_0^0 C_0^M =B_0^A C_A^M, \ \ \ \ \ \ \ \ \delta_0^A \DenoteBy C_0^0 B_0^A =C_0^M B_M^A, \ \ \ \ \ \ \ \ G_{00}\DenoteBy B_0^0 B_0^0=G_{MN}\varepsilon_0^M\varepsilon_0^N. \\
&\bar\varepsilon_M^0 \DenoteBy\bar B_0^0\bar C_M^0 =\bar B_A^0\bar C_M^A,\ \ \ \ \ \ \ \ \ \bar\delta_A^0 \DenoteBy\bar C_0^0\bar B_A^0 =\bar C_M^0\bar B_A^M, \ \ \ \ \ \ \ \ G^{00}\DenoteBy C_0^0 C_0^0=G^{MN}\bar\varepsilon^0_M\bar\varepsilon^0_N. \\
\end{aligned}$$
It is easy to verify that $dx_0=G_{00}dx^0$ and $\frac{d}{dx_0}= G^{00}\frac{d}{dx^0}$ are both true on $L$ by a simple calculation.

\SubSectionWithLabel{Evolution lemma}

We have the following two evolution lemmas. The affine connection representations of Yang-Mills equation, energy-momentum equation, and Dirac equation are dependent on them.

\sDefinition\label{DefinitionOfEquivalenceAndHomomorphism} $\forall p\in L$, suppose $T_p(M)$ and $T_p(L)$ are tangent spaces, $T_p^*(M)$ and $T_p^*(L)$ are cotangent spaces. The regular imbedding $\pi :L\to M,\ q\mapsto q$ induces the tangent map and the cotangent map
\begin{equation}
\LeftList
\pi_* :T_p(L)\to T_p(M),\ \ [\gamma_L]\mapsto [{\pi\circ\gamma_L}], \hfill\\
\pi^* :T_p^*(M)\to T_p^*(L),\ \ df\mapsto d(f \circ \pi). \hfill\\
\RightList
\end{equation}
Evidently, $\pi_*$ is an injection, and $\pi^*$ is a surjection. $\forall\frac{d}{dt_L}\in T_p(L),\ \frac{d}{dt}\in T_p(M),\ df\in T_p^*(M),\ df_L\in T_p^*(L)$, if and only if
\begin{equation}
  \frac{d}{dt} = \pi_*\left(\frac{d}{dt_L}\right),\ \ \ df_L = \pi^*(df)
\end{equation}
are true, we denote
\begin{equation}\frac{d}{dt}\DengJiaYu\frac{d}{dt_L},\ \ \ df\TongTaiYu df_L.\end{equation}
Then, we have the following two propositions that are evidently true.

\sProposition\label{PropositionConjugationofTangentAndCotangent} If $\frac{d}{dt}\DengJiaYu\frac{d}{dt_L}$ and $df\TongTaiYu df_L$, then 
\begin{equation} \label{EvolutionLemmaA}\LeftAngle{\frac{d}{dt},df}\RightAngle =\LeftAngle{\frac{d}{dt_L},df_L}\RightAngle.\end{equation}

\sProposition\label{PropositionOfEvolutionLemma} The following conclusions are true.
\begin{equation} \label{EvolutionLemmaB}
\LeftGathered
\begin{aligned}
  w^M\frac{\partial}{\partial x^M}&\DengJiaYu w^0\frac{d}{dx^0} &\IfAndOnlyIf \ \ w^M =w^0\varepsilon_0^M, \hfill\\
  w_M dx^M &\TongTaiYu w_0 dx^0 &\IfAndOnlyIf \ \ w_M \varepsilon_0^M =w_0, \hfill\\ 
\end{aligned}
\RightGathered\ \ \ \ \ \ 
\LeftGathered
\begin{aligned}
 \bar w_M\frac{\partial}{\partial x_M} &\DengJiaYu\bar w_0\frac{d}{dx_0} &\IfAndOnlyIf \ \ \bar w_M =\bar w_0\bar\varepsilon_M^0, \hfill\\
 \bar w^M dx_M &\TongTaiYu\bar w^0 dx_0 &\IfAndOnlyIf \ \ \bar w^M \bar\varepsilon_M^0 =\bar w^0. \hfill\\ 
\end{aligned}
\RightGathered
\end{equation}

\SubSectionWithLabel{On-shell evolution as a gradient}
 
Let $\mathbf{T}$ be a smooth $n$-order tensor field. The restriction on $(U,x^M)$ is $\mathbf{T}\DenoteBy t \left\{{\frac{\partial}{{\partial x }}\otimes dx }\right\}$, where $\left\{{\frac{\partial}{{\partial x }}\otimes dx }\right\}$ represents the tensor basis generated by several $\frac{\partial }{\partial x^M}$ and $dx^M$, and the tensor coefficients of $\mathbf{T}$ are concisely denoted by $t:U\to \mathbb{R}$. 

Let $D$ be a holonomic connection. Consider $D\mathbf{T}\DenoteBy t_{;Q} dx^Q \otimes\left\{{\frac{\partial}{{\partial x}}\otimes dx}\right\}$. Denote 
$$Dt \DenoteBy t_{;Q}dx^Q,\ \ \ \ \nabla t \DenoteBy t_{;Q}\frac{\partial}{\partial x_Q}.$$
$\forall p\in M$, the integral curve of $\nabla t$, that is, $L_p \DenoteBy \varphi_{\nabla t, p}$ , is a {\bf gradient line} of $\mathbf{T}$. It can be seen later that the above gradient operator $\nabla$ characterizes the {\bf on-shell evolution}.

For any evolution path $L$, let $U_L\DenoteBy U \cap L$. Denote ${t_L}\DenoteBy t|_{U_L}$ and ${t_L}_{;0} \DenoteBy t_{;Q}\varepsilon_0^Q$, as well as
$$D_L {t_L}\DenoteBy {t_L}_{;0}dx^0,\ \ \ \ \nabla_L {t_L}\DenoteBy {t_L}_{;0}\frac{d}{dx_0}.$$

\sProposition\label{TheoremOfActualEvolution} The following conclusions are evidently true. 

(i) $Dt \TongTaiYu D_L t_L$ if and only if $L$ is an arbitrary evolution path. 

(ii) $\nabla t \DengJiaYu \nabla_L t_L$ if and only if $L$ is a gradient line of $\mathbf{T}$.

\sRemark More generally, suppose there is a tensor $\mathbf{U}\DenoteBy u_{Q} dx^Q\otimes\left\{{\frac{\partial}{\partial x}\otimes dx}\right\}$. In such a notation, all the indices are concisely ignored except $Q$. $u_{Q} dx^Q$ uniquely determines a characteristic direction $u_Q \frac{\partial}{\partial x_Q}$.

If the system of 1-order linear partial differential equations $t_{;Q}=u_{Q}$ has a solution $t$, then it is true that $Dt=u_Q dx^Q$ and $\nabla t=u_Q \frac{\partial}{\partial x_Q}$. Thus, in the evolution direction $[L]=u_Q \frac{\partial}{\partial x_Q}$, the following conclusions are true. 
\begin{equation}Dt \TongTaiYu D_L t_L, \ \ \ \ \ \nabla t \DengJiaYu \nabla_L t_L,\end{equation}
where $D_L t_L \DenoteBy u_0 dx^0$, $\nabla_L t_L \DenoteBy u_0 \frac{d}{dx_0}$, and $u_0 \DenoteBy u_Q \varepsilon^Q_0$.

Now for any geometric property in the form of tensor $\mathbf{U}$, we are able to express its on-shell evolution in the form of $\nabla t$.

Next, two important on-shell evolutions are discussed in the following two sections. One is the on-shell evolution of the potential field of a reference-system. The other is the one that a general charge of a reference-system evolves in the potential field of another reference-system.

\SubSectionWithLabel{On-shell evolution of potential field and affine connection representation of Yang-Mills equation}

The table I of the article\cite{WuYang1975} proposes a famous correspondence between gauge field terminologies and fibre bundle terminologies. However, it does not find out the corresponding mathematical object to the source $J_\mu^K$. In this section, we give an answer to this problem, and show the affine connection representation of Yang-Mills equation. 

In order to obtain the general Yang-Mills equation with gravitation, we have to adopt holonomic connection to construct it. Suppose $\mathcal{F}$ evolves in $\mathcal{G}$ according to Definition \ref{DefinitionOfMotionAndInteraction}, that is, $\forall p\in M$, 
$$(U,\alpha^{A'})\xrightarrow{\mathfrak{f}(p)}(U,\xi^A)\xrightarrow{f(p)}(U,x^M)\xleftarrow{g(p)}(U,\zeta^A)\xleftarrow{\mathfrak{g}(p)}(U,\beta^{A'}).$$ 
We always take the following notations in the coordinate frame $(U,x^M)$.

(i) Let the holonomic connections, which are defined by Eq.(\ref{DefinitionFormulaOfHolonomicConnection}), of geometric manifolds $(M,\mathcal{F})$ and $(M,\mathcal{G})$ be $(\Gamma_\mathcal{F})_{NP}^M$ and $(\Gamma_\mathcal{G})_{NP}^M$, respectively. The colon ":" and the semicolon ";" are used to express the covariant derivatives on $(M,\mathcal{F})$ and $(M,\mathcal{G})$, respectively, e.g.:
$${u^Q}_{:P}=\frac{\partial u^Q}{\partial x^P}+ (\Gamma_\mathcal{F})_{HP}^Q u^H,\ \ \ \ \ {u^Q}_{;P}=\frac{\partial u^Q}{\partial x^P}+ (\Gamma_\mathcal{G})_{HP}^Q u^H.$$

(ii) Let the coefficients of curvature tensor of $(M,\mathcal{F})$ and $(M,\mathcal{G})$ be $K_{NPQ}^M$ and $R_{NPQ}^M$, respectively, i.e.
\begin{equation}\label{CurvatureTensorCoefficients}\LeftList
K_{NPQ}^M \DenoteBy\frac{\partial(\Gamma_\mathcal{F})_{NQ}^M}{\partial x^P}-\frac{\partial(\Gamma_\mathcal{F})_{NP}^M}{\partial x^Q}+(\Gamma_\mathcal{F})_{NQ}^H(\Gamma_\mathcal{F})_{HP}^M  -(\Gamma_\mathcal{F})_{NP}^H(\Gamma_\mathcal{F})_{HQ}^M, \hfill\\
R_{NPQ}^M \DenoteBy\frac{\partial(\Gamma_\mathcal{G})_{NQ}^M}{\partial x^P}-\frac{\partial(\Gamma_\mathcal{G})_{NP}^M}{\partial x^Q}+(\Gamma_\mathcal{G})_{NQ}^H(\Gamma_\mathcal{G})_{HP}^M  -(\Gamma_\mathcal{G})_{NP}^H(\Gamma_\mathcal{G})_{HQ}^M. \hfill\\
\RightList\end{equation}

Denote ${K_{NPQ}^M}^{:P}\DenoteBy (G_\mathcal{F})^{PP'}{K_{NPQ\,:P'}^M}$. On an arbitrary evolution path $L$, we define 
$$\rho_{N0}^M dx^0 \DenoteBy \pi^{*}\left({K_{NPQ}^M}^{:P}dx^Q\right) \in T^*(L).$$
Then, according to Definition \ref{DefinitionOfEquivalenceAndHomomorphism} and the evolution lemma of Proposition \ref{PropositionOfEvolutionLemma}, we obtain $\rho_{N0}^M = {K_{NPQ}^M}^{:P}\varepsilon_0^Q$ and 
$${K_{NPQ}^M}^{:P}dx^Q \TongTaiYu\rho_{N0}^M dx^0.$$
Let $\nabla t = {K_{NPQ}^M}^{:P}\frac{\partial}{\partial x_Q}$. Then, according to Proposition \ref{TheoremOfActualEvolution}, if and only if $\forall p\in M$, $[L_p]=\nabla t|_p$, we have 
$${K_{NPQ}^M}^{:P}\frac{\partial}{\partial x_Q}\DengJiaYu\rho_{N0}^M\frac{d}{{dx_0}}.$$
Applying the evolution lemma of Proposition \ref{PropositionOfEvolutionLemma} again, we obtain  
$${K_{NPQ}^M}^{:P}=\rho_{N0}^M\bar\varepsilon_Q^0.$$ 
Denote $j_{NQ}^M \DenoteBy\rho_{N0}^M\bar\varepsilon_Q^0$; then, if and only if $[L_p]=\nabla t|_p$ we have
\begin{equation}\label{YangMillsEquation}{K_{NPQ}^M}^{:P}=j_{NQ}^M,\end{equation}
which is said to be {\bf (affine) Yang-Mills equation} of $\mathcal{F}$. It contains effects of gravitation, which makes the gauge frames $(B_f)^A_M$ and $(C_f)^M_A$ have physical effects. According to Eq.(\ref{HolonomicCurvatureFrameTransformation}), we know Eq.(\ref{YangMillsEquation}) is coordinate covariant, and if gravitation is removed, it is also gauge covariant.

Thus, we have the following two results.

(i) The Yang-Mills equation originates from a geometric property in the direction $\nabla t$. In other words, the on-shell evolution of gauge field is described by the direction field $\nabla t$.

(ii) We obtain the mathematical origination of charge and current. We know that the evolution path $L$ is an imbedding submanifold of $M$. Thus, the charge $\rho_{N0}^M$ originates from the pull-back $\pi^{*}$ from $M$ to $L$, and the current $j_{NQ}^M$ originates from $\nabla t$ that is associated to $\rho_{N0}^M$. 

If we let $(M,f)$ be completely flat, i.e., $(B_f)^A_M=\delta^A_M$, $(C_f)^M_A=\delta^M_A$, then by calculation, we find $\rho_{N0}^M$ can still be non-vanishing. This shows that $\rho_{N0}^M$ originates from $(M,\mathfrak{f})$ ultimately.

\sDefinition\label{DefinitionOfGeneralCharge} We speak of the real-valued
\begin{equation}\rho_{MN0} \DenoteBy G_{MH}\rho_{N0}^{H}\end{equation}
as the {\bf field function of a general charge} or speak of it as a {\bf charge} of $\mathcal{F}$ for short.

\SubSectionWithLabel{On-shell evolution of general charge and affine connection representation of mass, energy, momentum, and action}\label{SectionOfGeneralEnergeMomentumEquation}\label{SectionOfAffineAction}

In order to be compatible with the affine connection representation of gauge fields, we also have to define mass, energy, momentum, and action in the form associated to affine connection. We are going to show them in this section and section \ref{Affine connection representation of classical spacetime evolution}. 

Let $\mathbf{F}_0 \DenoteBy \rho_{MN0} dx^M \otimes dx^N$. For the sake of simplicity, denote the charge $\rho_{MN0}$ of $\mathcal{F}$ by $\rho_{MN}$ concisely. Let $D$ be the holonomic connection of $(M,\mathcal{G})$; then,
$$D\mathbf{F}_0 \DenoteBy D\rho_{MN} \otimes dx^M \otimes dx^N,\ \ \ \ \ \ \ \ \nabla \mathbf{F}_0 \DenoteBy \nabla\rho_{MN}\otimes dx^M \otimes dx^N,$$
where $D\rho_{MN}\DenoteBy\rho_{MN;Q}dx^Q$ and $\nabla\rho_{MN}\DenoteBy\rho_{MN;Q}\frac{\partial}{\partial x_Q}$. According to Proposition \ref{TheoremOfActualEvolution}, if and only if $\forall p\in L$, the evolution direction is taken as $[L_p]=\nabla\rho_{MN}|_p$, we have
$$D\rho_{MN} \TongTaiYu D_L\rho_{MN}, \ \ \ \ \ \ \ \ \nabla\rho_{MN}\DengJiaYu \nabla_L \rho_{MN},$$
that is,
\begin{equation}
\rho_{MN;Q}dx^Q \TongTaiYu\rho_{MN;0}dx^0 ,\ \ \ \ \ \ \ \ \rho_{MN;Q}\frac{\partial}{\partial x_Q}\DengJiaYu\rho_{MN;0}\frac{d}{dx_0}.
\end{equation}

\sDefinition\label{DefinitionOfEnergyMomentum} For more convenience, the notation $\rho_{MN}$ is further abbreviated   as $\rho$. In affine connection representation, energy and momentum of $\rho$ are defined as
\begin{equation}\begin{aligned}
&E_0 \DenoteBy\rho_{;0}\DenoteBy\rho_{;Q}\varepsilon_0^Q,\ \ \ \ p_Q \DenoteBy\rho_{;Q},\ \ \ \ H_0 \DenoteBy\frac{d\rho}{dx^0},\ \ \ \ P_Q \DenoteBy\frac{\partial\rho}{\partial x^Q}, \\
&E^0 \DenoteBy\rho^{;0}\DenoteBy\rho^{;Q}\bar\varepsilon_Q^0, \ \ \ \ p^Q \DenoteBy\rho^{;Q},\ \ \ \ H^0 \DenoteBy\frac{{d\rho}}{dx_0},\ \ \ \ P^Q \DenoteBy\frac{\partial\rho}{\partial x_Q}. \\
\end{aligned}\end{equation}

\sProposition\label{PropositionOfGeneralEnergyMomentumEquations} At any point $p$ on $M$, the equation
\begin{equation}\label{EnergyMomentumEquation}E_0 E^0 =p_Q p^Q\end{equation}
holds if and only if the evolution direction $[L_p]=\nabla\rho|_p$. Eq.(\ref{EnergyMomentumEquation}) is the (affine) energy-momentum equation of $\rho$.

\sProof According to the above discussion, $\forall p \in M$, $[L_p]=\nabla\rho|_p$ is equivalent to
\begin{equation} \label{EqEnergyMomentum}
p_Q dx^Q \TongTaiYu E_0 dx^0,\ \ \ \ \ \ \ \ p_Q\frac{\partial}{\partial x_Q}\DengJiaYu E_0\frac{d}{dx_0}.
\end{equation}
Then, due to Proposition \ref{Evolution lemma}.\ref{PropositionConjugationofTangentAndCotangent}, we obtain the directional derivative in the gradient direction $\nabla\rho$:
$$\LeftAngle{p_Q\frac{\partial}{\partial x_Q},p_M dx^M}\RightAngle = \LeftAngle{E_0\frac{d}{dx_0},E_0 dx^0}\RightAngle,$$
i.e., $G^{QM}p_Q p_M = G^{00}E_0 E_0$, or $p_Q p^Q = E_0 E^0$.
\qed

\sProposition\label{PropositionOfTraditinoalDefinitionOfMomentum} At any point $p$ on $M$, the equations
\begin{equation} \label{TraditinoalDefinitionOfMomentum}
  p^Q =E^0\frac{dx^Q}{dx^0},\ \ \ p_Q =E_0\frac{dx_Q}{dx_0}
\end{equation}
hold if and only if the evolution direction $[L_p]=\nabla\rho|_p$.

\sProof Due to the evolution lemma of Proposition \ref{PropositionOfEvolutionLemma}, we immediately obtain Eq.(\ref{TraditinoalDefinitionOfMomentum}) from Eq.(\ref{EqEnergyMomentum}).
\qed

\noindent {\bf Remark.} In the gradient direction $\nabla\rho$, Eq.(\ref{TraditinoalDefinitionOfMomentum}) is consistent with the conventional formula
$$p=mv.$$ 
Thus, in affine connection representation, the energy-momentum equation and the conventional definition of momentum both originate from a geometric property in gradient direction. In other words, the on-shell evolution of the particle field $\rho$ is described by the gradient direction field $\nabla \rho$.

\sDefinition\label{AffineAction} Let $\mathcal{P}(b,a)$ be the totality of paths from $a$ to $b$. And suppose $L \in\mathcal{P}(b,a)$, and the evolution parameter $x^0$ satisfies $t_a \DenoteBy x^0(a)< x^0(b)\DenoteBy t_b$. The elementary affine action of $\rho$ is defined as
\begin{equation}\label{ElementaryAffineAction}\mathfrak{s}(L)\DenoteBy\int_{L}{D\rho}=\int_{L}{p_Q dx^Q}=\int_{t_a}^{t_b}{E_0 dx^0}.\end{equation}
Thus, $\delta \mathfrak{s}(L)=0$ if and only if $L$ is a gradient line of $\rho$.

In particular, in the case where $\mathcal{G}$ is orthogonal, we can also define action in the following way. 

On $(M,\mathcal{G})$ let there be Dirac algebras $\gamma^M$ and $\gamma_N$ such that
$$\gamma^M \gamma^N +\gamma^N \gamma^M  =2G^{MN},\ \ \ \ \ \ \gamma_M \gamma_N +\gamma_N \gamma_M  =2G_{MN},\ \ \ \ \ \ \gamma_M\gamma^M=1.$$
In a gradient direction of $\rho$, from Eq.(\ref{EnergyMomentumEquation}), we obtain that
$$\begin{aligned}
p_Q p^Q = E_0 E^0 &\ \IfAndOnlyIf\ \rho_{;Q}\rho^{;Q} = \rho_{;0}\rho^{;0} \\
&\ \IfAndOnlyIf\ G^{PQ}\rho_{;P}\rho_{;Q}=G^{00}\rho_{;0}\rho_{;0} \\
&\ \IfAndOnlyIf\ (\gamma^P \gamma^Q +\gamma^Q \gamma^P)\rho_{;P}\rho_{;Q}=2\rho_{;0}\rho_{;0} \\
&\ \IfAndOnlyIf\ (\gamma^P\rho_{;P})(\gamma^Q\rho_{;Q}) + (\gamma^Q\rho_{;Q})(\gamma^P\rho_{;P})=2\rho_{;0}\rho_{;0} \\
&\ \IfAndOnlyIf\ (\gamma^P\rho_{;P})^2 = (\rho_{;0})^2.
\end{aligned}$$
Take $\gamma^P\rho_{;P} = \rho_{;0}$ without loss of generality, and then, in the gradient direction of $\rho$, we have
\begin{equation}\label{GammaEvolution}\gamma^P\rho_{;P}dx^0 = \rho_{;0}dx^0 = \varepsilon^P_0 \rho_{;P}dx^0 = D\rho.\end{equation}
So we can take 
\begin{equation}\label{OrthogonalAffineAction}
s(L)\DenoteBy\int_{L}\left(\gamma^P\rho_{;P}dx^0 + D\rho\right) = \int_{t_a}^{t_b}\left(\gamma^P\rho_{;P} + \varepsilon^P_0\rho_{;P}\right)dx^0 = \int_{t_a}^{t_b}{\left(\gamma^P\rho_{;P} + E_0\right) dx^0}.
\end{equation}
Remark \ref{ApplicableForWaveAndField} and Remark \ref{WaveFunctionAndFieldFunctionDetail} explain the rationality of this definition. We have $s(L)=2\mathfrak{s}(L)$ in the gradient direction of $\rho$, so $\mathfrak{s}(L)$ and $s(L)$ are consistent.

\sRemark\label{ApplicableForWaveAndField} In the Minkowski coordinate frame of section \ref{Classical spacetime reference-system}, the evolution parameter $x^0$ is replaced by $\tilde x^\tau$; then, there still exists a concept of gradient direction $\tilde\nabla \tilde\rho$. Correspondingly, Eq.(\ref{ElementaryAffineAction}) and Eq.(\ref{OrthogonalAffineAction}) present as 
$$\tilde{\mathfrak{s}}(L)\DenoteBy\int_{L}{\tilde D\tilde\rho}=\int_{L}{\tilde p_\mu d\tilde x^\mu}=\int_{\tau_a}^{\tau_b}{\tilde m_\tau d\tilde x^\tau},\ \ \ \ \ \ \ \ \tilde s(L) = \int_{\tau_a}^{\tau_b}{\left(\gamma^\mu \tilde\rho_{;\mu} + \tilde m_\tau \right) d\tilde x^\tau},$$
where $\tilde m_\tau$ is the rest-mass and $\tilde x^\tau$ is the proper-time. 

\newpage

\sRemark\label{RemarkAffineConnectionRepresentationOfInteraction} Define the following notations.
$${[\rho \Gamma_G]}\DenoteBy\frac{\partial\rho_{MN}}{\partial x^G}-\rho_{MN;G}=\rho_{MH}\Gamma_{NG}^H +\rho_{HN}\Gamma_{MG}^H,\ \ \ \ \ \ \ \ \ {[\rho R_{PQ}]}\DenoteBy\rho_{MH}R_{NPQ}^H  +\rho_{HN}R_{MPQ}^H.$$
Then, through some calculations, we can obtain that
$$f_P \DenoteBy p_{P;0}=E_{0;P} - p_Q\varepsilon_{0;P}^Q  + [{\rho R_{PQ}}]\varepsilon_0^Q,$$
which is the affine connection representation of general Lorentz force equation. See Discussion \ref{DiscussionOfMinkowskiEnergyMomentumEquation} for further illustrations.

\SubSectionWithLabel{Inversion transformation in affine connection representation}

In affine conection representation, $CPT$ inversion is interpreted as a full inversion of coordinates and metrics. Let $i,j=1,2,3$ and $m,n=4,5,\cdots,\mathfrak{D}$.

Let the local coordinate representation of reference-system $k$ be ${x'}^j = -\delta_i^j x^i$, ${x'}^n =\delta_m^n x^m$; then, parity inversion can be represented as 
$$P\DenoteBy L_k: x^i \to -x^i, x^m \to x^m.$$ 
Let the local coordinate representation of reference-system $h$ be ${x'}^j =\delta_i^j x^i$, ${x'}^n = -\delta_m^n x^m$; then, charge conjugate inversion can be represented as 
$$C\DenoteBy L_h: x^i \to x^i, x^m \to -x^m.$$ 
Time coordinate inversion can be represented as 
$$T_0: x^0 \to -x^0.$$ 
{\bf Full inversion of coordinates} can be represented as 
\begin{equation}CPT_0: x^Q \to -x^Q, x^0 \to -x^0.\end{equation}

The positive or negative sign of metric marks two opposite directions of evolution. Let $N$ be a closed submanifold of $M$, and let its metric be $dx^{(N)}$. Denote the totality of closed submanifolds of $M$ by $\mathfrak{B}(M)$; then, {\bf full inversion of metrics} can be expressed as 
\begin{equation}T^{(M)}\DenoteBy\prod_{N\in\mathfrak{B}(M)}\left({dx^{(N)}\to - dx^{(N)}}\right).\end{equation}

Denote time inversion by 
$$T\DenoteBy T^{(M)}T_0,$$ 
and then, the joint transformation of the full inversion of coordinates $CPT_0$ and the full inversion of metrics $T^{(M)}$ is 
\begin{equation}(CPT_0)(T^{(M)})=CPT.\end{equation}
Summerize the above discussions; then, we have
$$\LeftList
\begin{aligned}
CPT_0&: x^Q \to -x^Q,\ x^0 \to -x^0,\ dx^Q \to dx^Q,\ dx^0 \to dx^0, \hfill\\
T^{(M)}&: x^Q \to x^Q,\ x^0 \to x^0,\ dx^Q \to -dx^Q,\ dx^0 \to -dx^0, \hfill\\
CPT&: x^Q \to -x^Q,\ x^0 \to -x^0,\ dx^Q \to -dx^Q,\ dx^0 \to -dx^0. \hfill\\
\end{aligned}
\RightList$$
The $CPT$ invariance in affine connection representation is very clear. Concretely, on $(M,\mathcal{G})$, we consider the $CPT$ transformation acting on $\mathcal{G}$. Denote $\displaystyle s\DenoteBy \int_L D\rho$ and $D_Pe^{is}\DenoteBy \left(\frac{\partial}{\partial x^P}-i[\rho\Gamma_P]\right)e^{is}$; then, through simple calculations, we obtain that
$$CPT: D\rho \to D\rho,\ \ \ \ D_Pe^{is} \to -D_Pe^{-is}.$$

\sRemark In quantum mechanics there is a complex conjugation in the time inversion of wave function $T: \psi(x,t) \to \psi^*(x,-t)$. In affine connection representation, we know the complex conjugation can be interpreted as a straightforward mathematical result of the full inversion of metrics $T^{(M)}$.

\SubSectionWithLabel{Two dual descriptions of gradient direction field}

\sDiscussion Let $X$ and $Y$ be non-vanishing smooth tangent vector fields on the manifold $M$. And let $L_Y$ be the Lie derivative operator induced by the one-parameter group of diffeomorphism $\varphi_Y$. Then, according to a well-known theorem\cite{Chern}, we obtain the Lie derivative equation 
\begin{equation}\label{LieEquation}[{X,Y}]=L_Y X.\end{equation}
Suppose $\forall p \in M$, $Y(p)$ is a unit-length vector, i.e., $\vert\vert Y(p) \vert\vert = 1$. Let the parameter of $\varphi_Y$ be $x^0$. Then, on the evolution path $L \DenoteBy \varphi_{Y,p}$, we have
\begin{equation}\label{YEquivalentEquation}Y \DengJiaYu \frac{d}{dx^0}.\end{equation}
Thus, Eq.(\ref{LieEquation}) can also be represented as 
\begin{equation}\label{PreHeisenbergEquation}[{X,Y}]=\frac{d}{dx^0} X.\end{equation}
On the other hand, $\forall df \in T(M)$ and $df_L\DenoteBy\pi^*(df)$, due to (\ref{YEquivalentEquation}) and Proposition \ref{Evolution lemma}.\ref{PropositionConjugationofTangentAndCotangent} we have $\LeftAngle{Y,df}\RightAngle =\LeftAngle{\frac{d}{dx^0} ,df_L}\RightAngle$, that is
\begin{equation}\label{PreSchrodingerEquation}Yf=\frac{d}{dx^0}f_L.\end{equation}

\sDefinition Let $H \DenoteBy \vert\vert\nabla\rho\vert\vert^{-1}\nabla\rho = \varepsilon^M_0 \frac{\partial}{\partial x^M} \DengJiaYu \frac{d}{dx^0}$. It is evident that $\forall p \in M$, $\vert\vert H(p) \vert\vert = 1$. If and only if taking $Y=H$, we speak of (\ref{PreHeisenbergEquation}) and (\ref{PreSchrodingerEquation}) as real-valued (affine) Heisenberg equation and (affine) Schr{\"{o}}dinger equation, respectively, that is,
\begin{equation}\label{HeisenbergEquationAndSchrodingerEquation}
{[{X,H}]}=\frac{d}{dx^0}X,\ \ \ \ Hf=\frac{d}{dx^0}f_L.
\end{equation}

\sDiscussion The above two equations both describe the gradient direction field and thereby reflect on-shell evolution. Such two dual descriptions of gradient direction show the real-valued affine connection representation of Heisenberg picture and Schr{\"{o}}dinger picture.

It is not hard to find out several different kinds of complex-valued representations of gradient direction. For examples, one is the affine Dirac equation in section \ref{Affine connection representation of Dirac equation}, and another is as follows. 

Let $\psi \DenoteBy fe^{is_L}$, where it is fine to take either $s_L \DenoteBy s(L)$ or $s_L \DenoteBy \mathfrak{s}(L)$ from Definition \ref{AffineAction}. According to Eq.(\ref{HeisenbergEquationAndSchrodingerEquation}), it is easy to obtain on $L$ that
\begin{equation}\label{ComplexValuedEquations}
{[{X,H}]}=\frac{d}{dx^0}X,\ \ \ \ H\psi=\frac{d\psi}{dx^0}.
\end{equation}
This is consistent with the conventional Heisenberg equation and Schr{\"{o}}dinger equation (taking the natural units that $\hbar=1,\ c=1$)
\begin{equation}\label{ConventinalEquations}
{[{X,-iH}]}=\frac{\partial}{\partial t}X,\ \ \ \ -iH\psi=\frac{\partial\psi}{\partial t},
\end{equation}
and they have a coordinate correspondence
$$\frac{\partial}{\partial (ix^k)} \leftrightarrow \frac{\partial}{\partial x^k},\ \ \ \frac{\partial}{\partial t} \leftrightarrow \frac{d}{dx^0}.$$
We know that $\frac{\partial}{\partial t} \leftrightarrow \frac{d}{dx^0}$ originates from the difference that the evolution parameter is $x^\tau$ or $x^0$. The imaginary unit $i$ originates from the difference between the regular coordinates $x^1,x^2,x^3,x^\tau$ and the Minkowski coordinates $x^1,x^2,x^3,x^0$. That is to say, the regular coordinates satisfy 
$$(dx^0)^2=(dx^1)^2+(dx^2)^2+(dx^3)^2+(dx^\tau)^2,$$ 
and the Minkowski coordinates satisfy
$$(dx^\tau)^2=(dx^0)^2-(dx^1)^2-(dx^2)^2-(dx^3)^2 = (dx^0)^2+(d(ix^1))^2+(d(ix^2))^2+(d(ix^3))^2.$$
This causes the appearance of the imaginary unit $i$ in the correspondence 
$$ix^k \leftrightarrow x^k.$$
So Eq.(\ref{ComplexValuedEquations}) and Eq.(\ref{ConventinalEquations}) have exactly the same essence, and their differences only come from different coordinate representations. 

The differences between coordinate representations have nothing to do with the geometric essence and the physical essence. We notice that the value of a gradient direction is dependent on geometry, but independent of that the equations are real-valued or complex-valued. Therefore, it is unnecessary for us to confine to such algebraic forms as real-valued or complex-valued forms, but we should focus on such geometric essence as gradient direction.

The essential virtue of complex-valued form is that it is applicable for describing the coherent superposition of propagator. However, this is independent of the above discussions, and we are going to discuss it in section \ref{Quantum evolution as a distribution of gradient directions}. 

\SubSectionWithLabel{Quantum evolution as a distribution of gradient directions}

From Proposition \ref{PropositionOfGeneralEnergyMomentumEquations}, we see that, in affine connection representation, the classical on-shell evolution is described by gradient direction. Then, naturally, quantum evolution should be described by the distribution of gradient directions. 

The distribution of gradient directions on a geometric manifold $(M,\mathcal{G})$ is effected by the bending shape of $(M,\mathcal{G})$; in other words, the distribution of gradient directions can be used to reflect the shape of $(M,\mathcal{G})$. This is the way that the quantum theory in affine connection representation describes physical reality. 

In order to know the full picture of physical reality, it is necessary to fully describe the shape of the geometric manifold. For a single observation,

(1) It is the reference-system, not a point, that is used to describe the physical reality, so the coordinate of an individual point is not enough to fully describe the location information about the physical reality.

(2) Through a single observation of momentum, we can only obtain information about an individual gradient direction; this cannot reflect the full picture of the shape of the geometric manifold. 

Quantum evolution provides us with a guarantee that we can obtain the distribution of gradient directions through multiple observations, so that we can describe the full picture of the shape of the geometric manifold. 

Next, we are going to carry out strict mathematical descriptions for the quantum evolution in affine connetion representation. 

\sDefinition\label{KernalOfCharge} Let $\rho$ be a geometric property on $M$, such as a charge of $\mathcal{F}$. Then $H\DenoteBy\nabla\rho$ is a gradient direction field of $\rho$ on $(M,\mathcal{G})$.

Let $\mathfrak{T}$ be the totality of all flat transformations $L_k$ defined in section \ref{Gauge transformation in affine connection representation}. $\forall T\in\mathfrak{T}$, the flat transformation $T:f\mapsto Tf$ induces a transformation $T^*:\rho\mapsto T^*\rho$. Denote 
$$|\rho|\DenoteBy\{{\rho_T\DenoteBy T^*\rho|{T\in\mathfrak{T}}}\},\ \ \ \ \ \ |H |\DenoteBy\{{H_T\DenoteBy \nabla\rho_T|{T\in\mathfrak{T}}}\}.$$
$\forall a\in M$, the restriction of $|H|$ at $a$ are denoted by $|H(a)|\DenoteBy\{{H_T(a)|{T\in\mathfrak{T}}}\}$.

We say $|H|$ is the {\bf total distribution} of the gradient direction field $H$. 

\sRemark When $T$ is fixed, $H_T$ can reflect the shape of $(M,\mathcal{G})$. When $a$ is fixed, the extension to $|H(a)|$ can reflect the shape of $(M,\mathcal{G})$. 

However, when $T$ and $a$ are both fixed, $H_T(a)$ is a fixed individual gradient direction, which cannot reflect the shape of $(M,\mathcal{G})$. In other words, if the momentum $p_T$ and the position $x_a$ of $\rho$ are both definitely observed, the physical reality $\mathcal{G}$ would be unknowable; therefore this is unacceptable. This is an embodiment of quantum uncertainty in affine connection representation.

\sDefinition Let $\varphi_H$ be the one-parameter group of diffeomorphisms corresponding to $H$. The parameter of $\varphi_H$ is $x^0$. $\forall a\in M$, according to Definition \ref{DefinitionOfEvolutionPath}, let $\varphi_{H,a}$ be the evolution path through $a$, such that $\varphi_{H,a}(0)=a$. $\forall t\in\mathbb{R}$, denote 
$$\varphi_{|H|,a}\DenoteBy\{{\varphi_{X,a}|{X\in |H|}}\},\ \ \ \ \ \ \varphi_{|H|,a}(t)\DenoteBy\{{\varphi_{X,a}(t)|{X\in |H|}}\}.$$
$\forall\Omega \subseteq\mathfrak{T}$, we also denote $|{H_\Omega}|\DenoteBy\{{H_T|{T\in\Omega}}\} \subseteq |H|$ and 
$$\varphi_{|{H_\Omega }|,a}\DenoteBy\{{\varphi_{X,a}|{X\in |{H_\Omega}|}}\} \subseteq \varphi_{|H|,a},\ \ \ \ \ \  \varphi_{|{H_\Omega}|,a}(t)\DenoteBy\{{\varphi_{X,a}(t)|{X\in |{H_\Omega}|}}\} \subseteq \varphi_{|H|,a}(t).$$
$\forall a\in M$, the restriction of $|{H_\Omega}|$ at $a$ are denoted by $|{H_\Omega (a)}|\DenoteBy\{{H_T(a)|{T\in\Omega}}\}$.

\sRemark\label{RemarkDistributionOfEvolution} At the beginning $t=0$, intuitively, the gradient directions $|{H(a)}|$ of $|\rho|$ start from $a$ and point to all directions around $a$ uniformly. If $(M,\mathcal{G})$ is not flat, when evolving to a certain time $t > 0$,  the distribution of gradient directions on $\varphi_{|H|,a}(t)$ is no longer as uniform as beginning. The following definition precisely characterizes this kind of ununiformity.

\sDefinition Let the transformation $L_{\mathcal{G}^{-1}}$ act on $\mathcal{G}$, then we obtain the trivial $e \DenoteBy L_{\mathcal{G}^{-1}}(\mathcal{G})$. Now $(M,\mathcal{G})$ is sent to a flat $(M,e)$, and the gradient direction field $|H|$ of $|\rho|$ on $(M,\mathcal{G})$ is sent to a gradient direction field $|O|$ of $|\rho|$ on $(M,e)$. Correspondingly, $\forall t\in\mathbb{R}$, $\varphi_{|H|,a}(t)$ is sent to $\varphi_{|O|,a}(t)$. In a word, $L_{\mathcal{G}^{-1}}$ induces the following two maps:
$$\mathcal{G}_*^{-1}:|H|\to |O|,\ \ \ \ \mathcal{G}_{**}^{-1}:\varphi_{|H|,a}\to\varphi_{|O|,a}.$$
$\forall T\in \mathfrak{T}$, deonte $\mathfrak{N} \DenoteBy \{N \in \mathfrak{T} \ |\ \det N = \det T \}$. Due to $\mathfrak{T}\DengJiaYu GL(\mathfrak{D},\mathbb{R})$, let $\mathfrak{U}$ be a neighborhood of $T$, with respect to the topology of $GL(\mathfrak{D},\mathbb{R})$. 

Take $\Omega = \mathfrak{N} \cap \mathfrak{U}$, then
$$|{O_\Omega}| = \mathcal{G}_*^{-1}\left(|{H_\Omega}|\right), \ \ \ \ \varphi_{|{O_\Omega}|,a} = \mathcal{G}_{**}^{-1}\left(\varphi_{|{H_\Omega}|,a}\right).$$
Let $\mu$ be a Borel measure on the manifold $M$. We know $\forall t\in\mathbb{R}$,
$$\varphi_{|{H_\mathfrak{N} }|,a}\left(t\right) \simeq \varphi_{|{O_\mathfrak{N} }|,a}\left(t\right) \simeq \mathbb{S}^{\mathfrak{D}-1}.$$
Thus, $\varphi_{|{H_\Omega }|,a}\left(t\right) \subseteq \varphi_{|{H_\mathfrak{N} }|,a}\left(t\right)$ and $\varphi_{|{O_\Omega }|,a}\left(t\right) \subseteq \varphi_{|{O_\mathfrak{N} }|,a}\left(t\right)$ are Borel sets, so they are measurable. Denote 
$$\mu_a \left(\varphi_{|{H_\Omega}|,a}\left(t\right)\right)\DenoteBy \mu\left(\mathcal{G}_{**}^{-1}\left(\varphi_{|{H_\Omega}|,a}\left(t\right)\right)\right) = \mu \left(\varphi_{|{O_\Omega}|,a}\left(t\right)\right).$$
When $\mathfrak{U} \to T$, we have $\Omega\to T$, $|{H_\Omega}| \to H_T$, $|{H_\Omega (a)}| \to H_T(a)$, and $\varphi_{|{H_\Omega }|,a}(t) \to b\DenoteBy\varphi_{H_T,a}(t)$.

For the sake of simplicity, denote $L \DenoteBy \varphi_{H_T,a}$. Thus, we have $a=L(0),\ b=L(t)$, and denote $p_a \DenoteBy [L_a] = H_T(a),\ p_b \DenoteBy [L_b] = H_T(b)$.

Because $\mu_a$ is absolutely continuous with respect to $\mu$, Radon-Nikodym theorem\cite{Nikodym1930} ensures the existence of the following limit. The Radon-Nikodym derivative
\begin{equation}
  W_L(b,a) \DenoteBy \frac{d\mu_a}{d\mu_b} \DenoteBy \sLimit_{\mathfrak{U} \to T}\frac{{\mu_a\left(\varphi_{|{H_\Omega }|,a}\left(t\right)\right)}}{{\mu\left(\varphi_{|{H_\Omega }|,a}\left(t\right)\right)}}=\sLimit_{\mathfrak{U} \to T}\frac{{\mu\left(\mathcal{G}_{**}^{-1}\left(\varphi_{|{H_\Omega }|,a}\left(t\right)\right)\right)}}{{\mu\left(\varphi_{|{H_\Omega }|,a}\left(t\right)\right)}} = \sLimit_{\mathfrak{U} \to T}\frac{{\mu\left(\varphi_{|{O_\Omega }|,a}\left(t\right)\right)}}{{\mu\left(\varphi_{|{H_\Omega }|,a}\left(t\right)\right)}}
\end{equation}
is said to be the {\bf distribution density of $|H|$ along $L$ in position representation}.

On a neighborhood $U$ of $a$, $\forall T\in\mathfrak{T}$, denote the normal section of $H_T(a)$ by $N_{H_T,a}$, that is, 
$$N_{H_T,a} \DenoteBy \{n \in U \ |\ H_T(a) \cdot (n-a) = 0 \},\ \ \ \ \ \ N_{H_T,a}(t) \DenoteBy \{\varphi_{H_T,x}(t)\ |\ x \in N_{H_T,a} \}.$$
Thus, $N_{H_T,a} = N_{H_T,a}(0)$ and $N_{H_T,b} \DenoteBy N_{H_T,a}(t)$. If $U \to a$, we have $N_{H_T,a} \to a$ and $N_{H_T,a}(t) \to b\DenoteBy\varphi_{H_T,a}(t)$. The Radon-Nikodym derivative
\begin{equation}
Z_L(b,a) \DenoteBy \frac{d\mu(a)}{d\mu(b)} \DenoteBy \sLimit_{U \to a} \frac{\mu(N_{H_T,a})}{\mu(N_{H_T,b})} = \sLimit_{U \to a} \frac{\mu(N_{H_T,a})}{\mu(N_{H_T,a}(t))}
\end{equation}
is said to be the {\bf distribution density of $|H|$ along $L$ in momentum representation}.

In a word, $W_L(b,a)$ and $Z_L(p_b,p_a)$ describe the density of the gradient lines that are adjacent to $b$ in two different ways. They have the following property that is evidently true.

\sProposition\label{PropositionSimplePropagate} Let $L$ be a gradient line. $\forall a,b,c \in L$ such that $L(x^0_a)=a$, $L(x^0_b)=b$, $L(x^0_c)=c$ and $x^0_b > x^0_c > x^0_a$; then,
\begin{equation}\label{SimplePropagate}
W_L(b,a) = W_L(b,c) W_L(c,a),\ \ \ \ \ \ Z_L(b,a) = Z_L(b,c) Z_L(c,a).
\end{equation}

\sDefinition If $L$ is a gradient line of some $\rho'\in|\rho|$, we also say $L$ is a gradient line of $|\rho|$.

\sRemark For any $a$ and $b$, it anyway makes sense to discuss the gradient line of $|\rho|$ from $a$ to $b$. It is because even if the gradient line of $\rho$ starting from $a$ does not pass through $b$, it just only needs to carry out a certain flat transformation $T$ defined in section \ref{Gauge transformation in affine connection representation} to obtain a $\rho'\DenoteBy T_* \rho$; thus, the gradient line of $\rho'$ starting from $a$ can just exactly pass through $b$. Due to $\rho, \rho' \in |\rho|$, we do not distinguish them, it is just fine to uniformly use $|\rho|$. Intuitively speaking, when $|\rho|$ takes two different initial momentums, $|\rho|$ presents as $\rho$ and $\rho'$, respectively.

\sDiscussion With the above preparations, we obtain a new way to describe the construction of the propagator strictly. 

For any path $L$ that starts at $a$ and ends at $b$, we denote $\displaystyle \vert\vert L \vert\vert \DenoteBy \int_L dx^0$ concisely. Let $\mathcal{P}(b,a)$ be the totality of all the paths from $a$ to $b$. Denote
$$\mathcal{P}(b,x^0_b;a,x^0_a) \DenoteBy \{L \ |\ L \in \mathcal{P}(b,a),\ \vert\vert L \vert\vert = x^0_b-x^0_a \}.$$
$\forall L \in \mathcal P(b,x^0_b;a,x^0_a)$, we can let $L(x^0_a)=a$ and $L(x^0_b)=b$ without loss of generality. Thus, $\mathcal P(b,x^0_b;a,x^0_a)$ is the totality of all the paths from $L(x^0_a)=a$ to $L(x^0_b)=b$.

Abstractly, the propagator is defined as the Green function of the evolution equation. Concretely, the propagator still needs a constructive definition. One method is the Feynman path integral
\begin{equation}\label{FeynmanDefinition}K(b,x^0_b;a,x^0_a) \DenoteBy \int_{\mathcal P(b,x^0_b;a,x^0_a)}{e^{is}dL}.\end{equation}
However, there are so many redundant paths in $\mathcal P(b,x^0_b;a,x^0_a)$ that (i) it is difficult to generally define a measure $dL$ on $\mathcal P(b,x^0_b;a,x^0_a)$, and (ii) it may cause unnecessary infinities when carrying out some calculations.

In order to solve this problem, we try to reduce the scope of summation from $\mathcal P(b,x^0_b;a,x^0_a)$ to $H(b,x^0_b;a,x^0_a)$, where $H(b,x^0_b;a,x^0_a)$ is the totality of all the gradient lines of $|\rho|$ from $L(x^0_a)=a$ to $L(x^0_b)=b$. Thus, Eq. (\ref{FeynmanDefinition}) is turned into
$$K(b,x^0_b;a,x^0_a) = \int_{H(b,x^0_b;a,x^0_a)}{\Psi(L) e^{is}dL}.$$
We notice that as long as we take the probability amplitude $\Psi(L)$ of the gradient line $L$ such that $[\Psi(L)]^2 = W_L(b,a)$ in position representation, or take $[\Psi(L)]^2 = Z_L(b,a)$ in momentum representation, it can exactly be consistent with the Copenhagen interpretation. This provides the following new constructive definition for the propagator. 

\sDefinition\label{DefinitionOfPropagator} Suppose $|\rho|$ is defined as Definition \ref{KernalOfCharge}, and denote $H \DenoteBy \nabla\rho$. 

Let $\mathcal{L}(b,a)$ be the totality of all the gradient lines of $|\rho|$ from $a$ to $b$. Denote 
$$H(b,x^0_b;a,x^0_a) \DenoteBy \{L \ |\ L \in \mathcal{L}(b,a),\ \vert\vert L \vert\vert = x^0_b-x^0_a \}.$$
Let $\mathcal{L}(p_b,p_a)$ be the totality of all the gradient lines of $|\rho|$, whose starting-direction is $p_a$ and ending-direction is $p_b$. Denote
$$H(p_b,x^0_b;p_a,x^0_a) \DenoteBy \{L \ |\ L \in \mathcal{L}(p_b,p_a),\ \vert\vert L \vert\vert = x^0_b-x^0_a \}.$$

Let $dL$ be a Borel measure on $H(b,x^0_b;a,x^0_a)$. In consideration of Remark \ref{WaveFunctionAndFieldFunctionDetail}, we let $s$ be the affine action $s(L)$ in Definition \ref{AffineAction}. We say the geometric property
\begin{equation} \label{StrictPositionPropagator}
K(b,x^0_b;a,x^0_a) \DenoteBy \int_{H(b,x^0_b;a,x^0_a)}{\sqrt{W_L(b,a)}e^{is}}dL\end{equation}
is the {\bf propagator} of $|\rho|$ from $(a,x^0_a)$ to $(b,x^0_b)$ {\bf in position representation}. If we let $dL$ be a Borel measure on $H(p_b,x^0_b;p_a,x^0_a)$, then we say
\begin{equation} \label{StrictMomentumPropagator}
\mathcal{K}(p_b,x^0_b;p_a,x^0_a) \DenoteBy \int_{H(p_b,x^0_b;p_a,x^0_a)}{\sqrt{Z_L(b,a)}e^{is}}dL\end{equation}
is the {\bf propagator} of $|\rho|$ from $(p_a,x^0_a)$ to $(p_b,x^0_b)$ {\bf in momentum representation}.

\sDiscussion Now (\ref{StrictPositionPropagator}) and (\ref{StrictMomentumPropagator}) are strictly defined, but the Feynman path integral (\ref{FeynmanDefinition}) has not been possessed of a strict mathematical definition until now. This makes it impossible at present to obtain (e.g., in position representation) a strict mathematical proof of 
\begin{equation}\label{FeynmanPropogatorDefinitionFormula}
\int_{\mathcal P(b,x^0_b;a,x^0_a)}{e^{is}dL} = \int_{H(b,x^0_b;a,x^0_a)}{\sqrt{W_L(b,a)}e^{is}}dL.
\end{equation}
Fortunately, the following two reasons make us believe that Eq.(\ref{FeynmanPropogatorDefinitionFormula}) is expected to be regarded as a strict definition of Feynman path integral; that is to say, the integral on the right-hand side of "$=$" can be regarded as the strict definition of the notation on the left-hand side of "$=$".

On the one hand, we notice that the distribution densities $W_L(b,a)$ and $Z_L(b,a)$ of gradient directions establish an association between probability interpretation and geometric interpretation of quantum evolution. Therefore, we can base on probability interpretation to intuitively consider both sides of "$=$" in Eq.(\ref{FeynmanPropogatorDefinitionFormula}) as the same thing.

On the other hand, on the condition of Proposition \ref{PropositionSimplePropagate}, denote $H(x^0_c) \DenoteBy \left\{L(x^0_c) \ \left| \ L \in H(b,x^0_b;a,x^0_a) \right. \right\}$; then,
\begin{equation}\label{FullPropogate}
K(b,x^0_b;a,x^0_a) = \int_{H(x^0_c)} K(b,x^0_b;c,x^0_c) K(c,x^0_c;a,x^0_a) dc.
\end{equation}
is expected to be provable according to Eq.(\ref{SimplePropagate}) and Eq.(\ref{StrictPositionPropagator}). However, to obtain a strict proof of Eq.(\ref{FullPropogate}) from Eq.(\ref{SimplePropagate}) and Eq.(\ref{StrictPositionPropagator}) is not a trivial mathematical problem, which is necessary but not easy, and needs more mathematical research.

\sDiscussion\label{QuantizationProcedure} The quantization methods of QFT are successful, and they are also applicable in affine connection representation, but in this paper, we do not discuss them.  We try to propose some more ideas to understand the quantization of field in affine connection representation.

\newpage

(1) If we take 
$$\displaystyle \mathfrak{s} = \int_L D\rho = \int_L p_Q dx^Q = \int_L E_0 dx^0$$
according to Definition \ref{AffineAction}, where $D$ is the holonomic connection of $(M,\mathcal{G})$, then consider the distribution of $H \DenoteBy \nabla \rho$, and we know that
$$K(b,x^0_b;a,x^0_a) \DenoteBy \int_{\nabla \rho(b,x^0_b;a,x^0_a)}{\sqrt{W_L(b,a)}e^{i\mathfrak{s}}}dL,\ \ \ \mathcal{K}(p_b,x^0_b;p_a,x^0_a) \DenoteBy \int_{\nabla \rho(p_b,x^0_b;p_a,x^0_a)}{\sqrt{Z_L(b,a)}e^{i\mathfrak{s}}}dL$$
describe the quantization of energy-momentum. Every gradient line in $\nabla \rho(b,x^0_b;a,x^0_a)$ corresponds to a set of eigenvalues of energy and momentum. This is consistent with conventional theories, and this inspires us to consider the following new ideas to carry out the quantization of charge and current of gauge field.

(2) In an analogous manner, if we take 
$$\displaystyle \mathfrak{s} = \int_L Dt = \int_L {K_{NPQ}^M}^{:P}dx^Q = \int_L \rho^M_{N0} dx^0,$$
according to section \ref{On-shell evolution of potential field and affine connection representation of Yang-Mills equation}, where $D$ is the holonomic connection of $(M,\mathcal{F})$, then consider the distribution of $H \DenoteBy \nabla t$,
$$K(b,x^0_b;a,x^0_a) \DenoteBy \int_{\nabla t(b,x^0_b;a,x^0_a)}{\sqrt{W_L(b,a)}e^{i\mathfrak{s}}}dL,\ \ \ \mathcal{K}(p_b,x^0_b;p_a,x^0_a) \DenoteBy \int_{\nabla t(p_b,x^0_b;p_a,x^0_a)}{\sqrt{Z_L(b,a)}e^{i\mathfrak{s}}}dL.$$
Denote $H(b,x^0_b ; x^0_c) \DenoteBy \left\{c\in M \ \left| \ \forall L \in H(b,x^0_b;c,x^0_c),\ \vert\vert L \vert\vert = x^0_b-x^0_c \right. \right\}$ and take $H=\nabla t$; then, the wave function $\psi(b,x^0_b)$ that is defined by the equation
\begin{equation}\psi(b,x^0_b) = \int_{H(b,x^0_b ; x^0_c)}K(b,x^0_b;c,x^0_c) \psi(c,x^0_c) dc\end{equation}
describes the quantization of charge and current. It should be emphasized that this is not the quantization of the energy-momentum of the field, but the quantization of the field itself, which presents as quantized (e.g., discrete) charges and currents. 

\SectionWithLabel{Affine connection representation of gauge fields in classical spacetime}

The new framework established in section \ref{The evolution in affine connection representation of gauge fields} is discussed in the $\mathfrak{D}$-dimensional general coordinate $x^M$, which is more general than the $(1+3)$-dimensional conventional Minkowski coordinate $x^\mu$. 

$(dx^0)^2=\Summation_{M=1}^{\mathfrak{D}}{(dx^M)^2}$ is the total metric of internal space and external space, $(dx^\tau)^2=\Summation_{m=4}^{\mathfrak{D}}{(dx^m)^2}$ is the metric of internal space.

(i) The evolution parameter of the $\mathfrak{D}$-dimensional general coordinate $x^M\ (M=1,2,\cdots, \mathfrak{D})$ is $x^0$. The parameter equation of an evolution path $L$ is represented as $x^M=x^M(x^0)$.

(ii) The evolution parameter of the $(1+3)$-dimensional Minkowski coordinate $x^\mu\ (\mu=0,1,2,3)$ is $x^\tau$. The parameter equation of $L$ is represented as $x^\mu=x^\mu(x^\tau)$.

\noindent The coordinate $x^\mu$ works on the $(1+3)$-dimensional classical spacetime submanifold defined as follows.

\SubSectionWithLabel{Classical spacetime submanifold}

Let there be a smooth tangent vector field $X$ on $(M,f)$. If $\forall p\in M$, $X(p)=b^A \left.\frac{\partial}{\partial\xi^A}\right|_p=c^M \left.\frac{\partial}{\partial x^M}\right|_p$ satisfies that $b^a$ are not all zero and $c^m$ are not all zero, where $a,m=r+1,\cdots,\mathfrak{D}$, then we say $X$ is {\bf internal-directed}. For any evolution path $L \DenoteBy \varphi_{X,p}$, we also say $L$ is {\bf internal-directed}.

Suppose $M=P \times N$, $\mathfrak{D}\DenoteBy dim M$ and $r\DenoteBy dim P=3$. $X$ is a smooth tangent vector field on $M$. Fix a point $o \in M$. If $X$ is internal-directed, then there exist a unique $(1+3)$-dimensional imbedding submanifold $\gamma :\tilde M\to M,\ p\mapsto p$ and a unique smooth tangent vector field $\tilde X$ on $\tilde M$ such that: 

(i) $P \times \{o\}$ is a closed submanifold of $\tilde M$. 

(ii) The tangent map $\gamma_* :T(\tilde M)\to T(M)$ satisfies that $\forall q\in \tilde M$, $\gamma_* :\tilde X(q)\mapsto X(q)$. 

\noindent Such an $\tilde M$ is said to be a {\bf classical spacetime submanifold}.

Let $\varphi_X:M \times \mathbb{R} \to M$ and $\varphi_{\tilde X}:\tilde M \times \mathbb{R} \to \tilde M$
be the one-parameter groups of diffeomorphisms corresponding to $X$ and $\tilde X$, respectively. Thus, we have 
$$\varphi_{\tilde X}={\varphi_X}|_{\tilde M \times \mathbb{R}}.$$
So the evolution in classical spacetime can be described by $\varphi_{\tilde X}$. It should be noticed that

\newpage

(i) $\tilde M$ inherits a part of geometric properties of $M$, but not all. The physical properties reflected by $\tilde M$ are incomplete.

(ii) The correspondence between $\tilde X$ and the restriction of $X$ to $\tilde M$ is one-to-one. For convenience, next we are not going to distinguish the notations $X$ and $\tilde X$ on $\tilde M$, but uniformly denote them by $X$.

(iii) An arbitrary path $\tilde L: T\to\tilde M,t\mapsto p$ on $\tilde M$ uniquely corresponds to a path $L\DenoteBy\gamma\circ\tilde L:T\to M,\ \ t\mapsto p$ on $M$. Evidently the image sets of $L$ and $\tilde L$ are the same, that is, $L(T)=\tilde L(T)$. For convenience, later we are not going to distinguish the notations $L$ and $\tilde L$ on $\tilde M$, but uniformly denote them by $L$.

\SubSectionWithLabel{Classical spacetime reference-system}

Let there be a geometric manifold $(M,f)$ and its classical spacetime submanifold $\tilde M$. And let $L \DenoteBy \varphi_{\tilde X,a}$ be an evolution path on $\tilde M$. Suppose $p\in L$ and $U$ is a coordinate neighborhood of $p$. According to Defnition \ref{DefinitionOfPathMapOfCoordinate}, suppose the $f(p)$ on $U$ and the $f_L(p)$ on $U_L\DenoteBy U \cap L$ satisfy that
\begin{equation}f(p)\ :\ \xi^A =\xi^A(x^M)=\xi^A(x^0),\ \ \ \xi^0 =\xi^0(x^0),\ \ \ \ \ \ \ \ A,M=1,2,\cdots ,\mathfrak{D}.\end{equation}
Thus, it is true that: 

(1) There exists a unique local reference-system $\tilde f(p)$ on $\tilde U\DenoteBy U \Jiao \tilde M$ such that
\begin{equation}\label{RegularCoordinateRepresentationDefinition}\tilde f(p)\ :\ \xi^U =\xi^U(x^K)=\xi^U(x^0),\ \ \ \xi^0 =\xi^0(x^0),\ \ \ \ \ \ \ \ U,K=1,2,3,\tau .\end{equation}

(2) If $L$ is internal-directed, then the above coordinate frames $(\tilde U,\xi^U)$ and $(\tilde U,x^K)$ of $\tilde f(p)$ uniquely determine the coordinate frames $(\tilde U,\tilde\xi^\alpha )$ and $(\tilde U,\tilde x^\mu )$ such that
\begin{equation}\label{MinkowskiCoordinateRepresentationDefinition}\tilde f(p)\ :\ \tilde\xi^\alpha  =\tilde\xi^\alpha \left(\tilde x^\mu \right)=\tilde\xi^\alpha \left(\tilde x^\tau \right),\ \ \ \tilde\xi^\tau  =\tilde\xi^\tau \left(\tilde x^\tau \right),\ \ \ \ \ \ \ \ \alpha ,\mu =0,1,2,3\end{equation}
and the coordinates satisfy 
$$\tilde\xi^s =\xi^s,\ \ \tilde\xi^\tau=\xi^\tau,\ \ \tilde\xi^0 =\xi^0,\ \ 
\tilde x^i = x^i,\ \ \tilde x^\tau = x^\tau,\ \ \tilde x^0 = x^0.$$
That is to say, $\tilde f(p)$ is just exactly the reference system in conventional sense, which has two different coordinate representations (\ref{RegularCoordinateRepresentationDefinition}) and (\ref{MinkowskiCoordinateRepresentationDefinition}). 

We speak of
$$\tilde f: \tilde M \to REF_{\tilde M},\ p \mapsto \tilde f(p) \in REF_p$$
as a {\bf classical spacetime reference-system}. Thus, inertial system can be strictly defined as follows, no need for Newton's first law. Suppose we have a geometric manifold $(\tilde M,\tilde g)$. $F_{\tilde g}$ is a transformation induced by $\tilde g$.

(1) If $\tilde\delta_{\alpha\beta}\tilde B_\mu^\alpha \tilde B_\nu^\beta  =\tilde\varepsilon_{\mu\nu}$, then $\tilde g$ is said to be {\bf (Lorentz) orthogonal}. In this case, $F_{\tilde g}$ is just exactly a local Lorentz transformation. 

(2) If $\tilde B_\mu^\alpha$ and $\tilde C_\alpha^\mu$ are constants on $\tilde M$, then $\tilde g$ is said to be {\bf flat}.

(3) If $\tilde g$ is both orthogonal and flat, then $\tilde g$ is said to be an {\bf inertial-system}. In this case, $F_{\tilde g}$ is just exactly a Lorentz transformation.

\sRemark\label{DiscussionOfMetricOfClassicalSpacetime} Due to 
$$\LeftList
 (d\tilde\xi^\tau )^2 =(d\xi^0)^2  -\Summation_{s=1}^3{(d\xi^s)^2}=\tilde \delta_{\alpha\beta}d\tilde\xi^\alpha d\tilde\xi^\beta = \tilde G_{\mu\nu} d\tilde x^\mu d\tilde x^\nu,\ \ \ \ \ \ \ \tilde G_{\mu\nu} \DenoteBy \tilde \delta_{\alpha\beta} \tilde B_\mu^\alpha \tilde B_\nu^\beta ,  \hfill\\
 (d\tilde x^\tau )^2 =(dx^0)^2  -\Summation_{i=1}^3{(dx^i)^2}=\tilde \varepsilon_{\mu\nu} d\tilde x^\mu d\tilde x^\nu = \tilde H_{\alpha\beta} d\tilde\xi^\alpha d\tilde\xi^\beta,\ \ \ \ \ \ \tilde H_{\alpha\beta} \DenoteBy \tilde \varepsilon_{\mu\nu} \tilde C_\alpha^\mu \tilde C_\beta^\nu, \hfill\\ 
\RightList$$
it is easy to know that $\tilde g$ is orthogonal if and only if $d\tilde\xi^\tau  =d\tilde x^\tau$, i.e., $\tilde G_{\tau\tau} \DenoteBy \tilde B^\tau_\tau \tilde B^\tau_\tau = 1,\ \ \tilde G^{\tau\tau} \DenoteBy \tilde C^\tau_\tau \tilde C^\tau_\tau = 1$. It is only in this case that we can denote $d\tilde\xi^\tau$ and $d\tilde x^\tau$ uniformly by $d\tau$; otherwise, we should be aware of the difference between $d\tilde\xi^\tau$ and $d\tilde x^\tau$ in non-trivial gravitational field. No matter whether $\tilde g$ is an inertial-system or not, and whether there is a non-trivial gravitation field or not, $(d\tilde\xi^\tau )^2 =(d\xi^0)^2  -\Summation_{s=1}^3{(d\xi^s)^2}$ and $(d\tilde x^\tau )^2 =(dx^0)^2  -\Summation_{i=1}^3{(dx^i)^2}$ are always both true in their respective coordinate frames.

\sRemark The evolution lemmas in section \ref{Evolution lemma} can be expressed in Minkowski coordinate as follows:

(i) If $\frac{d}{d\tilde t}\DengJiaYu\frac{d}{d\tilde t_L}$ and $d\tilde f\TongTaiYu d\tilde f_L$, then $\LeftAngle{\frac{d}{d\tilde t},d\tilde f}\RightAngle =\LeftAngle{\frac{d}{d\tilde t_L},d\tilde f_L}\RightAngle$.

(ii) The following conclusions are true.
$$\LeftList
  w^\mu \frac{\partial}{{\partial\tilde x^\mu }}\DengJiaYu w^\tau \frac{d}{{d\tilde x^\tau }}\IfAndOnlyIf \ w^\mu  =w^\tau \tilde\varepsilon_\tau^\mu,  \hfill\\
  w_\mu  d\tilde x^\mu  \TongTaiYu w_\tau  d\tilde x^\tau \ \ \IfAndOnlyIf \ \tilde\varepsilon_\tau^\mu  w_\mu  =w_\tau,  \hfill\\ 
\RightList\ \ \ \ \ \ 
\LeftList
 \bar w_\mu \frac{\partial}{{\partial\tilde x_\mu }}\DengJiaYu\bar w_\tau \frac{d}{{d\tilde x_\tau }}\IfAndOnlyIf \ \bar w_\mu  =\bar w_\tau \tilde{\bar\varepsilon}_\mu^\tau,  \hfill\\
 \bar w^\mu  d\tilde x_\mu  \TongTaiYu\bar w^\tau  d\tilde x_\tau  \ \ \IfAndOnlyIf \ \tilde{\bar\varepsilon}_\mu^\tau \bar w^\mu  =\bar w^\tau.  \hfill\\ 
\RightList$$

\SubSectionWithLabel{Affine connection representation of classical spacetime evolution}

Let $\tilde D$ be the holonomic connection on $(\tilde M,\tilde{\mathcal{G}})$, and denote ${\tilde t}_{L;\tau} \DenoteBy \tilde t_{;\sigma}\tilde \varepsilon_\tau^\sigma$; then, the absolute differential and gradient of section \ref{On-shell evolution as a gradient} can be expressed on $\tilde M$ in Minkowski coordinate as
$$\LeftList
\tilde D\tilde t \DenoteBy \tilde t_{;\sigma}d\tilde x^\sigma,\ \ \tilde D_L {\tilde t_L}\DenoteBy {\tilde t}_{L;\tau}d\tilde x^\tau,\\
\tilde \nabla \tilde t \DenoteBy \tilde t_{;\sigma}\frac{\partial}{\partial \tilde x_\sigma},\ \ \tilde \nabla_L {\tilde t_L}\DenoteBy {\tilde t}_{L;\tau}\frac{d}{d\tilde x_\tau}.\\
\RightList$$
Evidently, $\tilde D\tilde t \TongTaiYu \tilde D_L \tilde t_L$ if and only if $L$ is an arbitrary path. $\tilde \nabla \tilde t \DengJiaYu \tilde \nabla_L \tilde t_L$ if and only if $L$ is the gradient line. 

\sDefinition\label{DefinitionOfMinkowskiEnergyMomentum} Similar to section \ref{SectionOfGeneralEnergeMomentumEquation}, suppose a charge $\tilde\rho$ of $\tilde{\mathcal{F}}$ evolves on $(\tilde M,\tilde{\mathcal{G}})$. We have the following definitions.

(1) The geometric properties $\tilde m^\tau  \DenoteBy\tilde\rho^{;\tau}$ and $\tilde m_\tau  \DenoteBy\tilde\rho_{;\tau}$ are said to be the {\bf rest mass} of $\tilde\rho$.

(2) $\tilde p^\mu  \DenoteBy -\tilde\rho^{;\mu}$ and $\tilde p_\mu  \DenoteBy -\tilde\rho_{;\mu}$ are said to be the {\bf energy-momentum} of $\tilde\rho$, and $\tilde E^0 \DenoteBy \tilde\rho^{;0}$, $\tilde E_0 \DenoteBy \tilde\rho_{;0}$ are said to be the {\bf energy} of $\tilde\rho$.

(3) $\tilde M^\tau  \DenoteBy\frac{{d\tilde\rho}}{{d\tilde x_\tau }}$ and $\tilde M_\tau  \DenoteBy\frac{{d\tilde\rho}}{{d\tilde x^\tau }}$ are said to be the {\bf canonical rest mass} of $\tilde\rho$.

(4) $\tilde P^\mu  \DenoteBy -\frac{{\partial\tilde\rho}}{{\partial\tilde x_\mu }}$ and $\tilde P_\mu  \DenoteBy -\frac{{\partial\tilde\rho}}{{\partial\tilde x^\mu }}$ are said to be the {\bf canonical energy-momentum} of $\tilde\rho$, and $\tilde H^0 \DenoteBy  \frac{{\partial\tilde\rho}}{{\partial\tilde x_0 }}$, $\tilde H_0 \DenoteBy  \frac{{\partial\tilde\rho}}{{\partial\tilde x^0 }}$ are said to be the {\bf canonical energy} of $\tilde\rho$.

\sDiscussion\label{DiscussionOfMinkowskiEnergyMomentumEquation}\label{MinkowskiAffineAction} Similar to Proposition \ref{PropositionOfGeneralEnergyMomentumEquations}, $\forall p\in \tilde M$, if and only if the evolution direction $[L_p]=\tilde \nabla\tilde \rho|_p$, the directional derivative is 
$$\LeftAngle{\tilde m_\tau \frac{d}{{d\tilde x_\tau }},\tilde m_\tau  d\tilde x^\tau }\RightAngle =\LeftAngle{\tilde p_\mu \frac{\partial}{{\partial\tilde x_\mu }},\tilde p_\mu  d\tilde x^\mu }\RightAngle,$$
that is $\tilde G^{\tau\tau}\tilde m_\tau \tilde m_\tau  =\tilde G^{\mu\nu}\tilde p_\mu \tilde p_\nu$, or 
$$\tilde m_\tau \tilde m^\tau  =\tilde p_\mu \tilde p^\mu,$$
which is the affine connection representation of energy-momentum equation.

Similar to Proposition \ref{PropositionOfTraditinoalDefinitionOfMomentum}, according to the evolution lemma, $\forall p\in \tilde M$, if and only if the evolution direction $[L_p]=\tilde \nabla\tilde \rho|_p$, we have $\tilde p_\mu  =-\tilde m_\tau \frac{{d\tilde x_\mu }}{{d\tilde x_\tau }}$, that is $\tilde E_0=\tilde m_\tau\frac{{d\tilde x_0 }}{{d\tilde x_\tau }}=\tilde m_\tau\frac{{dx_0 }}{{dx_\tau }}$ and $\tilde p_i  =-\tilde m_\tau \frac{{d\tilde x_i }}{{d\tilde x_\tau }}=\tilde m_\tau \frac{{-d\tilde x_i }}{{d\tilde x_\tau }}=\tilde m_\tau \frac{{dx_i }}{{dx_\tau }}=\tilde E_0 \frac{{dx_i }}{{dx_0 }}$. This can also be regarded as the origin of $p=mv$.

Similar to Remark \ref{RemarkAffineConnectionRepresentationOfInteraction}, denote
$${[{\tilde\rho\tilde\Gamma_\omega }]}\DenoteBy\frac{{\partial\tilde\rho_{\mu\nu}}}{{\partial\tilde x^\omega }}-\tilde\rho_{\mu\nu ;\omega}=\tilde\rho_{\mu\chi}\tilde\Gamma_{\nu\omega}^\chi   +\tilde\rho_{\chi\nu}\tilde\Gamma_{\mu\omega}^\chi, \ \ \ \ \ \ \ \ \ \ \ \ {[{\tilde\rho\tilde R_{\rho\sigma}}]}\DenoteBy\tilde\rho_{\mu\chi}\tilde R_{\nu\rho\sigma}^\chi   +\tilde\rho_{\chi\nu}\tilde R_{\mu\rho\sigma}^\chi.$$
Then for the same reason as Remark \ref{RemarkAffineConnectionRepresentationOfInteraction}, based on Definition \ref{DefinitionOfMinkowskiEnergyMomentum}, we can strictly obtain
\begin{equation}\label{GeneralLorentzEquationUnSimplified}
\tilde f_\rho  \DenoteBy \tilde p_{\rho;\tau} = \tilde m_{\tau;\rho} -\tilde p_\sigma \tilde\varepsilon^\sigma_{\tau;\rho} + [{\tilde\rho\tilde R_{\rho\sigma}}]\tilde\varepsilon_\tau^\sigma.
\end{equation}
In the mass-point model, $\tilde m_{\tau;\rho}$ and $\tilde\varepsilon^\sigma_{\tau;\rho}$ do not make sense, so Eq.(\ref{GeneralLorentzEquationUnSimplified}) turns into 
$$\tilde f_\rho= [{\tilde\rho \tilde R_{\rho\sigma}}]\tilde\varepsilon_\tau^\sigma.$$
This is the affine connection representation of the force of interaction (e.g. the Lorentz force $ \pmb f = q \left( {\pmb E + \pmb v \times \pmb B} \right)\ $ or $\ f_\rho= j^\sigma F_{\rho\sigma}\ $ of the electrodynamics).

Similar to Definition \ref{AffineAction}, let $\tilde{\mathcal{P}}(b,a)$ be the totality of paths on $\tilde M$ from point $a$ to point $b$. And let $L\in\tilde{\mathcal{P}}(b,a)$, and parameter $\tilde x^\tau$ satisfy $\tau_a \DenoteBy\tilde x^\tau (a)<\tilde x^\tau (b)\DenoteBy\tau_b$. The affine connection representation of action in Minkowski coordinates can be defined as
\begin{equation} \label{EqMinkowskiAffineAction}
\tilde{\mathfrak{s}}(L)\DenoteBy\int_{L}{\tilde D\tilde\rho}=\int_{L}{\tilde p_\mu d\tilde x^\mu}=\int_{\tau_a}^{\tau_b}{\tilde m_\tau d\tilde x^\tau},\ \ \ \ \ \ \ \ \tilde s(L) \DenoteBy \int_{\tau_a}^{\tau_b}{\left(\gamma^\mu \tilde\rho_{;\mu} + \tilde m_\tau \right) d\tilde x^\tau}.
\end{equation}
There are more illustrations in Remark \ref{WaveFunctionAndFieldFunctionDetail}.

\SubSectionWithLabel{Affine connection representation of Dirac equation}

\sDiscussion Define Dirac algebras $\gamma^\mu$ and $\gamma^\alpha$ such that
$$\gamma^\mu=\tilde C_\alpha^\mu \gamma^\alpha,\ \ \ \ \ \ \ \ \gamma^\alpha \gamma^\beta   +\gamma^\beta \gamma^\alpha  =2\tilde\delta^{\alpha\beta},\ \ \ \ \ \ \ \ \gamma^\mu \gamma^\nu +\gamma^\nu \gamma^\mu  =2\tilde G^{\mu\nu}.$$
Suppose $(\tilde M,\tilde{\mathcal{G}})$ is orthogonal. According to Remark \ref{DiscussionOfMetricOfClassicalSpacetime}, $\tilde G^{\tau\tau}=1$. Due to Discussion \ref{DiscussionOfMinkowskiEnergyMomentumEquation}, in a gradient direction of $\tilde\rho \DenoteBy \tilde\rho_{\omega\nu}$, we have
$$\begin{aligned}
\tilde\rho_{;\mu}\tilde\rho^{;\mu}=\tilde\rho_{;\tau}\tilde\rho^{;\tau} 
&\ \IfAndOnlyIf \ \tilde G^{\mu\nu}\tilde\rho_{;\mu}\tilde\rho_{;\nu}=\tilde m_\tau ^2 \\
&\ \IfAndOnlyIf \ (\gamma^\mu \tilde\rho_{;\mu})(\gamma^\nu \tilde\rho_{;\nu})+(\gamma^\nu \tilde\rho_{;\nu})(\gamma^\mu \tilde\rho_{;\mu})= 2\tilde m_\tau ^2 \\
&\ \IfAndOnlyIf \ (\gamma^\mu \tilde\rho_{;\mu})(\gamma^\nu \tilde\rho_{;\nu})=\tilde m_\tau ^2 \\
&\ \IfAndOnlyIf \ (\gamma^\mu \tilde\rho_{;\mu})^2=\tilde m_\tau ^2.
\end{aligned}$$
Without loss of generality, take $\gamma^\mu \tilde\rho_{;\mu}= \tilde m_\tau$, that is
\begin{equation}\label{BasicRealValuedAffineDiracEquation}
\gamma^\mu \tilde\rho_{\omega\nu;\mu} = \tilde m_{\omega\nu\tau}\ .
\end{equation}
Next, denote
$${[g\tilde\Gamma_{\mu}]}^{\omega\nu} \DenoteBy \sum_\sigma \tilde G^{\nu\nu'}\tilde\Gamma_{\nu'\sigma\mu} + \sum_\kappa \tilde G^{\omega\omega'}\tilde\Gamma_{\omega'\kappa\mu}\ ,\ \ \ \ \ \ \ \ \tilde D_\mu^{\omega\nu} \DenoteBy \partial_\mu - {[g\tilde\Gamma_{\mu}]}^{\omega\nu}\ .$$
From Eq.(\ref{BasicRealValuedAffineDiracEquation}), it is obtained that
\begin{equation} \label{DiracDeduceProcedure}
\begin{aligned} 
\sum_{\omega,\nu} \gamma^\mu \tilde\rho_{\omega\nu;\mu} = \sum_{\omega,\nu} \tilde m_{\omega\nu\tau} 
&\ \IfAndOnlyIf \ \sum_{\omega,\nu} \gamma^\mu\left( \partial_\mu \rho_{\omega\nu} - \tilde\rho_{\omega\chi} \tilde\Gamma_{\nu\mu}^\chi - \tilde\rho_{\chi\nu} \tilde\Gamma_{\omega\mu}^\chi\right) = \sum_{\omega,\nu} \tilde m_{\omega\nu\tau} \\
&\ \IfAndOnlyIf \ \sum_{\omega,\nu} \gamma^\mu\left( \partial_\mu \rho_{\omega\nu} - \tilde\rho_{\omega\nu} \sum_\sigma \tilde\Gamma_{\sigma\mu}^\nu - \tilde\rho_{\omega\nu} \sum_\kappa \tilde\Gamma_{\kappa\mu}^\omega\right) = \sum_{\omega,\nu} \tilde m_{\omega\nu\tau} \\
&\ \IfAndOnlyIf \ \sum_{\omega,\nu} \gamma^\mu\left( \partial_\mu \rho_{\omega\nu} - \tilde\rho_{\omega\nu} {[g\tilde\Gamma_{\mu}]}^{\omega\nu}\right) = \sum_{\omega,\nu} \tilde m_{\omega\nu\tau} \\
&\ \IfAndOnlyIf \ \sum_{\omega,\nu} \gamma^\mu\left( \partial_\mu - {[g\tilde\Gamma_{\mu}]}^{\omega\nu} \right) \tilde\rho_{\omega\nu} = \sum_{\omega,\nu} \tilde m_{\omega\nu\tau}\ , \\
\end{aligned}
\end{equation}
that is,
\begin{equation} \label{RealValuedAffineDiracEquation}
\sum_{\omega,\nu} \gamma^\mu \tilde D_\mu^{\omega\nu} \tilde\rho_{\omega\nu} = \sum_{\omega,\nu} \tilde m_{\omega\nu\tau}\ ,\ \ \ \ \ \ \ \ \tilde D_\mu^{\omega\nu} \DenoteBy \partial_\mu - {[g\tilde\Gamma_{\mu}]}^{\omega\nu}\ .
\end{equation}
We speak of the real-valued Eq.(\ref{BasicRealValuedAffineDiracEquation}) and (\ref{RealValuedAffineDiracEquation}) as {\bf affine Dirac equations}.

\sDiscussion Next, we construct a kind of complex-valued representation of affine Dirac equation. The restriction of the charge $\tilde\rho_{\omega\nu}$ to $(\tilde U, \tilde x^\mu)$ is a function $\tilde\rho_{\omega\nu}(\tilde x^\mu)$ with respect to the coordinates $(\tilde x^\mu)\DenoteBy(\tilde x^0,\tilde x^1,\tilde x^2,\tilde x^3)$. Let 
$$\tilde{\rm P}_{\omega\nu}(\tilde x^0) \DenoteBy \int_{(\tilde x^1,\tilde x^2,\tilde x^3)}{\tilde\rho_{\omega\nu}(\tilde x^\mu) d^3\tilde x}.$$
Suppose a function $f_{\omega\nu} = f_{\omega\nu}(\tilde x^\mu)$ on $(\tilde U, \tilde x^\mu)$ satisfies that 
$$\tilde\rho_{\omega\nu}=(f_{\omega\nu})^2 \tilde{\rm P}_{\omega\nu},\ \ \ \ \ \ \ \ \int_{(\tilde x^1,\tilde x^2,\tilde x^3)}{(f_{\omega\nu})^2 d^3\tilde x}=1,\ \ \ \ \ \ \ \ \tilde\varepsilon^\mu_\tau \frac{{\partial f_{\omega\nu}}}{{\partial\tilde x^\mu }}= 0,\ \ \ \ \ \ \ \ \gamma^\mu \frac{{\partial f_{\omega\nu}}}{{\partial\tilde x^\mu }}= 0.$$
We define $\psi_{\omega\nu}$ and $\tilde{\mathbb{M}}_{\omega\nu\tau}$ in the following way.
$$\begin{aligned}
&\tilde y_{\omega\nu} \DenoteBy \int_L{d\tilde\rho_{\omega\nu}} = \int_L \frac{d\tilde\rho_{\omega\nu}}{d\tilde x^\tau}d\tilde x^\tau = \int_L \left( \frac{d(f_{\omega\nu}^2)}{d\tilde x^\tau}\tilde{\rm P}_{\omega\nu} + f_{\omega\nu}^2\frac{d\tilde{\rm P}_{\omega\nu}}{d\tilde x^\tau} \right) d\tilde x^\tau = f_{\omega\nu}^2 \int_L \frac{d\tilde{\rm P}_{\omega\nu}}{d\tilde x^\tau} d\tilde x^\tau \DenoteBy f_{\omega\nu}^2 \tilde Y_{\omega\nu}\ ,\\
&\psi_{\omega\nu} \DenoteBy f_{\omega\nu}e^{i\tilde Y_{\omega\nu}},\ \ \ \ \ \ \ \ \ \ \ \ \ \tilde m_{\omega\nu\tau} \DenoteBy \tilde\rho_{\omega\nu;\tau} = (f_{\omega\nu}^2)_{,\tau} \tilde{\rm P}_{\omega\nu} + f_{\omega\nu}^2 \tilde{\rm P}_{\omega\nu;\tau} = f_{\omega\nu}^2 \tilde{\rm P}_{\omega\nu;\tau} \DenoteBy f_{\omega\nu}^2 \tilde{\mathbb{M}}_{\omega\nu\tau}\ . \\
\end{aligned}$$
In the QFT propagator, we usually take $S$ in the path integral $\displaystyle \int e^{iS}\mathcal{D}\psi$ of a fermion in the form of
$$-\int \left(i\bar\psi\gamma^\mu D_\mu\psi-\bar\psi\tilde{\mathbb{M}}_\tau\psi\right) d^4 \tilde x,$$
where $S$ and $d^4 \tilde x$ are both covariant. We believe that the external spatial integral $\displaystyle \int_{(\tilde x^1,\tilde x^2,\tilde x^3)} d^3 \tilde x$ is not an essential part for evolution, so for the sake of simplicity, we do not take into account the external spatial part $\displaystyle \int_{(\tilde x^1,\tilde x^2,\tilde x^3)} d^3 \tilde x$, but only consider the evolution part $\displaystyle \int_L d\tilde x^0$. Meanwhile, in order to remain the covariance, $\displaystyle \int_L d\tilde x^0$ has to be replaced by $\displaystyle \int_L d\tilde x^\tau$. Thus, in affine connection representation of gauge fields, we shall consider an action in the form of
$$-\int_L \left(i\bar\psi\gamma^\mu D_\mu\psi-\bar\psi\tilde{\mathbb{M}}_\tau\psi\right) d\tilde x^\tau.$$
Concretely speaking, denote 
$$\tilde D_{\omega\nu\mu} \DenoteBy \frac{\partial}{{\partial\tilde x^\mu }}- i{[{\tilde{\rm P}\tilde\Gamma_\mu }]}_{\omega\nu},\ \ \ \ \ \ {[\tilde{\rm P}\tilde\Gamma_{\mu}]}_{\omega\nu} \DenoteBy \sum_\sigma \tilde{\rm P}_{\omega\nu}\tilde\Gamma^\nu_{\sigma\mu} + \sum_\kappa \tilde{\rm P}_{\omega\nu}\tilde\Gamma^\omega_{\kappa\mu}\ . $$
From Eq.(\ref{EqMinkowskiAffineAction}), we have 
$$\displaystyle \tilde s_{\omega\nu}(L) \DenoteBy \int_L {\left(\gamma^\mu \tilde\rho_{\omega\nu;\mu} + \tilde m_{\omega\nu\tau} \right) d\tilde x^\tau}.$$ 
And from Eq.(\ref{DiracDeduceProcedure}) we know $\displaystyle \sum_{\omega,\nu}\gamma^\mu \tilde\rho_{\omega\nu;\mu} = \sum_{\omega,\nu}\gamma^\mu\tilde D_\mu^{\omega\nu} \tilde\rho_{\omega\nu}$. Then it is obtained that
\begin{equation}\label{EqTotalAction}
\begin{aligned}
\tilde s_{\tilde\rho}(L) &\DenoteBy \sum_{\omega,\nu} \tilde s_{\omega\nu}(L) 
= \int_L \sum_{\omega,\nu}{\left(\gamma^\mu \tilde\rho_{\omega\nu;\mu} + \tilde m_{\omega\nu\tau} \right)} d\tilde x^\tau 
=\int_L \sum_{\omega,\nu} \left(\gamma^\mu \tilde D_\mu^{\omega\nu} \tilde\rho_{\omega\nu} + \tilde m_{\omega\nu\tau} \right) d\tilde x^\tau \\
&=\int_L \sum_{\omega,\nu} \left(\gamma^\mu\left( \partial_\mu\tilde\rho_{\omega\nu} - {[g\tilde\Gamma_{\mu}]}^{\omega\nu}\tilde\rho_{\omega\nu} \right) + \tilde m_{\omega\nu\tau} \right) d\tilde x^\tau \\
&=\int_L \sum_{\omega,\nu} \left(\gamma^\mu\left( \partial_\mu\tilde y_{\omega\nu} - {[g\tilde\Gamma_{\mu}]}^{\omega\nu}\tilde\rho_{\omega\nu} \right) + \tilde m_{\omega\nu\tau} \right) d\tilde x^\tau \\
&=\int_L \sum_{\omega,\nu} \left(\gamma^\mu\left( \partial_\mu (f_{\omega\nu}^2 \tilde Y_{\omega\nu}) - {[g\tilde\Gamma_{\mu}]}^{\omega\nu}f_{\omega\nu}^2 \tilde{\rm P}_{\omega\nu} \right) + f_{\omega\nu}^2 \tilde{\mathbb{M}}_{\omega\nu\tau} \right) d\tilde x^\tau \\
&=\int_L \sum_{\omega,\nu} \left(\gamma^\mu\left( \partial_\mu \tilde Y_{\omega\nu} - {[\tilde{\rm P}\tilde\Gamma_{\mu}]}_{\omega\nu} \right) f_{\omega\nu}^2 + f_{\omega\nu}^2 \tilde{\mathbb{M}}_{\omega\nu\tau} \right) d\tilde x^\tau \\
&=\int_L \sum_{\omega,\nu} \left(f_{\omega\nu}e^{-i\tilde Y_{\omega\nu}} \gamma^\mu \left( f_{\omega\nu}e^{i\tilde Y_{\omega\nu}} \partial_\mu \tilde Y_{\omega\nu} - {[\tilde{\rm P}\tilde\Gamma_{\mu}]}_{\omega\nu} f_{\omega\nu}e^{i\tilde Y_{\omega\nu}} \right) + f_{\omega\nu}e^{-i\tilde Y_{\omega\nu}} \tilde{\mathbb{M}}_{\omega\nu\tau}f_{\omega\nu}e^{i\tilde Y_{\omega\nu}} \right) d\tilde x^\tau \\
&=\int_L \sum_{\omega,\nu} \left(-\bar\psi_{\omega\nu} i\gamma^\mu \left( e^{i\tilde Y_{\omega\nu}}\partial_\mu f_{\omega\nu} + f_{\omega\nu}e^{i\tilde Y_{\omega\nu}} i\partial_\mu \tilde Y_{\omega\nu} - i{[\tilde{\rm P}\tilde\Gamma_{\mu}]}_{\omega\nu} \psi_{\omega\nu} \right) + \bar\psi_{\omega\nu} \tilde{\mathbb{M}}_{\omega\nu\tau}\psi_{\omega\nu} \right) d\tilde x^\tau \\
&=\int_L \sum_{\omega,\nu} \left(-\bar\psi_{\omega\nu} i\gamma^\mu \left( \partial_\mu (f_{\omega\nu}e^{i\tilde Y_{\omega\nu}}) - i{[\tilde{\rm P}\tilde\Gamma_{\mu}]}_{\omega\nu} \psi_{\omega\nu} \right) + \bar\psi_{\omega\nu} \tilde{\mathbb{M}}_{\omega\nu\tau}\psi_{\omega\nu} \right) d\tilde x^\tau \\
&=\int_L \sum_{\omega,\nu} \left(-\bar\psi_{\omega\nu} i\gamma^\mu \left( \partial_\mu  - i{[\tilde{\rm P}\tilde\Gamma_{\mu}]}_{\omega\nu} \right) \psi_{\omega\nu} + \bar\psi_{\omega\nu} \tilde{\mathbb{M}}_{\omega\nu\tau}\psi_{\omega\nu} \right) d\tilde x^\tau \\
&=\int_L \sum_{\omega,\nu} \left(-\bar\psi_{\omega\nu} i\gamma^\mu \tilde D_{\omega\nu\mu} \psi_{\omega\nu} + \bar\psi_{\omega\nu} \tilde{\mathbb{M}}_{\omega\nu\tau}\psi_{\omega\nu} \right) d\tilde x^\tau 
=-\int_L \sum_{\omega,\nu} \bar\psi_{\omega\nu} \left(i\gamma^\mu \tilde D_{\omega\nu\mu} - \tilde{\mathbb{M}}_{\omega\nu\tau}\right) \psi_{\omega\nu} d\tilde x^\tau .\\
\end{aligned}
\end{equation}
Thus, we have obtained a complex-valued representation of gradient direction of $\tilde\rho_{\omega\nu}$. 

\sRemark\label{WaveFunctionAndFieldFunctionDetail} From the above discussion, we know in the gradient direction of $\rho_{\omega\nu}$, that
$$-\sum_{\omega,\nu} \bar\psi_{\omega\nu} i\gamma^\mu \tilde D_{\omega\nu\mu} \psi_{\omega\nu} d\tilde x^\tau= \sum_{\omega,\nu} \tilde D \tilde\rho_{\omega\nu}.$$
This shows that $s(L)$ and $\tilde s(L)$ in Definition \ref{AffineAction} and Remark \ref{ApplicableForWaveAndField} are indeed applicable for constructing propagator by $e^{is(L)}$ and $e^{i\tilde s(L)}$ in affine connection representation of gauge fields. Therefore, the idea in Discussion \ref{QuantizationProcedure} is reasonable.

\SubSectionWithLabel{From classical spacetime back to full-dimensional space}

\sDiscussion\label{DiscussionFromClassicalToFull} Now there is a problem. $(\tilde M, \tilde{\mathcal{F}})$ and $(\tilde M, \tilde{\mathcal{G}})$ cannot totally reflect the geometric properties of internal space of $(M,\mathcal{F})$ and $(M,\mathcal{G})$. 

Concretely speaking, in the previous section, we discuss the affine Dirac equation $\gamma^\mu \tilde\rho_{\omega\nu;\mu} = \tilde m_{\omega\nu\tau}$ on $(\tilde M, \tilde{\mathcal{G}})$. Similar to section \ref{On-shell evolution of potential field and affine connection representation of Yang-Mills equation}, we have the affine Yang-Mills equation $\left.{\tilde K_{\nu\rho\sigma}^\mu}\right.^{:\rho}=\tilde\rho_{\nu}^\mu \gamma_\sigma$ on $(\tilde M, \tilde{\mathcal{F}})$. Suppose there is no gravitational field, then the remaining non-vanishing equations are just only 
$$\gamma^\mu \tilde\rho_{00;\mu} = \tilde m_{00\tau},\ \ \ \ \ \ \ \ \left.{\tilde K_{0\rho\sigma}^0}\right.^{:\rho}=\tilde\rho_{0}^0 \gamma_\sigma \ .$$
There are multiple internal charges 
$$\rho_{mn} \ (m,n=4,5,\cdots,\mathfrak{D})$$ 
on $(M,\mathcal{F})$. We intend to use these $\rho_{mn}$ to describe leptons and hadrons. However, via encapsulation of classical spacetime, $(\tilde M,\tilde{\mathcal{F}})$ remains only one internal charge $\tilde\rho_{00}$, and it falls short. It is impossible for the only one real-valued field function $\tilde\rho_{00}$ to describe so many leptons and hadrons. 

On the premise of not abandoning the $(1+3)$-dimensional spacetime, if we want to describe gauge fields, there is a method that to use some non-coordinate abstract degrees of freedom on the phase of $e^{iT_a \theta^a}$ of a complex-valued field function $\psi$. This way is effective, but not natural. It is not satisfactory for a theory to adopt a coordinate representation for external space but a non-coordinate representation for internal space. 

A logically more natural way is required to abandon the framework of $(1+3)$-dimensional spacetime $(\tilde M,\tilde{\mathcal{F}})$ and $(\tilde M, \tilde{\mathcal{G}})$. We should put internal space and external space together to describe their unified geometry with the same spatial frame. On $(M,\mathcal{F})$ and $(M,\mathcal{G})$, there are enough real-valued field functions $\rho_{mn}$ to describe leptons and hadrons and enough internal components $\left[mnP\right]$ of affine connection to describe gauge potentials.

Therefore, only on the full-dimensional $(M,\mathcal{F})$ and $(M,\mathcal{G})$ can total advantages of affine connection representation of gauge fields be brought into full play and thereby show complete details of geometric properties of gauge field. So we are going to stop the discussions about the classical spacetime $\tilde M$, but to focus on the full-dimensional manifold $M$. 

\sDiscussion On $M$, due to $\Gamma_{MNP} = \frac{1}{2}\left(\left[MNP\right] + \left\{MNP\right\}\right)$, $\left[MNP\right]=\delta_{AD}B_M^D \left(\frac{\partial B_N^A}{\partial x^P} + \left(^A_{BP}\right) B^B_N \right)$ and $ G_{MN} = \delta_{AB} B^A_M  B^B_N$, we know that gauge field and gravitational field can both be described by spatial frames $B^A_M$ and $C^M_A$ in a reference-system. Reference-system is the common origination of gauge field and gravitational field. The invariance under reference-system transformation is the common origination of gauge covariance and general covariance. 

We adopt the components $\left[mnP\right]$ of $\left[MNP\right]$ with $m,n\in \{4,5,\cdots,\mathfrak{D}\}$ to describe the gauge potentials of typical gauge fields such as electromagnetic, weak, and strong interaction fields and adopt the components $\rho_{mn}$ of $\rho_{MN}$ with $m,n\in \{4,5,\cdots,\mathfrak{D}\}$ to describe the charges of leptons and hadrons. The physical meanings of the other components of $\rho_{MN}$ and $\left[MNP\right]$ are not clear at present; maybe they could be used to describe dark matters and their interactions. 

On orthogonal $(M,\mathcal{G})$ and $(M,\mathcal{F})$, there are full-dimensional field equations, i.e., affine Dirac equation and affine Yang-Mills equation
\begin{equation}
\gamma^P \rho_{MN;P} =\rho_{MN;0},\ \ \ \ \ \ \ \ {K_{NPQ}^M}^{:P}=\rho_N^M \gamma_Q,
\end{equation}
which reflect the on-shell evolution directions $\nabla \rho$ and $\nabla t$, respectively. Their quantum evolutions are described by the propagators in Definition \ref{DefinitionOfPropagator} or Discussion \ref{QuantizationProcedure}.

\sDiscussion\label{PropositionOfGaugeTransformation} On an orthogonal $(M,\mathcal{G})$, Eq.(\ref{EqTotalAction}) presents as a full-dimensional action 
\begin{equation}\label{GeneralActionOfSummationOfRho}
s_{\rho}(L) = \int_L \sum_{M,N}{\left(\gamma^P \rho_{MN;P} + \varepsilon^P_0 \rho_{MN;P} \right)} d x^0 = -i\int_L \sum_{M,N}\bar\psi_{MN} \left(\gamma^P D_{MNP} + \varepsilon^P_0 D_{MNP} \right) \psi_{MN} d x^0.
\end{equation}
If and only if $L_k: g \to g'$ is an orthogonal transformation, $L_k$ sends $s_{\rho}(L)$ to
$$s'_{\rho}(L) = \int_L \sum_{M,N}{\left(\gamma^{P'} \rho_{MN;P'} + \varepsilon^{P'}_{0'} \rho_{MN;P'} \right)} d x^{0'} = -i\int_L \sum_{M,N}\bar\psi'_{MN} \left(\gamma^{P'} D'_{MNP'} + \varepsilon^{P'}_{0'} D'_{MNP'} \right) \psi'_{MN} d x^{0'},$$
where $\rho_{MN}$ is determined by the reference-system $\mathfrak{f}\circ f$ but not $\mathfrak{g}\circ g$, so $\rho_{MN}$ does not vary with the transformation $L_k: g \to g'$. We see that in affine connection representation of gauge fields, the gauge transformations $\psi \mapsto \psi'$ and $D \mapsto D'$ essentially boil down to the reference-system transformation $L_k$.

\noindent {\bf Remark 1.} For a general $(M,\mathcal{G})$, $\mathcal{G}$ is not necessarily orthogonal, so the corresponding action should be described by 
$$s_{MN}(L) = \int_L \left( B^0_0 \gamma^P \rho_{MN;P} + \varepsilon^P_0 \rho_{MN;P} \right) dx^0.$$
In this general case, Definition \ref{AffineAction} and the method in Discussion \ref{QuantizationProcedure} are also available and effective, where we take 
$$s_{MN}(L) = \int_L D\rho_{MN}.$$

\noindent {\bf Remark 2.} We see that the real-valued representation of action is more concise than the complex-valued representation of action. Hence, it is more convenient to adopt real-valued representations for field function, field equation, and action. 

In the following sections, we are going to use $\left[MNP\right]$ to show the affine connection representations of electromagnetic, weak, and strong interaction fields and to adopt the real-valued representation $\rho_{MN;P}$ to discuss the interactions between gauge fields and elementary particles. They are based on the following definition. 

\sDefinition\label{DefinitionOfTypicalGaugeField} Let $M=P \times N$, $r\DenoteBy dim P=3$ and $\mathfrak{D}\DenoteBy dim M=5\ \mbox{or}\ 6\ \mbox{or}\ 8$. Consider $\mathcal{F}=\mathfrak{f}\circ f$ and $\mathcal{G}=\mathfrak{g}\circ g$ that are defined by Eq.(\ref{DefinitionFormulaOfEvolution}), that is, $\forall p \in M$,
$$(U,\alpha^{A'})\xrightarrow{\mathfrak{f}(p)}(U,\xi^A)\xrightarrow{f(p)}(U,x^M)\xleftarrow{g(p)}(U,\zeta^A)\xleftarrow{\mathfrak{g}(p)}(U,\beta^{A'})$$
and furthermore, let
\begin{equation}\label{CoordinateRepresentationOfTypicalGaugeReferenceSysmtem}\begin{aligned}
&f(p):\ \ \xi^a =\xi^a(x^m),\ \ \ \ \xi^s =\delta_i^s x^i \ ; \ \ \ \ \ \ \ \ \ \mathfrak{f}(p):\ \ \alpha^{a'} =\alpha^{a'}(\xi^a),\ \ \ \ \alpha^{s'} =\delta_s^{s'} \xi^s \ ; \\
&g(p):\ \ \zeta^a =\zeta^a(x^m),\ \ \ \ \zeta^s =\delta_i^s x^i \ ; \ \ \ \ \ \ \ \ \ \mathfrak{g}(p):\ \ \beta^{a'} =\beta^{a'}(\zeta^a),\ \ \ \ \beta^{s'} =\delta_s^{s'} \zeta^s \ ; \\
\end{aligned}\end{equation}
($s',s,i=1,2,3; \ \ a',a,m,n=4,5,\cdots,\mathfrak{D}$) and both of $\mathcal{F}$ and $\mathcal{G}$ satisfy 
\begin{equation} \label{InternalStandardConditions}
\mbox{(i) }G_{mn}=const,\ \ \ \ \mbox{(ii) when }m\neq n, G_{mn}=0.
\end{equation}
In the above extremely simplified case, we use $\mathcal{F}$ and $\mathcal{G}$ to show electromagnetic, weak, and strong interactions without gravitation. 

\SectionWithLabel{Affine connection representation of the gauge field of weak-electromagnetic interaction}

\sDefinition\label{DefinitionOfWeakElectromagneticReferenceSystem} Suppose $(M,\mathcal{F})$ and $(M,\mathcal{G})$ conform to Definition \ref{DefinitionOfTypicalGaugeField}. Let $\mathfrak{D}=r+2=5$ and both of $\mathcal{F}$ and $\mathcal{G}$ satisfy
$$G^{(\mathfrak{D}-1)(\mathfrak{D}-1)}=G^{\mathfrak{D}\mathfrak{D}}.$$ 
Thus, $\mathcal{F}$ and $\mathcal{G}$ can describe weak and electromagnetic interactions. 

\sProposition Let the holonomic connection of $(M,\mathcal{F})$ be $\Gamma_{NP}^M$ and $\Gamma_{MNP}$. And let the coefficients of curvature tensor of $(M,\mathcal{F})$ be $K_{NPQ}^M$ and $K_{MNPQ}$. Denote
$$\begin{aligned}
&\LeftGathered
  B_P \DenoteBy\frac{1}{{\sqrt 2}}\left(\Gamma_{\mathfrak{D}\mathfrak{D}P} +\Gamma_{(\mathfrak{D}-1)(\mathfrak{D}-1)P}\right),\hfill\\
  A^3_P \DenoteBy\frac{1}{{\sqrt 2}}\left(\Gamma_{\mathfrak{D}\mathfrak{D}P} -\Gamma_{(\mathfrak{D}-1)(\mathfrak{D}-1)P}\right),\hfill\\
\RightGathered\ \ \ \ \ \ \ \ \ \ \ \ 
\LeftGathered
  B_{PQ}\DenoteBy\frac{1}{{\sqrt 2}}\left(K_{\mathfrak{D}\mathfrak{D}PQ} + K_{(\mathfrak{D}-1)(\mathfrak{D}-1)PQ}\right),\hfill\\
  F^3_{PQ}\DenoteBy\frac{1}{{\sqrt 2}}\left(K_{\mathfrak{D}\mathfrak{D}PQ} - K_{(\mathfrak{D}-1)(\mathfrak{D}-1)PQ}\right),\hfill\\ 
\RightGathered\ \ \ \\
&\LeftGathered
  A^1_P \DenoteBy\frac{1}{{\sqrt 2}}\left(\Gamma_{(\mathfrak{D}-1)\mathfrak{D}P} +\Gamma_{\mathfrak{D}(\mathfrak{D}-1)P}\right),\hfill\\
  A^2_P \DenoteBy\frac{1}{{\sqrt 2}}\left(\Gamma_{(\mathfrak{D}-1)\mathfrak{D}P} -\Gamma_{\mathfrak{D}(\mathfrak{D}-1)P}\right),\hfill\\ 
\RightGathered\ \ \ \ \ \ \ \ \ \ \ \ 
\LeftGathered
  F^1_{PQ}\DenoteBy\frac{1}{{\sqrt 2}}\left(K_{(\mathfrak{D}-1)\mathfrak{D}PQ} + K_{\mathfrak{D}(\mathfrak{D}-1)PQ}\right),\hfill\\
  F^2_{PQ}\DenoteBy\frac{1}{{\sqrt 2}}\left(K_{(\mathfrak{D}-1)\mathfrak{D}PQ} - K_{\mathfrak{D}(\mathfrak{D}-1)PQ}\right).\hfill\\ 
\RightGathered
\end{aligned}$$
And denote $g\DenoteBy \sqrt{\left(G^{(\mathfrak{D}-1)(\mathfrak{D}-1)}\right)^2+\left(G^{\mathfrak{D}\mathfrak{D}}\right)^2}$. Thus, the following equations hold spontaneously.
$$\LeftList
B_{PQ}=\frac{{\partial B_Q}}{{\partial x^P}}-\frac{{\partial B_P}}{{\partial x^Q}}, \hfill\\
F^3_{PQ}=\frac{{\partial A^3_Q}}{{\partial x^P}}-\frac{{\partial A^3_P}}{{\partial x^Q}}+g\left(A^1_P A^2_Q  - A^2_P A^1_Q\right), \hfill\\
F^1_{PQ}=\frac{{\partial A^1_Q}}{{\partial x^P}}-\frac{{\partial A^1_P}}{{\partial x^Q}}+g\left(A^2_P A^3_Q  - A^3_P A^2_Q\right), \hfill\\
F^2_{PQ}=\frac{{\partial A^2_Q}}{{\partial x^P}}-\frac{{\partial A^2_P}}{{\partial x^Q}}+g\left(A^1_P A^3_Q - A^3_P A^1_Q\right). \hfill\\
\RightList$$

\sProof Due to Eq.(\ref{CoordinateRepresentationOfTypicalGaugeReferenceSysmtem}) it is obtained that the semi-metric of $(M,\mathfrak{f})$ satisfies
$$(B_\mathfrak{f})_a^{s'} =0,\ \ (C_\mathfrak{f})_{s'}^a =0,\ \ (B_\mathfrak{f})_s^{a'} =0,\ \ (C_\mathfrak{f})_{a'}^s =0,\ \ (B_\mathfrak{f})_s^{s'} =\delta_s^{s'},\ \ (C_\mathfrak{f})_{s'}^s =\delta_{s'}^s.$$
Then, compute $(^A_{BC})_\mathfrak{f} \DenoteBy \frac{1}{2} (C_\mathfrak{f})_{A'}^A \left(\frac{\partial (B_\mathfrak{f})_B^{A'}}{\partial \xi^C} + \frac{\partial (B_\mathfrak{f})_C^{A'}}{\partial \xi^B}\right)$ and we obtain
$$\left(^s_{BC}\right)_\mathfrak{f}=0,\ \ \ \ \ \left(^a_{tu}\right)_\mathfrak{f} =0,\ \ \ \ \ \left(^a_{bC}\right)_\mathfrak{f} \neq 0; \ \ \ \ \ s,t,u=1,2,3;\ \ \ \ a,b=4,5,\cdots,\mathfrak{D};\ \ \ \ A,B,C=1,2,\cdots,\mathfrak{D}\ .$$
It is obtained from Eq.(\ref{CoordinateRepresentationOfTypicalGaugeReferenceSysmtem}) again that the semi-metric of $(M,f)$ satisfies 
$$B_m^s =0,\ \ C_s^m =0,\ \ B_i^a =0,\ \ C_a^i =0,\ \ B_i^s =\delta_i^s,\ \ C_s^i =\delta_s^i.$$
Let $s',t',i,j,k=1,2,3;\ a',b',m,n,p=4,5,\cdots,\mathfrak{D}$. Compute the metric of $(M,\mathcal{F})$ and we obtain
$$\LeftList
\LeftGathered
  G_{ij}=\delta_{s't'}B_i^{s'} B_j^{t'}  +\delta_{a'b'}B_i^{a'} B_j^{b'} =\delta_{s't'}\delta_i^{s'}\delta_j^{t'} =\delta_{ij},\hfill\\
  G_{in}=\delta_{s't'}B_i^{s'} B_n^{t'}  +\delta_{a'b'}B_i^{a'} B_n^{b'} =0,\hfill\\
  G_{mj}=\delta_{s't'}B_m^{s'} B_j^{t'}  +\delta_{a'b'}B_m^{a'} B_j^{b'} =0,\hfill\\
  G_{mn}=B_m^{\mathfrak{D}-1}B_n^{\mathfrak{D}-1} + B_m^\mathfrak{D}B_n^\mathfrak{D} = const,\hfill\\ 
\RightGathered\ \ \LeftGathered
  G^{ij}=\delta^{s't'}C_{s'}^i C_{t'}^j =\delta^{s't'}\delta_{s'}^i\delta_{t'}^j =\delta^{ij},\hfill\\
  G^{in}=\delta^{s't'}C_{s'}^i C_{t'}^n =0,\hfill\\
  G^{mj}=\delta^{s't'}C_{s'}^m C_{t'}^j =0,\hfill\\
  G^{mn}=C_{\mathfrak{D}-1}^m C_{\mathfrak{D}-1}^n  + C_\mathfrak{D}^m C_\mathfrak{D}^n = const. \hfill\\ 
\RightGathered\hfill\\
\RightList$$
Compute the holonomic connection of $\mathcal{F}$ according to $\Gamma_{NP}^M \DenoteBy\frac{1}{2}\left([^M_{NP}] + \{^M_{NP}\}\right) = \frac{1}{2}\left(C_A^M \frac{\partial B_N^A}{\partial x^P} + C_A^M \left(^A_{BP}\right)_\mathfrak{f} B_N^B\right)$, and it is obtained that 
\begin{equation}\label{CalculationResultofWeakElectroPotentials}\LeftGathered
 \Gamma_{NP}^i =0,\hfill\\
 \Gamma_{jk}^m =0,\hfill\\
 \Gamma_{nP}^m =\frac{1}{2}\left(C_a^m\frac{{\partial B_n^a}}{{\partial x^P}}+ C_a^m \left(^a_{bP}\right)_\mathfrak{f} B_n^b\right),\hfill\\
 \Gamma_{Np}^m =\frac{1}{2}\left(C_a^m\frac{{\partial B_N^a}}{{\partial x^p}}+ C_a^m \left(^a_{Bp}\right)_\mathfrak{f} B_N^B\right),\hfill\\ 
\RightGathered\ \ \ \ \ \ \ \ \ \ \ \ \ \ \ \ \ \ \ 
\LeftGathered
 \Gamma_{iNP} =G_{iM'}\Gamma_{NP}^{M'}=G_{ii'}\Gamma_{NP}^{i'}=0,\hfill\\
 \Gamma_{mjk} =G_{mM'}\Gamma_{jk}^{M'}=G_{mm'}\Gamma_{jk}^{m'}=0,\hfill\\
 \Gamma_{mnP} =\frac{1}{2}\delta_{ab}B_m^b\left(\frac{{\partial B_n^a}}{{\partial x^P}}+ \left(^a_{bP}\right)_\mathfrak{f} B_n^b\right),\hfill\\
 \Gamma_{mNp} =\frac{1}{2}\delta_{ab}B_m^b\left(\frac{{\partial B_N^a}}{{\partial x^p}}+ \left(^a_{Bp}\right)_\mathfrak{f} B_N^B\right).\hfill\\ 
\RightGathered\end{equation}
Compute the coefficients of curvature of $\mathcal{F}$, that is, 
$$\LeftList
K_{nPQ}^m \DenoteBy\frac{{\partial\Gamma_{nQ}^m}}{{\partial x^P}}-\frac{{\partial\Gamma_{nP}^m}}{{\partial x^Q}}+\Gamma_{HP}^m\Gamma_{nQ}^H  -\Gamma_{nP}^H\Gamma_{HQ}^m,\ \ \ \ \ \ \ \ K_{mnPQ}\DenoteBy G_{mM'}K_{nPQ}^{M'} = G_{mm'}K_{nPQ}^{m'},\hfill\\
\RightList$$
then we obtain
$$\LeftList
  K_{(\mathfrak{D}-1)(\mathfrak{D}-1)PQ}=\frac{{\partial\Gamma_{(\mathfrak{D}-1)(\mathfrak{D}-1)Q}}}{{\partial x^P}}-\frac{{\partial\Gamma_{(\mathfrak{D}-1)(\mathfrak{D}-1)P}}}{{\partial x^Q}}+G^{\mathfrak{D}\mathfrak{D}}\left(\Gamma_{(\mathfrak{D}-1)\mathfrak{D}P}\Gamma_{\mathfrak{D}(\mathfrak{D}-1)Q}-\Gamma_{\mathfrak{D}(\mathfrak{D}-1)P}\Gamma_{(\mathfrak{D}-1)\mathfrak{D}Q}\right),\hfill\\
  K_{\mathfrak{D}(\mathfrak{D}-1)PQ}=\frac{{\partial\Gamma_{\mathfrak{D}(\mathfrak{D}-1)Q}}}{{\partial x^P}}-\frac{{\partial\Gamma_{\mathfrak{D}(\mathfrak{D}-1)P}}}{{\partial x^Q}}+G^{\mathfrak{D}\mathfrak{D}}\left(\Gamma_{\mathfrak{D}\mathfrak{D}P}\Gamma_{\mathfrak{D}(\mathfrak{D}-1)Q}-\Gamma_{\mathfrak{D}(\mathfrak{D}-1)P}\Gamma_{\mathfrak{D}\mathfrak{D}Q}\right)\hfill\\
\ \ \ \ \ \ \ \ \ \ \ \ \ \ \ \ \ \ \ \ \ \ +G^{(\mathfrak{D}-1)(\mathfrak{D}-1)}\left(\Gamma_{\mathfrak{D}(\mathfrak{D}-1)P}\Gamma_{(\mathfrak{D}-1)(\mathfrak{D}-1)Q}-\Gamma_{(\mathfrak{D}-1)(\mathfrak{D}-1)P}\Gamma_{\mathfrak{D}(\mathfrak{D}-1)Q}\right),\hfill\\
  K_{(\mathfrak{D}-1)\mathfrak{D}PQ}=\frac{{\partial\Gamma_{(\mathfrak{D}-1)\mathfrak{D}Q}}}{{\partial x^P}}-\frac{{\partial\Gamma_{(\mathfrak{D}-1)\mathfrak{D}P}}}{{\partial x^Q}}+G^{\mathfrak{D}\mathfrak{D}}\left(\Gamma_{(\mathfrak{D}-1)\mathfrak{D}P}\Gamma_{\mathfrak{D}\mathfrak{D}Q}-\Gamma_{\mathfrak{D}\mathfrak{D}P}\Gamma_{(\mathfrak{D}-1)\mathfrak{D}Q}\right)\hfill\\
  \ \ \ \ \ \ \ \ \ \ \ \ \ \ \ \ \ \ \ \ \ \ +G^{(\mathfrak{D}-1)(\mathfrak{D}-1)}\left(\Gamma_{(\mathfrak{D}-1)(\mathfrak{D}-1)P}\Gamma_{(\mathfrak{D}-1)\mathfrak{D}Q}-\Gamma_{(\mathfrak{D}-1)\mathfrak{D}P}\Gamma_{(\mathfrak{D}-1)(\mathfrak{D}-1)Q}\right),\hfill\\
  K_{\mathfrak{D}\mathfrak{D}PQ}=\frac{{\partial\Gamma_{\mathfrak{D}\mathfrak{D}Q}}}{{\partial x^P}}-\frac{{\partial\Gamma_{\mathfrak{D}\mathfrak{D}P}}}{{\partial x^Q}}+G^{(\mathfrak{D}-1)(\mathfrak{D}-1)}\left(\Gamma_{\mathfrak{D}(\mathfrak{D}-1)P}\Gamma_{(\mathfrak{D}-1)\mathfrak{D}Q}-\Gamma_{(\mathfrak{D}-1)\mathfrak{D}P}\Gamma_{\mathfrak{D}(\mathfrak{D}-1)Q}\right).\hfill\\
\RightList$$
Hence,
$$\begin{aligned}
B_{PQ}&\DenoteBy\frac{1}{{\sqrt 2}}\left(K_{\mathfrak{D}\mathfrak{D}PQ} + K_{(\mathfrak{D}-1)(\mathfrak{D}-1)PQ}\right)\\
&=\frac{1}{{\sqrt 2}}\frac{{\partial\left(\Gamma_{\mathfrak{D}\mathfrak{D}Q} +\Gamma_{(\mathfrak{D}-1)(\mathfrak{D}-1)Q}\right)}}{{\partial x^P}}-\frac{1}{{\sqrt 2}}\frac{{\partial\left(\Gamma_{\mathfrak{D}\mathfrak{D}P} +\Gamma_{(\mathfrak{D}-1)(\mathfrak{D}-1)P}\right)}}{{\partial x^Q}}
=\frac{{\partial B_Q}}{{\partial x^P}}-\frac{{\partial B_P}}{{\partial x^Q}}. \\
F^3_{PQ}&\DenoteBy\frac{1}{{\sqrt 2}}\left(K_{\mathfrak{D}\mathfrak{D}PQ} - K_{(\mathfrak{D}-1)(\mathfrak{D}-1)PQ}\right)\\
&=\frac{1}{{\sqrt 2}}\left(\frac{{\partial\Gamma_{\mathfrak{D}\mathfrak{D}Q}}}{{\partial x^P}}-\frac{{\partial\Gamma_{\mathfrak{D}\mathfrak{D}P}}}{{\partial x^Q}}+G^{(\mathfrak{D}-1)(\mathfrak{D}-1)}\left(\Gamma_{\mathfrak{D}(\mathfrak{D}-1)P}\Gamma_{(\mathfrak{D}-1)\mathfrak{D}Q}-\Gamma_{(\mathfrak{D}-1)\mathfrak{D}P}\Gamma_{\mathfrak{D}(\mathfrak{D}-1)Q}\right)\right)\\
&-\frac{1}{{\sqrt 2}}\left(\frac{{\partial\Gamma_{(\mathfrak{D}-1)(\mathfrak{D}-1)Q}}}{{\partial x^P}}-\frac{{\partial\Gamma_{(\mathfrak{D}-1)(\mathfrak{D}-1)P}}}{{\partial x^Q}}+G^{\mathfrak{D}\mathfrak{D}}\left(\Gamma_{(\mathfrak{D}-1)\mathfrak{D}P}\Gamma_{\mathfrak{D}(\mathfrak{D}-1)Q}-\Gamma_{\mathfrak{D}(\mathfrak{D}-1)P}\Gamma_{(\mathfrak{D}-1)\mathfrak{D}Q}\right)\right)\\
&=\frac{{\partial A^3_Q}}{{\partial x^P}}-\frac{{\partial A^3_P}}{{\partial x^Q}}+g\left(\Gamma_{\mathfrak{D}(\mathfrak{D}-1)P}\Gamma_{(\mathfrak{D}-1)\mathfrak{D}Q} -\Gamma_{(\mathfrak{D}-1)\mathfrak{D}P}\Gamma_{\mathfrak{D}(\mathfrak{D}-1)Q}\right)\\
&=\frac{{\partial A^3_Q}}{{\partial x^P}}-\frac{{\partial A^3_P}}{{\partial x^Q}}+g\left(A^1_P A^2_Q  - A^2_P A^1_Q\right).\\
\end{aligned}$$
Then, $F^1_{PQ}$ and $F^2_{PQ}$ can also be computed similarly.
\qed

\sRemark\label{RemarkGWS} Comparing the above conclusion and $U(1) \times SU(2)$ principal bundle theory, we know this proposition shows that the reference-system $\mathcal{F}$ indeed can describe weak and electromagnetic field.

The following proposition shows an advantage of affine connection representation, that is, affine connection representation spontaneously implies the chiral asymmetry of neutrinos, but $U(1) \times SU(2)$ principal bundle connection representation cannot imply it spontaneously. 

\sDefinition\label{DefinitionOfLeptons} According to Definition \ref{DefinitionOfGeneralCharge}, let the charges of the above reference-system $\mathcal{F}$ be $\rho_{mn}$, where $m,n\in \{\mathfrak{D}-1, \mathfrak{D}\} = \{4,5\}$. Then, $l\DenoteBy(\rho_{(\mathfrak{D}-1)(\mathfrak{D}-1)},\ \rho_{\mathfrak{D}\mathfrak{D}})^T$ is said to be an {\bf electric charged lepton}, $\nu \DenoteBy(\rho_{\mathfrak{D}(\mathfrak{D}-1)},\ \rho_{(\mathfrak{D}-1)\mathfrak{D}})^T$ is said to be a {\bf neutrino}. $l$ and $\nu$ are collectively denoted by $L$. Thus, $\frac{1}{\sqrt 2}\left(1, 1\right)L$ is said to be a {\bf left-handed lepton}, $\frac{1}{\sqrt 2}\left(1, -1\right)L$ is said to be a {\bf right-handed lepton}, denoted by
\begin{equation}\label{LeptonsDefinition}\LeftGathered
l_L \DenoteBy\frac{1}{{\sqrt 2}}\left(\rho_{\left(\mathfrak{D}-1\right)\left(\mathfrak{D}-1\right)} +\rho_{\mathfrak{D}\mathfrak{D}}\right),\hfill\\
l_R \DenoteBy\frac{1}{{\sqrt 2}}\left(\rho_{\left(\mathfrak{D}-1\right)\left(\mathfrak{D}-1\right)} -\rho_{\mathfrak{D}\mathfrak{D}}\right),\hfill\\
\RightGathered\ \ \ 
\LeftGathered
\nu_L \DenoteBy\frac{1}{{\sqrt 2}}\left(\rho_{\mathfrak{D}\left(\mathfrak{D}-1\right)} +\rho_{\left(\mathfrak{D}-1\right)\mathfrak{D}}\right),\hfill\\
\nu_R \DenoteBy\frac{1}{{\sqrt 2}}\left(\rho_{\mathfrak{D}\left(\mathfrak{D}-1\right)} -\rho_{\left(\mathfrak{D}-1\right)\mathfrak{D}}\right).\hfill\\
\RightGathered\end{equation}
Denote $(\Gamma_\mathcal{G})_{MNP}$ by $\Gamma_{MNP}$ concisely. Then, we define on $(M,\mathcal{G})$ that
\begin{equation}\label{WeakAndElectromagneticPotentialsDefinition}\LeftGathered
  Z_P \DenoteBy\frac{1}{{\sqrt 2}}\left(\Gamma_{(\mathfrak{D}-1)(\mathfrak{D}-1)P} +\Gamma_{\mathfrak{D}\mathfrak{D}P}\right),\hfill\\
  A_P \DenoteBy\frac{1}{{\sqrt 2}}\left(\Gamma_{(\mathfrak{D}-1)(\mathfrak{D}-1)P} -\Gamma_{\mathfrak{D}\mathfrak{D}P}\right),\hfill\\ 
\RightGathered\ \ \ 
\LeftGathered
  W^1_P \DenoteBy\frac{1}{{\sqrt 2}}\left(\Gamma_{(\mathfrak{D}-1)\mathfrak{D}P} +\Gamma_{\mathfrak{D}(\mathfrak{D}-1)P}\right),\hfill\\
  W^2_P \DenoteBy\frac{1}{{\sqrt 2}}\left(\Gamma_{(\mathfrak{D}-1)\mathfrak{D}P} -\Gamma_{\mathfrak{D}(\mathfrak{D}-1)P}\right),\hfill\\ 
\RightGathered\end{equation}
and say $A_P$ is {\bf (affine) electromagnetic potential}, while $Z_P$, $W^1_P$, and $W^2_P$ are {\bf (affine) weak gauge potentials}.

\sProposition\label{PropositionOfLeptonEvolutionForm} If $(M,\mathcal{G})$ satisfies the symmetry condition $\Gamma_{(\mathfrak{D}-1)\mathfrak{D}P}=\Gamma_{\mathfrak{D}(\mathfrak{D}-1)P}$, then the geometric properties $l$ and $\nu$ of $\mathcal{F}$ satisfy the following conclusions on $(M,\mathcal{G})$,
\begin{equation}\label{LeptonsEvolutionSimple}\LeftList
\LeftGathered
  l_{L;P}=\partial_P l_L  -g l_L Z_P  -g l_R A_P  -g\nu_L W^1_P \ , \hfill\\
  l_{R;P}=\partial_P l_R  -g l_R Z_P  -g l_L A_P \ , \hfill\\ 
  \nu_{L;P}=\partial_P\nu_L  -g\nu_L Z_P  -g l_L W^1_P \ , \hfill\\
  \nu_{R;P}=\partial_P\nu_R  -g\nu_R Z_P \ . \hfill\\ 
\RightGathered\\
\RightList\end{equation}

\sProof Let $H\in \{1,2,3,4,5\}$, $h\in \{4,5\}$. It follows from Eq.(\ref{CalculationResultofWeakElectroPotentials}) that
$$\begin{aligned}
\rho_{mn;P}&=\partial_P\rho_{mn} -\rho_{Hn}\Gamma_{mP}^H  -\rho_{mH}\Gamma_{nP}^H \\
&=\partial_P\rho_{mn} -\rho_{hn}\Gamma_{mP}^h  -\rho_{mh}\Gamma_{nP}^h.\\
\end{aligned}$$
Then, Eq.(\ref{LeptonsDefinition}) and Eq.(\ref{WeakAndElectromagneticPotentialsDefinition}) lead to Eq.(\ref{LeptonsEvolutionSimple}).
\qed

\newpage

\sRemark\label{DiscussionWaysOfGeometryAndGroup} From the above proposition, we see that some constraint conditions make the general linear group $GL(2,\mathbb{R})$ broken to a smaller group $S$, i.e.,
$$GL(2,\mathbb{R})
 \xrightarrow{\ \ \ \ G_{(\mathfrak{D}-1)(\mathfrak{D}-1)}=G_{\mathfrak{D}\mathfrak{D}},\ \ \ \ \Gamma_{(\mathfrak{D}-1)\mathfrak{D}P}=\Gamma_{\mathfrak{D}(\mathfrak{D}-1)P} \ \ \ \ }
S,$$
so that the chiral asymmetry of leptons arises in Eq.(\ref{LeptonsEvolutionSimple}) spontaneously.

\sRemark \label{HiggsExplanation} Proposition \ref{PropositionOfLeptonEvolutionForm} shows that:

(1) In affine connection representation of gauge fields, the coupling constant $g$ is possessed of a geometric meaning, that is in fact the metric of internal space. But it does not have such a clear geometric meaning in $U(1) \times SU(2)$ principal bundle connection representation. 

(2) At the most fundamental level, the coupling constant of $Z_P$ and that of $A_P$ are equal, i.e., 
$$g_Z=g_A=g.$$ 
Suppose there is a kind of medium. $Z$ boson and photon move in it. Suppose $Z$ field has interaction with the medium, but electromagnetic field $A$ has no interaction with the medium. Thus, we have coupling constants 
$$\tilde g_Z \neq g_A = g$$
in the medium, and the Weinberg angle arises. 

It is quite reasonable to consider a Higgs boson as a zero-spin pair of neutrinos, because in the Lagrangian, Higgs boson only couples with $Z$ field and $W$ field but does not couple with electromagnetic field and gluon field. If so, Higgs boson would lose its fundamentality and it would not have enough importance in a theory at the most fundamental level. 

(3) The mixing of three generations of leptons does not appear in Proposition \ref{PropositionOfLeptonEvolutionForm}, but it can spontaneously arise in Proposition \ref{PropositionOfMNSMixing} due to the affine connection representation of the gauge field that is given by Definition \ref{DefinitionOfElectroWeakStrongUnifiedReferenceSystem}. 

\SectionWithLabel{Affine connection representation of the gauge field of strong interaction}

\sDefinition\label{DefinitionOfTypicalStrongReferenceSystem} Suppose $(M,\mathcal{F})$ and $(M,\mathcal{G})$ conform to Definition \ref{DefinitionOfTypicalGaugeField}. Let $\mathfrak{D}=r+3=6$ and both of $\mathcal{F}$ and $\mathcal{G}$ satisfy
$$G^{(\mathfrak{D}-2)(\mathfrak{D}-2)}=G^{(\mathfrak{D}-1)(\mathfrak{D}-1)}=G^{\mathfrak{D}\mathfrak{D}}.$$ 
Thus, $\mathcal{F}$ and $\mathcal{G}$ can describe strong interaction. 

\sDefinition\label{DefinitionOfColorCharge}\label{DefinitionOfStrongPotential} According to Definition \ref{DefinitionOfGeneralCharge}, let the charges of $\mathcal{F}$ be $\rho_{mn}$, where $m,n=4,5,\cdots,\mathfrak{D}$. Define 
$$\LeftGathered
  d_1 \DenoteBy(\rho_{(\mathfrak{D}- 2)(\mathfrak{D}- 2)},\ \ \rho_{(\mathfrak{D}-1)(\mathfrak{D}-1)})^T,\hfill\\
  d_2 \DenoteBy(\rho_{(\mathfrak{D}-1)(\mathfrak{D}-1)},\ \ \rho_{\mathfrak{D}\mathfrak{D}})^T,\hfill\\
  d_3 \DenoteBy(\rho_{\mathfrak{D}\mathfrak{D}},\ \ \rho_{(\mathfrak{D}- 2)(\mathfrak{D}- 2)})^T,\hfill\\ 
\RightGathered\ \ \ 
\LeftGathered
  u_1 \DenoteBy(\rho_{(\mathfrak{D}- 2)(\mathfrak{D}-1)},\ \ \rho_{(\mathfrak{D}-1)(\mathfrak{D}- 2)})^T,\hfill\\
  u_2 \DenoteBy(\rho_{(\mathfrak{D}-1)\mathfrak{D}},\ \ \rho_{\mathfrak{D}(\mathfrak{D}-1)})^T,\hfill\\
  u_3 \DenoteBy(\rho_{\mathfrak{D}(\mathfrak{D}- 2)},\ \ \rho_{(\mathfrak{D}- 2)\mathfrak{D}})^T.\hfill\\ 
\RightGathered$$
We say $d_1$ and $u_1$ are {\bf red color charges}, $d_2$ and $u_2$ are {\bf blue color charges}, and $d_3$ and $u_3$ are {\bf green color charges}. Then, $d_1$, $d_2$, $d_3$ are said to be {\bf down-type color charges}, and $u_1$, $u_2$, $u_3$ are said to be {\bf up-type color charges}. Their left-handed and right-handed charges are
$$\LeftGathered
  d_{1L} \DenoteBy\frac{1}{{\sqrt 2}}\left(\rho_{\left(\mathfrak{D}- 2\right)\left(\mathfrak{D}- 2\right)} +\rho_{\left(\mathfrak{D}-1\right)\left(\mathfrak{D}-1\right)}\right),\hfill\\
  d_{2L} \DenoteBy\frac{1}{{\sqrt 2}}\left(\rho_{\left(\mathfrak{D}-1\right)\left(\mathfrak{D}-1\right)} +\rho_{\mathfrak{D}\mathfrak{D}}\right),\hfill\\
  d_{3L} \DenoteBy\frac{1}{{\sqrt 2}}\left(\rho_{\mathfrak{D}\mathfrak{D}} +\rho_{\left(\mathfrak{D}- 2\right)\left(\mathfrak{D}- 2\right)}\right),\hfill\\ 
\RightGathered\ \ \ 
\LeftGathered
  d_{1R} \DenoteBy\frac{1}{{\sqrt 2}}\left(\rho_{\left(\mathfrak{D}- 2\right)\left(\mathfrak{D}- 2\right)} -\rho_{\left(\mathfrak{D}-1\right)\left(\mathfrak{D}-1\right)}\right),\hfill\\
  d_{2R} \DenoteBy\frac{1}{{\sqrt 2}}\left(\rho_{\left(\mathfrak{D}-1\right)\left(\mathfrak{D}-1\right)} -\rho_{\mathfrak{D}\mathfrak{D}}\right),\hfill\\
  d_{3R} \DenoteBy\frac{1}{{\sqrt 2}}\left(\rho_{\mathfrak{D}\mathfrak{D}} -\rho_{\left(\mathfrak{D}- 2\right)\left(\mathfrak{D}- 2\right)}\right),\hfill\\ 
\RightGathered$$
$$\LeftGathered
  u_{1L}\DenoteBy\frac{1}{{\sqrt 2}}\left(\rho_{\left(\mathfrak{D}- 2\right)\left(\mathfrak{D}-1\right)} +\rho_{\left(\mathfrak{D}-1\right)\left(\mathfrak{D}- 2\right)}\right),\hfill\\
  u_{2L}\DenoteBy\frac{1}{{\sqrt 2}}\left(\rho_{\left(\mathfrak{D}-1\right)\mathfrak{D}} +\rho_{\mathfrak{D}\left(\mathfrak{D}-1\right)}\right),\hfill\\
  u_{3L}\DenoteBy\frac{1}{{\sqrt 2}}\left(\rho_{\mathfrak{D}\left(\mathfrak{D}- 2\right)} +\rho_{\left(\mathfrak{D}- 2\right)\mathfrak{D}}\right),\hfill\\ 
\RightGathered\ \ \ 
\LeftGathered
  u_{1R}\DenoteBy\frac{1}{{\sqrt 2}}\left(\rho_{\left(\mathfrak{D}- 2\right)\left(\mathfrak{D}-1\right)} -\rho_{\left(\mathfrak{D}-1\right)\left(\mathfrak{D}- 2\right)}\right),\hfill\\
  u_{2R}\DenoteBy\frac{1}{{\sqrt 2}}\left(\rho_{\left(\mathfrak{D}-1\right)\mathfrak{D}} -\rho_{\mathfrak{D}\left(\mathfrak{D}-1\right)}\right),\hfill\\
  u_{3R}\DenoteBy\frac{1}{{\sqrt 2}}\left(\rho_{\mathfrak{D}\left(\mathfrak{D}- 2\right)} -\rho_{\left(\mathfrak{D}- 2\right)\mathfrak{D}}\right).\hfill\\ 
\RightGathered$$
On $(M,\mathcal{G})$ we denote
$$\begin{aligned}
g_s&\DenoteBy\sqrt{ \left(G^{(\mathfrak{D}-1)(\mathfrak{D}-1)}\right)^2+\left(G^{\mathfrak{D}\mathfrak{D}}\right)^2}\\
&=\sqrt{ \left(G^{(\mathfrak{D}-1)(\mathfrak{D}-1)}\right)^2+\left(G^{(\mathfrak{D}-2)(\mathfrak{D}-2)}\right)^2}\hfill\\
&=\sqrt{ \left(G^{(\mathfrak{D}-2)(\mathfrak{D}-2)}\right)^2+\left(G^{\mathfrak{D}\mathfrak{D}}\right)^2}.\hfill\\
\end{aligned}$$
$$\LeftGathered
  U^1_P \DenoteBy\frac{1}{{\sqrt 2}}\left(\Gamma_{(\mathfrak{D}- 2)(\mathfrak{D}- 2)P} +\Gamma_{(\mathfrak{D}-1)(\mathfrak{D}-1)P}\right),\hfill\\
  V^1_P \DenoteBy\frac{1}{{\sqrt 2}}\left(\Gamma_{(\mathfrak{D}- 2)(\mathfrak{D}- 2)P} -\Gamma_{(\mathfrak{D}-1)(\mathfrak{D}-1)P}\right),\hfill\\ 
\RightGathered\ \ \ 
\LeftGathered
  X^{23}_P \DenoteBy\frac{1}{{\sqrt 2}}\left(\Gamma_{(\mathfrak{D}- 2)(\mathfrak{D}-1)P} +\Gamma_{(\mathfrak{D}-1)(\mathfrak{D}- 2)P}\right),\hfill\\
  Y^{23}_P \DenoteBy\frac{1}{{\sqrt 2}}\left(\Gamma_{(\mathfrak{D}- 2)(\mathfrak{D}-1)P} -\Gamma_{(\mathfrak{D}-1)(\mathfrak{D}- 2)P}\right),\hfill\\ 
\RightGathered$$
$$\LeftGathered
  U^2_P \DenoteBy\frac{1}{{\sqrt 2}}\left(\Gamma_{(\mathfrak{D}-1)(\mathfrak{D}-1)P} +\Gamma_{\mathfrak{D}\mathfrak{D}P}\right),\hfill\\
  V^2_P \DenoteBy\frac{1}{{\sqrt 2}}\left(\Gamma_{(\mathfrak{D}-1)(\mathfrak{D}-1)P} -\Gamma_{\mathfrak{D}\mathfrak{D}P}\right),\hfill\\ 
\RightGathered\ \ \ 
\LeftGathered
  X^{31}_P \DenoteBy\frac{1}{{\sqrt 2}}\left(\Gamma_{(\mathfrak{D}-1)\mathfrak{D}P} +\Gamma_{\mathfrak{D}(\mathfrak{D}-1)P}\right),\hfill\\
  Y^{31}_P \DenoteBy\frac{1}{{\sqrt 2}}\left(\Gamma_{(\mathfrak{D}-1)\mathfrak{D}P} -\Gamma_{\mathfrak{D}(\mathfrak{D}-1)P}\right),\hfill\\ 
\RightGathered$$
$$\LeftGathered
  U^3_P \DenoteBy\frac{1}{{\sqrt 2}}\left(\Gamma_{\mathfrak{D}\mathfrak{D}P} +\Gamma_{(\mathfrak{D}- 2)(\mathfrak{D}- 2)P}\right),\hfill\\
  V^3_P \DenoteBy\frac{1}{{\sqrt 2}}\left(\Gamma_{\mathfrak{D}\mathfrak{D}P} -\Gamma_{(\mathfrak{D}- 2)(\mathfrak{D}- 2)P}\right),\hfill\\ 
\RightGathered\ \ \ 
\LeftGathered
  X^{12}_P \DenoteBy\frac{1}{{\sqrt 2}}\left(\Gamma_{\mathfrak{D}(\mathfrak{D}- 2)P} +\Gamma_{(\mathfrak{D}- 2)\mathfrak{D}P}\right),\hfill\\
  Y^{12}_P \DenoteBy\frac{1}{{\sqrt 2}}\left(\Gamma_{\mathfrak{D}(\mathfrak{D}- 2)P} -\Gamma_{(\mathfrak{D}- 2)\mathfrak{D}P}\right).\hfill\\ 
\RightGathered$$
We notice that there are just only three independent ones in $U^1_P$, $U^2_P$, $U^3_P$, $V^1_P$, $V^2_P$, and $V^3_P$. Without loss of generality, let
$$\LeftGathered
  R_P \DenoteBy a_R U^1_P  + b_R U^2_P  + c_R U^3_P, \hfill\\
  S_P \DenoteBy a_S U^1_P  + b_S U^2_P  + c_S U^3_P, \hfill\\
  T_P \DenoteBy a_T U^1_P  + b_T U^2_P  + c_T U^3_P, \hfill\\ 
\RightGathered\ \ \ 
\LeftGathered
  U^1_P \DenoteBy\alpha_R R_P  +\alpha_S S_P  +\alpha_T T_P, \hfill\\
  U^2_P \DenoteBy\beta_R R_P  +\beta_S S_P  +\beta_T T_P, \hfill\\
  U^3_P \DenoteBy\gamma_R R_P  +\gamma_S S_P  +\gamma_T T_P, \hfill\\ 
\RightGathered$$
where the coefficients matrix is non-singular. Thus, it is not hard to find the following proposition true.

\sProposition Let $\lambda_a\ (a=1,2,\cdots ,8)$ be the Gell-Mann matrices, and $T_a \DenoteBy\frac{1}{2}\lambda_a$ the generators of $SU(3)$ group. When $(M,\mathcal{G})$ satisfies the symmetry condition $\Gamma_{(\mathfrak{D}- 2)(\mathfrak{D}- 2)P} +\Gamma_{(\mathfrak{D}-1)(\mathfrak{D}-1)P} +\Gamma_{\mathfrak{D}\mathfrak{D}P}=0$, denote
$$
A_P \DenoteBy\frac{1}{2}
\left(\begin{array}{*{20}c}
  {A^{11}_P}&{A^{12}_P}&{A^{13}_P}\\
  {A^{21}_P}&{A^{22}_P}&{A^{23}_P}\\
  {A^{31}_P}&{A^{32}_P}&{A^{33}_P}\\
\end{array}\right),$$
where
$$\LeftList
  A^{11}_P \DenoteBy S_P  +\frac{1}{{\sqrt 6}}T_P,\ \ A^{12}_P \DenoteBy X^{12}_P  - iY^{12}_P,\ \ \ \ \ \ A^{13}_P \DenoteBy X^{31}_P  - iY^{31}_P, \hfill\\
  A^{21}_P \DenoteBy X^{12}_P  + iY^{12}_P,\ \ \ A^{22}_P \DenoteBy - S_P  +\frac{1}{{\sqrt 6}}T_P,\ \ A^{23}_P \DenoteBy X^{23}_P  - iY^{23}_P,  \hfill\\
  A^{31}_P \DenoteBy X^{31}_P  + iY^{31}_P,\ \ \ A^{32}_P \DenoteBy X^{23}_P  + iY^{23}_P,\ \ \ \ \ \ \ A^{33}_P \DenoteBy  -\frac{2}{{\sqrt 6}}T_P. \hfill\\
\RightList
$$
Thus, $A_P=T_a A^a_P$ if and only if 
$$\begin{aligned}
&A^1_P = X^{12}_P,\ \ \ \ A^2_P = Y^{12}_P,\ \ \ \ A^3_P = S_P,\ \ \ \ A^4_P = X^{31}_P, \\
&A^5_P = Y^{31}_P,\ \ \ \ A^6_P = X^{23}_P,\ \ \ \ A^7_P = Y^{23}_P,\ \ \ \ A^8_P = T_P. \\
\end{aligned}$$

\sRemark On the one hand, the above proposition shows that Definition \ref{DefinitionOfTypicalStrongReferenceSystem} is an affine connection representation of strong interaction field. It does not define the gauge potentials as abstractly as that in principal $SU(3)$-bundle theory but endows gauge potentials with concrete geometric constructions.

On the other hand, the above proposition implies that if we take appropriate symmetry conditions, the algebraic properties of $SU(3)$ group can be described by the transformation group $GL(3,\mathbb R)$ of internal space of $\mathcal{G}$. In other words, the exponential map 
$$exp: \ GL(3,\mathbb R)\to U(3),\ \ [B_m^a]\mapsto e^{iT_a^m B_m^a}$$
defines a homomorphism, and $SU(3)$ is a subgroup of $U(3)$. Therefore, Definition \ref{DefinitionOfTypicalStrongReferenceSystem} is compatible with $SU(3)$ theory.

\newpage

\SectionWithLabel{Affine connection representation of the unified gauge field}

\sDefinition\label{DefinitionOfElectroWeakStrongUnifiedReferenceSystem} Suppose $(M,\mathcal{F})$ and $(M,\mathcal{G})$ conform to Definition \ref{DefinitionOfTypicalGaugeField}. Let $\mathfrak{D}=r+5=8$ and both of $\mathcal{F}$ and $\mathcal{G}$ satisfy
$$G^{(\mathfrak{D}-4)(\mathfrak{D}-4)}=G^{(\mathfrak{D}-3)(\mathfrak{D}-3)},\ \ \ \ G^{(\mathfrak{D}-2)(\mathfrak{D}-2)}=G^{(\mathfrak{D}-1)(\mathfrak{D}-1)}=G^{\mathfrak{D}\mathfrak{D}}.$$ 
Thus, $\mathcal{F}$ and $\mathcal{G}$ can describe the unified field of electromagnetic, weak, and strong interactions.

\sDefinition\label{DefinitionOfWeakElectroAndColorCharge}\label{DefinitionOfWeakElectroStrongPotential} According to Definition \ref{DefinitionOfGeneralCharge}, let the charges of $\mathcal{F}$ be $\rho_{mn}$, where $m,n=4,5,\cdots,\mathfrak{D}$. Define 
$$\LeftGathered
  l\DenoteBy\left(\rho_{\left(\mathfrak{D}- 4\right)\left(\mathfrak{D}- 4\right)},\ \rho_{\left(\mathfrak{D}- 3\right)\left(\mathfrak{D}- 3\right)}\right)^T \ ,\hfill\\
  d_1 \DenoteBy(\rho_{(\mathfrak{D}- 2)(\mathfrak{D}- 2)},\ \ \rho_{(\mathfrak{D}-1)(\mathfrak{D}-1)})^T \ ,\hfill\\
  d_2 \DenoteBy(\rho_{(\mathfrak{D}-1)(\mathfrak{D}-1)},\ \ \rho_{\mathfrak{D}\mathfrak{D}})^T \ ,\hfill\\
  d_3 \DenoteBy(\rho_{\mathfrak{D}\mathfrak{D}},\ \ \rho_{(\mathfrak{D}- 2)(\mathfrak{D}- 2)})^T \ ,\hfill\\ 
\RightGathered\ \ \ 
\LeftGathered
  \nu \DenoteBy\left(\rho_{\left(\mathfrak{D}- 3\right)\left(\mathfrak{D}- 4\right)},\ \rho_{\left(\mathfrak{D}- 4\right)\left(\mathfrak{D}- 3\right)}\right)^T \ ,\hfill\\
  u_1 \DenoteBy(\rho_{(\mathfrak{D}- 2)(\mathfrak{D}-1)},\ \ \rho_{(\mathfrak{D}-1)(\mathfrak{D}- 2)})^T \ ,\hfill\\
  u_2 \DenoteBy(\rho_{(\mathfrak{D}-1)\mathfrak{D}},\ \ \rho_{\mathfrak{D}(\mathfrak{D}-1)})^T \ ,\hfill\\
  u_3 \DenoteBy(\rho_{\mathfrak{D}(\mathfrak{D}- 2)},\ \ \rho_{(\mathfrak{D}- 2)\mathfrak{D}})^T \ .\hfill\\ 
\RightGathered$$
And Denote
$$\LeftGathered
  l_L \DenoteBy\frac{1}{{\sqrt 2}}\left(\rho_{\left(\mathfrak{D}- 4\right)\left(\mathfrak{D}- 4\right)} +\rho_{\left(\mathfrak{D}- 3\right)\left(\mathfrak{D}- 3\right)}\right),\hfill\\
  l_R \DenoteBy\frac{1}{{\sqrt 2}}\left(\rho_{\left(\mathfrak{D}- 4\right)\left(\mathfrak{D}- 4\right)} -\rho_{\left(\mathfrak{D}- 3\right)\left(\mathfrak{D}- 3\right)}\right),\hfill\\ 
\RightGathered\ \ \ 
\LeftGathered
 \nu_L \DenoteBy\frac{1}{{\sqrt 2}}\left(\rho_{\left(\mathfrak{D}- 3\right)\left(\mathfrak{D}- 4\right)} +\rho_{\left(\mathfrak{D}- 4\right)\left(\mathfrak{D}- 3\right)}\right),\hfill\\
 \nu_R \DenoteBy\frac{1}{{\sqrt 2}}\left(\rho_{\left(\mathfrak{D}- 3\right)\left(\mathfrak{D}- 4\right)} -\rho_{\left(\mathfrak{D}- 4\right)\left(\mathfrak{D}- 3\right)}\right),\hfill\\ 
\RightGathered$$
$$\LeftGathered
  d_{1L} \DenoteBy\frac{1}{{\sqrt 2}}\left(\rho_{\left(\mathfrak{D}- 2\right)\left(\mathfrak{D}- 2\right)} +\rho_{\left(\mathfrak{D}-1\right)\left(\mathfrak{D}-1\right)}\right)\hfill,\\
  d_{2L} \DenoteBy\frac{1}{{\sqrt 2}}\left(\rho_{\left(\mathfrak{D}-1\right)\left(\mathfrak{D}-1\right)} +\rho_{\mathfrak{D}\mathfrak{D}}\right),\hfill\\
  d_{3L} \DenoteBy\frac{1}{{\sqrt 2}}\left(\rho_{\mathfrak{D}\mathfrak{D}} +\rho_{\left(\mathfrak{D}- 2\right)\left(\mathfrak{D}- 2\right)}\right),\hfill\\ 
\RightGathered\ \ \ 
\LeftGathered
  d_{1R} \DenoteBy\frac{1}{{\sqrt 2}}\left(\rho_{\left(\mathfrak{D}- 2\right)\left(\mathfrak{D}- 2\right)} -\rho_{\left(\mathfrak{D}-1\right)\left(\mathfrak{D}-1\right)}\right),\hfill\\
  d_{2R} \DenoteBy\frac{1}{{\sqrt 2}}\left(\rho_{\left(\mathfrak{D}-1\right)\left(\mathfrak{D}-1\right)} -\rho_{\mathfrak{D}\mathfrak{D}}\right),\hfill\\
  d_{3R} \DenoteBy\frac{1}{{\sqrt 2}}\left(\rho_{\mathfrak{D}\mathfrak{D}} -\rho_{\left(\mathfrak{D}- 2\right)\left(\mathfrak{D}- 2\right)}\right),\hfill\\ 
\RightGathered$$
$$\LeftGathered
  u_{1L} \DenoteBy\frac{1}{{\sqrt 2}}\left(\rho_{\left(\mathfrak{D}- 2\right)\left(\mathfrak{D}-1\right)} +\rho_{\left(\mathfrak{D}-1\right)\left(\mathfrak{D}- 2\right)}\right),\hfill\\
  u_{2L} \DenoteBy\frac{1}{{\sqrt 2}}\left(\rho_{\left(\mathfrak{D}-1\right)\mathfrak{D}} +\rho_{\mathfrak{D}\left(\mathfrak{D}-1\right)}\right),\hfill\\
  u_{3L} \DenoteBy\frac{1}{{\sqrt 2}}\left(\rho_{\mathfrak{D}\left(\mathfrak{D}- 2\right)} +\rho_{\left(\mathfrak{D}- 2\right)\mathfrak{D}}\right),\hfill\\ 
\RightGathered\ \ \ 
\LeftGathered
  u_{1R} \DenoteBy\frac{1}{{\sqrt 2}}\left(\rho_{\left(\mathfrak{D}- 2\right)\left(\mathfrak{D}-1\right)} -\rho_{\left(\mathfrak{D}-1\right)\left(\mathfrak{D}- 2\right)}\right),\hfill\\
  u_{2R} \DenoteBy\frac{1}{{\sqrt 2}}\left(\rho_{\left(\mathfrak{D}-1\right)\mathfrak{D}} -\rho_{\mathfrak{D}\left(\mathfrak{D}-1\right)}\right),\hfill\\
  u_{3R} \DenoteBy\frac{1}{{\sqrt 2}}\left(\rho_{\mathfrak{D}\left(\mathfrak{D}- 2\right)} -\rho_{\left(\mathfrak{D}- 2\right)\mathfrak{D}}\right).\hfill\\ 
\RightGathered$$
On $(M,\mathcal{G})$ we denote
$$\LeftGathered
g\DenoteBy\sqrt{ \left(G^{(\mathfrak{D}-4)(\mathfrak{D}-4)}\right)^2+\left(G^{(\mathfrak{D}-3)(\mathfrak{D}-3)}\right)^2},\hfill\\
\begin{aligned}
g_s&\DenoteBy\sqrt{ \left(G^{(\mathfrak{D}-1)(\mathfrak{D}-1)}\right)^2+\left(G^{\mathfrak{D}\mathfrak{D}}\right)^2}=\sqrt{ \left(G^{(\mathfrak{D}-1)(\mathfrak{D}-1)}\right)^2+\left(G^{(\mathfrak{D}-2)(\mathfrak{D}-2)}\right)^2}\\
&=\sqrt{ \left(G^{(\mathfrak{D}-2)(\mathfrak{D}-2)}\right)^2+\left(G^{\mathfrak{D}\mathfrak{D}}\right)^2},
\end{aligned}\hfill\\
\RightGathered$$
$$\LeftGathered
  Z_P \DenoteBy\frac{1}{{\sqrt 2}}(\Gamma_{(\mathfrak{D}- 4)(\mathfrak{D}- 4)P} +\Gamma_{(\mathfrak{D}- 3)(\mathfrak{D}- 3)P}),\hfill\\
  A_P \DenoteBy\frac{1}{{\sqrt 2}}(\Gamma_{(\mathfrak{D}- 4)(\mathfrak{D}- 4)P} -\Gamma_{(\mathfrak{D}- 3)(\mathfrak{D}- 3)P}),\hfill\\ 
\RightGathered\ \ \ 
\LeftGathered
  W^1_P \DenoteBy\frac{1}{{\sqrt 2}}(\Gamma_{(\mathfrak{D}- 4)(\mathfrak{D}- 3)P} +\Gamma_{(\mathfrak{D}- 3)(\mathfrak{D}- 4)P}),\hfill\\
  W^2_P \DenoteBy\frac{1}{{\sqrt 2}}(\Gamma_{(\mathfrak{D}- 4)(\mathfrak{D}- 3)P} -\Gamma_{(\mathfrak{D}- 3)(\mathfrak{D}- 4)P}),\hfill\\ 
\RightGathered$$
$$\LeftGathered
  U^1_P \DenoteBy\frac{1}{{\sqrt 2}}\left(\Gamma_{(\mathfrak{D}- 2)(\mathfrak{D}- 2)P} +\Gamma_{(\mathfrak{D}-1)(\mathfrak{D}-1)P}\right),\hfill\\
  V^1_P \DenoteBy\frac{1}{{\sqrt 2}}\left(\Gamma_{(\mathfrak{D}- 2)(\mathfrak{D}- 2)P} -\Gamma_{(\mathfrak{D}-1)(\mathfrak{D}-1)P}\right),\hfill\\ 
\RightGathered\ \ \ 
\LeftGathered
  X^{23}_P \DenoteBy\frac{1}{{\sqrt 2}}\left(\Gamma_{(\mathfrak{D}- 2)(\mathfrak{D}-1)P} +\Gamma_{(\mathfrak{D}-1)(\mathfrak{D}- 2)P}\right),\hfill\\
  Y^{23}_P \DenoteBy\frac{1}{{\sqrt 2}}\left(\Gamma_{(\mathfrak{D}- 2)(\mathfrak{D}-1)P} -\Gamma_{(\mathfrak{D}-1)(\mathfrak{D}- 2)P}\right),\hfill\\ 
\RightGathered$$
$$\LeftGathered
  U^2_P \DenoteBy\frac{1}{{\sqrt 2}}\left(\Gamma_{(\mathfrak{D}-1)(\mathfrak{D}-1)P} +\Gamma_{\mathfrak{D}\mathfrak{D}P}\right),\hfill\\
  V^2_P \DenoteBy\frac{1}{{\sqrt 2}}\left(\Gamma_{(\mathfrak{D}-1)(\mathfrak{D}-1)P} -\Gamma_{\mathfrak{D}\mathfrak{D}P}\right),\hfill\\ 
\RightGathered\ \ \ 
\LeftGathered
  X^{31}_P \DenoteBy\frac{1}{{\sqrt 2}}\left(\Gamma_{(\mathfrak{D}-1)\mathfrak{D}P} +\Gamma_{\mathfrak{D}(\mathfrak{D}-1)P}\right),\hfill\\
  Y^{31}_P \DenoteBy\frac{1}{{\sqrt 2}}\left(\Gamma_{(\mathfrak{D}-1)\mathfrak{D}P} -\Gamma_{\mathfrak{D}(\mathfrak{D}-1)P}\right),\hfill\\ 
\RightGathered$$
$$\LeftGathered
  U^3_P \DenoteBy\frac{1}{{\sqrt 2}}\left(\Gamma_{\mathfrak{D}\mathfrak{D}P} +\Gamma_{(\mathfrak{D}- 2)(\mathfrak{D}- 2)P}\right),\hfill\\
  V^3_P \DenoteBy\frac{1}{{\sqrt 2}}\left(\Gamma_{\mathfrak{D}\mathfrak{D}P} -\Gamma_{(\mathfrak{D}- 2)(\mathfrak{D}- 2)P}\right),\hfill\\ 
\RightGathered\ \ \ 
\LeftGathered
  X^{12}_P \DenoteBy\frac{1}{{\sqrt 2}}\left(\Gamma_{\mathfrak{D}(\mathfrak{D}- 2)P} +\Gamma_{(\mathfrak{D}- 2)\mathfrak{D}P}\right),\hfill\\
  Y^{12}_P \DenoteBy\frac{1}{{\sqrt 2}}\left(\Gamma_{\mathfrak{D}(\mathfrak{D}- 2)P} -\Gamma_{(\mathfrak{D}- 2)\mathfrak{D}P}\right).\hfill\\ 
\RightGathered$$

\newpage

\sDiscussion We know from section \ref{Gauge transformation in affine connection representation} that the gauge frame matrix $[B^a_m]\in GL(5,\mathbb{R}),\ \ (a,m=4,5,\cdots,8)$; therefore, when $B^a_m$ are without any constraints, we can obtain a $GL(5,\mathbb{R})$ gauge theory. In consideration of the fact that the exponential map
$$exp:\ GL(5,\mathbb R)\to U(5),\ \ [B_m^a]\mapsto e^{iT_a^m B_m^a}$$
defines a homomorphism and $U(1)\times SU(2) \times SU(3)$ is a subgroup of $U(5)$. So there must exist some constraint conditions of $B^a_m$ to make $GL(5,\mathbb{R})$ reduced to $U(1)\times SU(2) \times SU(3)$, i.e., 
$$GL(5,\mathbb{R}) \xrightarrow{\ \ \ \ \mbox{constraint conditions of $B^a_m$}\ \ \ \ } U(1)\times SU(2) \times SU(3).$$
More generally, suppose we do not know what the symmetry that can exactly describe "the real world" is, we just denote it by $S$; then, the map
$$GL(5,\mathbb{R}) \xrightarrow{\ \ \ \ \mbox{constraint conditions of $B^a_m$}\ \ \ \ } S$$
makes us be able to turn the problem of seeking for $S$ into the problem of seeking for a set of constraint conditions of $B^a_m$. "To describe $S$" and "to describe the constraint conditions of $B^a_m$" are equivalent to each other.

Because gauge potentials $\Gamma_{mnP}$ and particle fields $\rho_{mn}$ are both constructed from the gauge frame field $B^a_m$, clearly here, it is more flexible and convenient "to describe the constraint conditions of $B^a_m$" than "to describe $S$".

Next, we have no idea what the best constraint conditions look like, but we can try to define a set of constraint conditions to see what can be obtained.

\sDefinition\label{DefinitionOfUnifiedSymmetryCondition} Similar to Remark \ref{DiscussionWaysOfGeometryAndGroup}, we define the constraint conditions as follows.

(1) 1st basic conditions:
$$\LeftGathered
G^{(\mathfrak{D}-4)(\mathfrak{D}-4)}=G^{(\mathfrak{D}-3)(\mathfrak{D}-3)},\hfill\\
G^{(\mathfrak{D}-2)(\mathfrak{D}-2)}=G^{(\mathfrak{D}-1)(\mathfrak{D}-1)}=G^{\mathfrak{D}\mathfrak{D}},\hfill\\
\RightGathered$$

(2) 2nd basic conditions:
$$\LeftGathered
\Gamma_{(\mathfrak{D}- 3)(\mathfrak{D}- 4)P}=\Gamma_{(\mathfrak{D}- 4)(\mathfrak{D}- 3)P},\hfill\\
\Gamma_{(\mathfrak{D}- 2)(\mathfrak{D}- 2)P}+\Gamma_{(\mathfrak{D}- 1)(\mathfrak{D}- 1)P}+\Gamma_{\mathfrak{D}\mathfrak{D}P}=0,\hfill\\
\RightGathered$$

(3) 1st conditions of PMNS mixing of leptons:
$$\LeftGathered
 \Gamma_{(\mathfrak{D}- 4)P}^{\mathfrak{D}- 2}=c_{\mathfrak{D}- 3}^{\mathfrak{D}- 2}\Gamma_{(\mathfrak{D}- 4)P}^{\mathfrak{D}- 3},\hfill\\
 \Gamma_{(\mathfrak{D}- 4)P}^{\mathfrak{D}-1}=c_{\mathfrak{D}- 3}^{\mathfrak{D}-1}\Gamma_{(\mathfrak{D}- 4)P}^{\mathfrak{D}- 3},\hfill\\
 \Gamma_{(\mathfrak{D}- 4)P}^\mathfrak{D}=c_{\mathfrak{D}- 3}^\mathfrak{D}\Gamma_{(\mathfrak{D}- 4)P}^{\mathfrak{D}- 3},\hfill\\ 
\RightGathered\ \ \ 
\LeftGathered
 \Gamma_{(\mathfrak{D}- 3)P}^{\mathfrak{D}- 2}=c_{\mathfrak{D}- 4}^{\mathfrak{D}- 2}\Gamma_{(\mathfrak{D}- 3)P}^{\mathfrak{D}- 4},\hfill\\
 \Gamma_{(\mathfrak{D}- 3)P}^{\mathfrak{D}-1}=c_{\mathfrak{D}- 4}^{\mathfrak{D}-1}\Gamma_{(\mathfrak{D}- 3)P}^{\mathfrak{D}- 4},\hfill\\
 \Gamma_{(\mathfrak{D}- 3)P}^\mathfrak{D}=c_{\mathfrak{D}- 4}^\mathfrak{D}\Gamma_{(\mathfrak{D}- 3)P}^{\mathfrak{D}- 4},\hfill\\ 
\RightGathered\ \ \ 
\LeftGathered
  c_{\mathfrak{D}- 3}^{\mathfrak{D}- 2}=c_{\mathfrak{D}- 4}^{\mathfrak{D}- 2},\hfill\\
  c_{\mathfrak{D}- 3}^{\mathfrak{D}-1}=c_{\mathfrak{D}- 4}^{\mathfrak{D}-1},\hfill\\
  c_{\mathfrak{D}- 3}^\mathfrak{D}=c_{\mathfrak{D}- 4}^\mathfrak{D},\hfill\\ 
\RightGathered$$

(4) 2nd conditions of PMNS mixing of leptons:
$$\LeftGathered
 \rho_{(\mathfrak{D}- 2)(\mathfrak{D}- 3)}=\rho_{(\mathfrak{D}- 2)(\mathfrak{D}- 4)},\hfill\\
 \rho_{(\mathfrak{D}-1)(\mathfrak{D}- 3)}=\rho_{(\mathfrak{D}-1)(\mathfrak{D}- 4)},\hfill\\
 \rho_{\mathfrak{D}(\mathfrak{D}- 3)}=\rho_{\mathfrak{D}(\mathfrak{D}- 4)},\hfill\\ 
\RightGathered\ \ \ 
\LeftGathered
 \rho_{(\mathfrak{D}- 3)(\mathfrak{D}- 2)}=\rho_{(\mathfrak{D}- 4)(\mathfrak{D}- 2)},\hfill\\
 \rho_{(\mathfrak{D}- 3)(\mathfrak{D}-1)}=\rho_{(\mathfrak{D}- 4)(\mathfrak{D}-1)},\hfill\\
 \rho_{(\mathfrak{D}- 3)\mathfrak{D}}=\rho_{(\mathfrak{D}- 4)\mathfrak{D}},\hfill\\ 
\RightGathered$$

(5) 1st conditions of CKM mixing of quarks:
$$\LeftGathered
 \Gamma_{(\mathfrak{D}- 2)P}^{\mathfrak{D}- 3}=c_{\mathfrak{D}- 2}^{\mathfrak{D}- 4}\Gamma_{(\mathfrak{D}- 4)P}^{\mathfrak{D}- 3},\hfill\\
 \Gamma_{(\mathfrak{D}-1)P}^{\mathfrak{D}- 3}=c_{\mathfrak{D}-1}^{\mathfrak{D}- 4}\Gamma_{(\mathfrak{D}- 4)P}^{\mathfrak{D}- 3},\hfill\\
 \Gamma_{\mathfrak{D}P}^{\mathfrak{D}- 3}=c_\mathfrak{D}^{\mathfrak{D}- 4}\Gamma_{(\mathfrak{D}- 4)P}^{\mathfrak{D}- 3},\hfill\\ 
\RightGathered\ \ \ 
\LeftGathered
 \Gamma_{(\mathfrak{D}- 2)P}^{\mathfrak{D}- 4}=c_{\mathfrak{D}- 2}^{\mathfrak{D}- 3}\Gamma_{(\mathfrak{D}- 3)P}^{\mathfrak{D}- 4},\hfill\\
 \Gamma_{(\mathfrak{D}-1)P}^{\mathfrak{D}- 4}=c_{\mathfrak{D}-1}^{\mathfrak{D}- 3}\Gamma_{(\mathfrak{D}- 3)P}^{\mathfrak{D}- 4},\hfill\\
 \Gamma_{\mathfrak{D}P}^{\mathfrak{D}- 4}=c_\mathfrak{D}^{\mathfrak{D}- 3}\Gamma_{(\mathfrak{D}- 3)P}^{\mathfrak{D}- 4},\hfill\\ 
\RightGathered\ \ \ 
\LeftGathered
  c_{\mathfrak{D}- 2}^{\mathfrak{D}- 4}=c_{\mathfrak{D}-1}^{\mathfrak{D}- 4}=c_\mathfrak{D}^{\mathfrak{D}- 4},\hfill\\
  c_{\mathfrak{D}- 2}^{\mathfrak{D}- 3}=c_{\mathfrak{D}-1}^{\mathfrak{D}- 3}=c_\mathfrak{D}^{\mathfrak{D}- 3},\hfill\\ 
\RightGathered$$

(6) 2nd conditions of CKM mixing of quarks:
$$\LeftGathered
 \rho_{(\mathfrak{D}- 2)(\mathfrak{D}- 3)}=\rho_{(\mathfrak{D}-1)(\mathfrak{D}- 3)}=\rho_{\mathfrak{D}(\mathfrak{D}- 3)},\hfill\\
 \rho_{(\mathfrak{D}- 2)(\mathfrak{D}- 4)}=\rho_{(\mathfrak{D}-1)(\mathfrak{D}- 4)}=\rho_{\mathfrak{D}(\mathfrak{D}- 4)},\hfill\\ 
\RightGathered\ \ \ 
\LeftGathered
 \rho_{(\mathfrak{D}- 3)(\mathfrak{D}- 2)}=\rho_{(\mathfrak{D}- 3)(\mathfrak{D}-1)}=\rho_{(\mathfrak{D}- 3)\mathfrak{D}},\hfill\\
 \rho_{(\mathfrak{D}- 4)(\mathfrak{D}- 2)}=\rho_{(\mathfrak{D}- 4)(\mathfrak{D}-1)}=\rho_{(\mathfrak{D}- 4)\mathfrak{D}},\hfill\\ 
\RightGathered$$
where $c_n^m$ are constants.

\newpage

\sProposition\label{PropositionOfMNSMixing} When $(M,\mathcal{F})$ and $(M,\mathcal{G})$ satisfy the symmetry conditions (1)(2)(3)(4) of Definition \ref{DefinitionOfUnifiedSymmetryCondition}, denote
$$\begin{aligned}
  l' &\DenoteBy \left( \rho_{\left(\mathfrak{D}- 4\right)\left(\mathfrak{D}- 4\right)} +\frac{{c_{\mathfrak{D}- 4}^{\mathfrak{D}- 2}}}{2}\left(\rho_{\left(\mathfrak{D}- 2\right)\left(\mathfrak{D}- 4\right)} +\rho_{\left(\mathfrak{D}- 4\right)\left(\mathfrak{D}- 2\right)}\right) \right. \\
     &\ \ \ \ \ \ \ \ \ \ \ \ \ \ \ \ \ \ \ \ \ \ +\frac{{c_{\mathfrak{D}- 4}^{\mathfrak{D}-1}}}{2}\left(\rho_{\left(\mathfrak{D}-1\right)\left(\mathfrak{D}- 4\right)} +\rho_{\left(\mathfrak{D}- 4\right)\left(\mathfrak{D}-1\right)}\right)+\frac{{c_{\mathfrak{D}- 4}^\mathfrak{D}}}{2}\left(\rho_{\mathfrak{D}\left(\mathfrak{D}- 4\right)} +\rho_{\left(\mathfrak{D}- 4\right)\mathfrak{D}}\right), \\
     &\ \ \ \ \ \ \ \ \rho_{\left(\mathfrak{D}- 3\right)\left(\mathfrak{D}- 3\right)} +\frac{{c_{\mathfrak{D}- 3}^{\mathfrak{D}- 2}}}{2}\left(\rho_{\left(\mathfrak{D}- 2\right)\left(\mathfrak{D}- 3\right)} +\rho_{\left(\mathfrak{D}- 3\right)\left(\mathfrak{D}- 2\right)}\right) \\
     &\ \ \ \ \ \ \ \ \ \ \ \ \ \ \ \ \ \ \ \ \ \ \left. +\frac{{c_{\mathfrak{D}- 3}^{\mathfrak{D}-1}}}{2}\left(\rho_{\left(\mathfrak{D}-1\right)\left(\mathfrak{D}- 3\right)} +\rho_{\left(\mathfrak{D}- 3\right)\left(\mathfrak{D}-1\right)}\right)+\frac{{c_{\mathfrak{D}- 3}^\mathfrak{D}}}{2}\left(\rho_{\mathfrak{D}\left(\mathfrak{D}- 3\right)} +\rho_{\left(\mathfrak{D}- 3\right)\mathfrak{D}}\right) \right)^T,\\
  \nu' &\DenoteBy \left( \rho_{\left(\mathfrak{D}- 3\right)\left(\mathfrak{D}- 4\right)} +\frac{{c_{\mathfrak{D}- 3}^{\mathfrak{D}- 2}}}{2}\left(\rho_{\left(\mathfrak{D}- 2\right)\left(\mathfrak{D}- 4\right)} +\rho_{\left(\mathfrak{D}- 4\right)\left(\mathfrak{D}- 2\right)}\right) \right. \\
     &\ \ \ \ \ \ \ \ \ \ \ \ \ \ \ \ \ \ \ \ \ \ +\frac{{c_{\mathfrak{D}- 3}^{\mathfrak{D}-1}}}{2}\left(\rho_{\left(\mathfrak{D}-1\right)\left(\mathfrak{D}- 4\right)} +\rho_{\left(\mathfrak{D}- 4\right)\left(\mathfrak{D}-1\right)}\right)+\frac{{c_{\mathfrak{D}- 3}^\mathfrak{D}}}{2}\left(\rho_{\mathfrak{D}\left(\mathfrak{D}- 4\right)} +\rho_{\left(\mathfrak{D}- 4\right)\mathfrak{D}}\right), \\
     &\ \ \ \ \ \ \ \ \rho_{\left(\mathfrak{D}- 4\right)\left(\mathfrak{D}- 3\right)} +\frac{{c_{\mathfrak{D}- 4}^{\mathfrak{D}- 2}}}{2}\left(\rho_{\left(\mathfrak{D}- 2\right)\left(\mathfrak{D}- 3\right)} +\rho_{\left(\mathfrak{D}- 3\right)\left(\mathfrak{D}- 2\right)}\right) \\
     &\ \ \ \ \ \ \ \ \ \ \ \ \ \ \ \ \ \ \ \ \ \ \left. +\frac{{c_{\mathfrak{D}- 4}^{\mathfrak{D}-1}}}{2}\left(\rho_{\left(\mathfrak{D}-1\right)\left(\mathfrak{D}- 3\right)} +\rho_{\left(\mathfrak{D}- 3\right)\left(\mathfrak{D}-1\right)}\right)+\frac{{c_{\mathfrak{D}- 4}^\mathfrak{D}}}{2}\left(\rho_{\mathfrak{D}\left(\mathfrak{D}- 3\right)} +\rho_{\left(\mathfrak{D}- 3\right)\mathfrak{D}}\right) \right)^T.
\end{aligned}$$
Then, the geometric properties $l$ and $\nu$ of $\mathcal{F}$ satisfy the following conclusions on $(M,\mathcal{G})$.
\begin{equation}\label{ElectronMNSMixing}\LeftGathered
  l_{L;P}=\partial_P l_L  -g l_L Z_P  -g l_R A_P  -g\nu '_L W^1_P, \hfill\\
  l_{R;P}=\partial_P l_R  -g l_R Z_P  -g l_L A_P, \hfill\\
 \nu_{L;P}=\partial_P\nu_L  -g\nu_L Z_P  -g l'_L W^1_P, \hfill\\
 \nu_{R;P}=\partial_P\nu_R  -g\nu_R Z_P. \hfill\\ 
\RightGathered\end{equation}

\sProof First, we compute the covariant differential of $\rho_{mn}$ of $\mathcal{F}$.
$$\begin{aligned}
\rho_{mn;P}&=\partial_P\rho_{mn} -\rho_{Hn}\Gamma_{mP}^H  -\rho_{mH}\Gamma_{nP}^H\\
&=\partial_P\rho_{mn} -\rho_{(\mathfrak{D}- 4)n}\Gamma_{mP}^{\mathfrak{D}- 4} -\rho_{(\mathfrak{D}- 3)n}\Gamma_{mP}^{\mathfrak{D}- 3} -\rho_{(\mathfrak{D}- 2)n}\Gamma_{mP}^{\mathfrak{D}- 2} -\rho_{(\mathfrak{D}-1)n}\Gamma_{mP}^{\mathfrak{D}-1} -\rho_{\mathfrak{D}n}\Gamma_{mP}^\mathfrak{D}\\
&-\rho_{m(\mathfrak{D}- 4)}\Gamma_{nP}^{\mathfrak{D}- 4} -\rho_{m(\mathfrak{D}- 3)}\Gamma_{nP}^{\mathfrak{D}- 3} -\rho_{m(\mathfrak{D}- 2)}\Gamma_{nP}^{\mathfrak{D}- 2} -\rho_{m(\mathfrak{D}-1)}\Gamma_{nP}^{\mathfrak{D}-1} -\rho_{m\mathfrak{D}}\Gamma_{nP}^\mathfrak{D}.\\ 
\end{aligned}$$
According to Definition \ref{DefinitionOfWeakElectroStrongPotential} and Definition \ref{DefinitionOfUnifiedSymmetryCondition}, by calculation, we obtain that
$$\LeftList
\begin{aligned}
l_{L;P}&=\partial_P l_L  -g l_L Z_P  -g l_R A_P  -g\nu_L W^1_P \\
&-\frac{1}{2}\left[{c_{\mathfrak{D}- 4}^{\mathfrak{D}- 2}\left(\rho_{\left(\mathfrak{D}- 2\right)\left(\mathfrak{D}- 3\right)} +\rho_{\left(\mathfrak{D}- 3\right)\left(\mathfrak{D}- 2\right)}\right)+ c_{\mathfrak{D}- 3}^{\mathfrak{D}- 2}\left(\rho_{\left(\mathfrak{D}- 2\right)\left(\mathfrak{D}- 4\right)} +\rho_{\left(\mathfrak{D}- 4\right)\left(\mathfrak{D}- 2\right)}\right)}\right]\frac{g}{\sqrt 2}W^1_P \\
\end{aligned}\hfill\\
\begin{aligned}
&-\frac{1}{2}\left[{c_{\mathfrak{D}- 4}^{\mathfrak{D}-1}\left(\rho_{\left(\mathfrak{D}-1\right)\left(\mathfrak{D}- 3\right)} +\rho_{\left(\mathfrak{D}- 3\right)\left(\mathfrak{D}-1\right)}\right)+ c_{\mathfrak{D}- 3}^{\mathfrak{D}-1}\left(\rho_{\left(\mathfrak{D}-1\right)\left(\mathfrak{D}- 4\right)} +\rho_{\left(\mathfrak{D}- 4\right)\left(\mathfrak{D}-1\right)}\right)}\right]\frac{g}{\sqrt 2}W^1_P \\
&-\frac{1}{2}\left[{c_{\mathfrak{D}- 4}^\mathfrak{D}\left(\rho_{\mathfrak{D}\left(\mathfrak{D}- 3\right)} +\rho_{\left(\mathfrak{D}- 3\right)\mathfrak{D}}\right)+ c_{\mathfrak{D}- 3}^\mathfrak{D}\left(\rho_{\mathfrak{D}\left(\mathfrak{D}- 4\right)} +\rho_{\left(\mathfrak{D}- 4\right)\mathfrak{D}}\right)}\right]\frac{g}{\sqrt 2}W^1_P, \\ 
l_{R;P}&=\partial_P l_R  -g l_R Z_P  -g l_L A_P,\\
\end{aligned}\hfill\\
\RightList$$
$$\LeftList
\begin{aligned}
\nu_{L;P}&=\partial_P\nu_L  -g\nu_L Z_P  -g l_L W^1_P \\
&-\frac{1}{2}\left[{c_{\mathfrak{D}- 4}^{\mathfrak{D}- 2}\left(\rho_{\left(\mathfrak{D}- 2\right)\left(\mathfrak{D}- 4\right)} +\rho_{\left(\mathfrak{D}- 4\right)\left(\mathfrak{D}- 2\right)}\right)+ c_{\mathfrak{D}- 3}^{\mathfrak{D}- 2}\left(\rho_{\left(\mathfrak{D}- 3\right)\left(\mathfrak{D}- 2\right)} +\rho_{\left(\mathfrak{D}- 2\right)\left(\mathfrak{D}- 3\right)}\right)}\right]\frac{g}{\sqrt 2}W^1_P \\
&-\frac{1}{2}\left[{c_{\mathfrak{D}- 4}^{\mathfrak{D}-1}\left(\rho_{\left(\mathfrak{D}-1\right)\left(\mathfrak{D}- 4\right)} +\rho_{\left(\mathfrak{D}- 4\right)\left(\mathfrak{D}-1\right)}\right)+ c_{\mathfrak{D}- 3}^{\mathfrak{D}-1}\left(\rho_{\left(\mathfrak{D}- 3\right)\left(\mathfrak{D}-1\right)} +\rho_{\left(\mathfrak{D}-1\right)\left(\mathfrak{D}- 3\right)}\right)}\right]\frac{g}{\sqrt 2}W^1_P \\
&-\frac{1}{2}\left[{c_{\mathfrak{D}- 4}^\mathfrak{D}\left(\rho_{\mathfrak{D}\left(\mathfrak{D}- 4\right)} +\rho_{\left(\mathfrak{D}- 4\right)\mathfrak{D}}\right)+ c_{\mathfrak{D}- 3}^\mathfrak{D}\left(\rho_{\left(\mathfrak{D}- 3\right)\mathfrak{D}} +\rho_{\mathfrak{D}\left(\mathfrak{D}- 3\right)}\right)}\right]\frac{g}{\sqrt 2}W^1_P, \\ 
\end{aligned}\hfill\\
\begin{aligned}
\nu_{R;P}&=\partial_P\nu_R  -g\nu_R Z_P.
\end{aligned}\hfill\\
\RightList$$
Then, according to definitions of $l'$ and $\nu'$, we obtain that
$$\begin{aligned}
  l'_L &=l_L  \\
&+\frac{{c_{\mathfrak{D}- 4}^{\mathfrak{D}- 2}}}{{2\sqrt 2}}\left(\rho_{\left(\mathfrak{D}- 2\right)\left(\mathfrak{D}- 4\right)} +\rho_{\left(\mathfrak{D}- 4\right)\left(\mathfrak{D}- 2\right)}\right)+\frac{{c_{\mathfrak{D}- 4}^{\mathfrak{D}-1}}}{{2\sqrt 2}}\left(\rho_{\left(\mathfrak{D}-1\right)\left(\mathfrak{D}- 4\right)} +\rho_{\left(\mathfrak{D}- 4\right)\left(\mathfrak{D}-1\right)}\right)+\frac{{c_{\mathfrak{D}- 4}^\mathfrak{D}}}{{2\sqrt 2}}\left(\rho_{\mathfrak{D}\left(\mathfrak{D}- 4\right)} +\rho_{\left(\mathfrak{D}- 4\right)\mathfrak{D}}\right)\hfill\\
&+\frac{{c_{\mathfrak{D}- 3}^{\mathfrak{D}- 2}}}{{2\sqrt 2}}\left(\rho_{\left(\mathfrak{D}- 2\right)\left(\mathfrak{D}- 3\right)} +\rho_{\left(\mathfrak{D}- 3\right)\left(\mathfrak{D}- 2\right)}\right)+\frac{{c_{\mathfrak{D}- 3}^{\mathfrak{D}-1}}}{{2\sqrt 2}}\left(\rho_{\left(\mathfrak{D}-1\right)\left(\mathfrak{D}- 3\right)} +\rho_{\left(\mathfrak{D}- 3\right)\left(\mathfrak{D}-1\right)}\right)+\frac{{c_{\mathfrak{D}- 3}^\mathfrak{D}}}{{2\sqrt 2}}\left(\rho_{\mathfrak{D}\left(\mathfrak{D}- 3\right)} +\rho_{\left(\mathfrak{D}- 3\right)\mathfrak{D}}\right), \\
 \nu '_L &=\nu_L  \\
&+\frac{{c_{\mathfrak{D}- 3}^{\mathfrak{D}- 2}}}{{2\sqrt 2}}\left(\rho_{\left(\mathfrak{D}- 2\right)\left(\mathfrak{D}- 4\right)} +\rho_{\left(\mathfrak{D}- 4\right)\left(\mathfrak{D}- 2\right)}\right)+\frac{{c_{\mathfrak{D}- 3}^{\mathfrak{D}-1}}}{{2\sqrt 2}}\left(\rho_{\left(\mathfrak{D}-1\right)\left(\mathfrak{D}- 4\right)} +\rho_{\left(\mathfrak{D}- 4\right)\left(\mathfrak{D}-1\right)}\right)+\frac{{c_{\mathfrak{D}- 3}^\mathfrak{D}}}{{2\sqrt 2}}\left(\rho_{\mathfrak{D}\left(\mathfrak{D}- 4\right)} +\rho_{\left(\mathfrak{D}- 4\right)\mathfrak{D}}\right)\hfill\\
&+\frac{{c_{\mathfrak{D}- 4}^{\mathfrak{D}- 2}}}{{2\sqrt 2}}\left(\rho_{\left(\mathfrak{D}- 2\right)\left(\mathfrak{D}- 3\right)} +\rho_{\left(\mathfrak{D}- 3\right)\left(\mathfrak{D}- 2\right)}\right)+\frac{{c_{\mathfrak{D}- 4}^{\mathfrak{D}-1}}}{{2\sqrt 2}}\left(\rho_{\left(\mathfrak{D}-1\right)\left(\mathfrak{D}- 3\right)} +\rho_{\left(\mathfrak{D}- 3\right)\left(\mathfrak{D}-1\right)}\right)+\frac{{c_{\mathfrak{D}- 4}^\mathfrak{D}}}{{2\sqrt 2}}\left(\rho_{\mathfrak{D}\left(\mathfrak{D}- 3\right)} +\rho_{\left(\mathfrak{D}- 3\right)\mathfrak{D}}\right).\hfill\\ 
\end{aligned}$$
Substitute them into the previous equations, and we obtain that
$$\begin{aligned}
&\LeftGathered
  l_{L;P}=\partial_P l_L  -g l_L Z_P  -g l_R A_P  -g\nu '_L W^1_P, \hfill\\
  l_{R;P}=\partial_P l_R  -g l_R Z_P  -g l_L A_P, \hfill\\
\RightGathered\\
&\LeftGathered
 \nu_{L;P}=\partial_P\nu_L  -g\nu_L Z_P  -g l'_L W^1_P, \hfill\\
 \nu_{R;P}=\partial_P\nu_R  -g\nu_R Z_P.\ \ \ \ \ \ \ \ \ \ \ \ \ \ \ \ \  \hfill\\ 
\RightGathered\\
\end{aligned}$$
\qed

\sRemark\label{RemarkThreeGenerationOfLeptons} The above proposition shows the geometric origin of PMNS mixing of weak interaction. In affine connection representation of gauge fields, PMNS mixing arises as a geometric property on manifold.

In conventional physics, $e$, $\mu$ and $\tau$ have just only ontological differences, but they have no difference in mathematical connotation. By contrast, Proposition \ref{PropositionOfMNSMixing} tells us that leptons of three generations should be constructed by different linear combinations of $\{\rho_{pq},\rho_{qp}\}_{p=4,5; \ q=6,7,8}$. Thus, $e$, $\mu$, and $\tau$ may have concrete and  distinguishable mathematical connotations. For example, let $a_\mu$, $b_\mu$, ${a_\mu}_n^m$, ${b_\mu}_n^m$, $a_\tau$, $b_\tau$, ${a_\tau}_n^m$, ${b_\tau}_n^m$ be constants; then, we might suppose that
$$\LeftList
\LeftGathered
  e\DenoteBy l=(\rho_{(\mathfrak{D}- 4)(\mathfrak{D}- 4)},\ \ \rho_{(\mathfrak{D}- 3)(\mathfrak{D}- 3)})^T,\hfill\\
 \nu_e \DenoteBy\nu =(\rho_{(\mathfrak{D}- 3)(\mathfrak{D}- 4)},\ \ \rho_{(\mathfrak{D}- 4)(\mathfrak{D}- 3)})^T.\hfill\\ 
\RightGathered \hfill\\
\LeftGathered
\begin{aligned}
 \mu \DenoteBy a_\mu e +\frac{1}{2}&\left({a_\mu} _{\mathfrak{D}- 4}^{\mathfrak{D}- 2}\rho_{(\mathfrak{D}- 2)(\mathfrak{D}- 4)} + {a_\mu} _{\mathfrak{D}- 4}^{\mathfrak{D}-1}\rho_{(\mathfrak{D}-1)(\mathfrak{D}- 4)} + {a_\mu} _{\mathfrak{D}- 4}^\mathfrak{D}\rho_{\mathfrak{D}(\mathfrak{D}- 4)},\right.\\
&\left.{a_\mu} _{\mathfrak{D}- 3}^{\mathfrak{D}- 2}\rho_{(\mathfrak{D}- 2)(\mathfrak{D}- 3)} + {a_\mu} _{\mathfrak{D}- 3}^{\mathfrak{D}-1}\rho_{(\mathfrak{D}-1)(\mathfrak{D}- 3)} + {a_\mu} _{\mathfrak{D}- 3}^\mathfrak{D}\rho_{\mathfrak{D}(\mathfrak{D}- 3)}\right)^T. \hfill\\
\end{aligned}\hfill\\ 
\begin{aligned}
 \nu_\mu \DenoteBy b_\mu \nu_e +\frac{1}{2}&\left({b_\mu} _{\mathfrak{D}- 3}^{\mathfrak{D}- 2}\rho_{(\mathfrak{D}- 2)(\mathfrak{D}- 4)} + {b_\mu} _{\mathfrak{D}- 3}^{\mathfrak{D}-1}\rho_{(\mathfrak{D}-1)(\mathfrak{D}- 4)} + {b_\mu} _{\mathfrak{D}- 3}^\mathfrak{D}\rho_{\mathfrak{D}(\mathfrak{D}- 4)},\right.\\
&\left.{b_\mu} _{\mathfrak{D}- 4}^{\mathfrak{D}- 2}\rho_{(\mathfrak{D}- 2)(\mathfrak{D}- 3)} + {b_\mu} _{\mathfrak{D}- 4}^{\mathfrak{D}-1}\rho_{(\mathfrak{D}-1)(\mathfrak{D}- 3)} + {b_\mu} _{\mathfrak{D}- 4}^\mathfrak{D}\rho_{\mathfrak{D}(\mathfrak{D}- 3)}\right)^T. \hfill\\ 
\end{aligned}\hfill\\ 
\RightGathered \hfill\\
\LeftGathered
\begin{aligned}
 \tau \DenoteBy {a_\tau} \mu  +\frac{1}{2}&\left({a_\tau} _{\mathfrak{D}- 4}^{\mathfrak{D}- 2}\rho_{(\mathfrak{D}- 4)(\mathfrak{D}- 2)} + {a_\tau} _{\mathfrak{D}- 4}^{\mathfrak{D}-1}\rho_{(\mathfrak{D}- 4)(\mathfrak{D}-1)} + {a_\tau} _{\mathfrak{D}- 4}^\mathfrak{D}\rho_{(\mathfrak{D}- 4)\mathfrak{D}},\right.\\ 
&\left.{a_\tau} _{\mathfrak{D}- 3}^{\mathfrak{D}- 2}\rho_{(\mathfrak{D}- 3)(\mathfrak{D}- 2)} + {a_\tau} _{\mathfrak{D}- 3}^{\mathfrak{D}-1}\rho_{(\mathfrak{D}- 3)(\mathfrak{D}-1)} + {a_\tau} _{\mathfrak{D}- 3}^\mathfrak{D}\rho_{(\mathfrak{D}- 3)\mathfrak{D}}\right)^T. \hfill\\
\end{aligned}\hfill\\ 
\begin{aligned}
 \nu_\tau  \DenoteBy {b_\tau} \nu_\mu +\frac{1}{2}&\left({b_\tau} _{\mathfrak{D}- 3}^{\mathfrak{D}- 2}\rho_{(\mathfrak{D}- 4)(\mathfrak{D}- 2)} + {b_\tau} _{\mathfrak{D}- 3}^{\mathfrak{D}-1}\rho_{(\mathfrak{D}- 4)(\mathfrak{D}-1)} + {b_\tau} _{\mathfrak{D}- 3}^\mathfrak{D}\rho_{(\mathfrak{D}- 4)\mathfrak{D}},\right.\\ 
&\left.{b_\tau} _{\mathfrak{D}- 4}^{\mathfrak{D}- 2}\rho_{(\mathfrak{D}- 3)(\mathfrak{D}- 2)} + {b_\tau} _{\mathfrak{D}- 4}^{\mathfrak{D}-1}\rho_{(\mathfrak{D}- 3)(\mathfrak{D}-1)} + {b_\tau} _{\mathfrak{D}- 4}^\mathfrak{D}\rho_{(\mathfrak{D}- 3)\mathfrak{D}}\right)^T. \hfill\\ 
\end{aligned} \hfill\\ 
\RightGathered \hfill\\
\RightList$$

\sProposition\label{PropositionOfCKMMixing} When $(M,\mathcal{F})$ and $(M,\mathcal{G})$ satisfy the symmetry conditions (1)(2)(5)(6) of Definition \ref{DefinitionOfUnifiedSymmetryCondition}, denote
$$\LeftList
\begin{aligned}
  d'_{1L}&\DenoteBy\frac{1}{{2\sqrt 2}}c_{\mathfrak{D}-1}^{\mathfrak{D}- 3}(\rho_{(\mathfrak{D}- 4)(\mathfrak{D}- 2)} +\rho_{(\mathfrak{D}- 2)(\mathfrak{D}- 4)})+\frac{1}{{2\sqrt 2}}c_{\mathfrak{D}- 2}^{\mathfrak{D}- 3}(\rho_{(\mathfrak{D}- 4)(\mathfrak{D}-1)} +\rho_{(\mathfrak{D}-1)(\mathfrak{D}- 4)})\\
&+\frac{1}{{2\sqrt 2}}c_{\mathfrak{D}-1}^{\mathfrak{D}- 4}(\rho_{(\mathfrak{D}- 3)(\mathfrak{D}- 2)} +\rho_{(\mathfrak{D}- 2)(\mathfrak{D}- 3)})+\frac{1}{{2\sqrt 2}}c_{\mathfrak{D}- 2}^{\mathfrak{D}- 4}(\rho_{(\mathfrak{D}- 3)(\mathfrak{D}-1)} +\rho_{(\mathfrak{D}-1)(\mathfrak{D}- 3)}),\\
\end{aligned}\hfill\\
\RightList$$
$$\LeftList
\begin{aligned}
  d'_{2L}&\DenoteBy\frac{1}{{2\sqrt 2}}c_\mathfrak{D}^{\mathfrak{D}- 3}(\rho_{(\mathfrak{D}- 4)(\mathfrak{D}-1)} +\rho_{(\mathfrak{D}-1)(\mathfrak{D}- 4)})+\frac{1}{{2\sqrt 2}}c_{\mathfrak{D}-1}^{\mathfrak{D}- 3}(\rho_{(\mathfrak{D}- 4)\mathfrak{D}} +\rho_{\mathfrak{D}(\mathfrak{D}- 4)})\\
&+\frac{1}{{2\sqrt 2}}c_\mathfrak{D}^{\mathfrak{D}- 4}(\rho_{(\mathfrak{D}- 3)(\mathfrak{D}-1)} +\rho_{(\mathfrak{D}-1)(\mathfrak{D}- 3)})+\frac{1}{{2\sqrt 2}}c_{\mathfrak{D}-1}^{\mathfrak{D}- 4}(\rho_{(\mathfrak{D}- 3)\mathfrak{D}} +\rho_{\mathfrak{D}(\mathfrak{D}- 3)}),\\
\end{aligned}\hfill\\
\begin{aligned}
  d'_{3L}&\DenoteBy\frac{1}{{2\sqrt 2}}c_{\mathfrak{D}- 2}^{\mathfrak{D}- 3}(\rho_{(\mathfrak{D}- 4)\mathfrak{D}} +\rho_{\mathfrak{D}(\mathfrak{D}- 4)})+\frac{1}{{2\sqrt 2}}c_\mathfrak{D}^{\mathfrak{D}- 3}(\rho_{(\mathfrak{D}- 4)(\mathfrak{D}- 2)} +\rho_{(\mathfrak{D}- 2)(\mathfrak{D}- 4)})\\
&+\frac{1}{{2\sqrt 2}}c_{\mathfrak{D}- 2}^{\mathfrak{D}- 4}(\rho_{(\mathfrak{D}- 3)\mathfrak{D}} +\rho_{\mathfrak{D}(\mathfrak{D}- 3)})+\frac{1}{{2\sqrt 2}}c_\mathfrak{D}^{\mathfrak{D}- 4}(\rho_{(\mathfrak{D}- 3)(\mathfrak{D}- 2)} +\rho_{(\mathfrak{D}- 2)(\mathfrak{D}- 3)}),\\ 
\end{aligned}\hfill\\
\begin{aligned}
  u'_{1L}&\DenoteBy\frac{1}{{2\sqrt 2}}c_{\mathfrak{D}- 2}^{\mathfrak{D}- 3}(\rho_{(\mathfrak{D}- 4)(\mathfrak{D}- 2)} +\rho_{(\mathfrak{D}- 2)(\mathfrak{D}- 4)})+\frac{1}{{2\sqrt 2}}c_{\mathfrak{D}- 2}^{\mathfrak{D}- 4}(\rho_{(\mathfrak{D}- 3)(\mathfrak{D}- 2)} +\rho_{(\mathfrak{D}- 2)(\mathfrak{D}- 3)})\hfill\\
&+\frac{1}{{2\sqrt 2}}c_{\mathfrak{D}-1}^{\mathfrak{D}- 3}(\rho_{(\mathfrak{D}- 4)(\mathfrak{D}-1)} +\rho_{(\mathfrak{D}-1)(\mathfrak{D}- 4)})+\frac{1}{{2\sqrt 2}}c_{\mathfrak{D}-1}^{\mathfrak{D}- 4}(\rho_{(\mathfrak{D}- 3)(\mathfrak{D}-1)} +\rho_{(\mathfrak{D}-1)(\mathfrak{D}- 3)}),\hfill\\
\end{aligned}\hfill\\
\begin{aligned}
  u'_{2L}&\DenoteBy\frac{1}{{2\sqrt 2}}c_{\mathfrak{D}-1}^{\mathfrak{D}- 3}(\rho_{(\mathfrak{D}- 4)(\mathfrak{D}-1)} +\rho_{(\mathfrak{D}-1)(\mathfrak{D}- 4)})+\frac{1}{{2\sqrt 2}}c_{\mathfrak{D}-1}^{\mathfrak{D}- 4}(\rho_{(\mathfrak{D}- 3)(\mathfrak{D}-1)} +\rho_{(\mathfrak{D}-1)(\mathfrak{D}- 3)})\hfill\\
&+\frac{1}{{2\sqrt 2}}c_\mathfrak{D}^{\mathfrak{D}- 3}(\rho_{(\mathfrak{D}- 4)\mathfrak{D}} +\rho_{\mathfrak{D}(\mathfrak{D}- 4)})+\frac{1}{{2\sqrt 2}}c_\mathfrak{D}^{\mathfrak{D}- 4}(\rho_{(\mathfrak{D}- 3)\mathfrak{D}} +\rho_{\mathfrak{D}(\mathfrak{D}- 3)}),\hfill\\
\end{aligned}\hfill\\
\begin{aligned}
  u'_{3L}&\DenoteBy\frac{1}{{2\sqrt 2}}c_\mathfrak{D}^{\mathfrak{D}- 3}(\rho_{(\mathfrak{D}- 4)\mathfrak{D}} +\rho_{\mathfrak{D}(\mathfrak{D}- 4)})+\frac{1}{{2\sqrt 2}}c_\mathfrak{D}^{\mathfrak{D}- 4}(\rho_{(\mathfrak{D}- 3)\mathfrak{D}} +\rho_{\mathfrak{D}(\mathfrak{D}- 3)})\hfill\\
&+\frac{1}{{2\sqrt 2}}c_{\mathfrak{D}- 2}^{\mathfrak{D}- 3}(\rho_{(\mathfrak{D}- 4)(\mathfrak{D}- 2)} +\rho_{(\mathfrak{D}- 2)(\mathfrak{D}- 4)})+\frac{1}{{2\sqrt 2}}c_{\mathfrak{D}- 2}^{\mathfrak{D}- 4}(\rho_{(\mathfrak{D}- 3)(\mathfrak{D}- 2)} +\rho_{(\mathfrak{D}- 2)(\mathfrak{D}- 3)}).\hfill\\ 
\end{aligned}\hfill\\
\RightList$$
Then the geometric properties $d_1$, $d_2$, $d_3$, $u_1$, $u_2$, $u_3$ of $\mathcal{F}$ satisfy the following conclusions on $(M,\mathcal{G})$.
$$\LeftList
\begin{aligned}
  d_{1L;P}&=\partial_P d_{1L}  -g_s d_{1L} U^1_P  +g_s d_{2L} V^1_P  -g_s d_{3L} V^1_P \\
&-g_s u_{1L} X^{23}_P  -\frac{g_s}{2}u_{2L} X^{31}_P  +\frac{g_s}{2}u_{2L} Y^{31}_P  -\frac{g_s}{2}u_{3L} X^{12}_P  -\frac{g_s}{2}u_{3L} Y^{12}_P  -g u'_{1L}W^1_P, \\
\end{aligned}\hfill\\
\begin{aligned}
  d_{2L;P}&=\partial_P d_{2L} -g_s d_{2L} U^2_P  +g_s d_{3L} V^2_P  -g_s d_{1L} V^2_P \\
&-g_s u_{2L} X^{31}_P  -\frac{g_s}{2}u_{3L} X^{12}_P  +\frac{g_s}{2}u_{3L} Y^{12}_P  -\frac{g_s}{2}u_{1L} X^{23}_P  -\frac{g_s}{2}u_{1L} Y^{23}_P  -g u'_{2L}W^1_P, \\
\end{aligned}\hfill\\
\begin{aligned}
  d_{3L;P}&=\partial_P d_{3L}  -g_s d_{3L} U^3_P  +g_s d_{1L} V^3_P  -g_s d_{2L} V^3_P \\
&-g_s u_{3L} X^{12}_P  -\frac{g_s}{2}u_{1L} X^{23}_P  +\frac{g_s}{2}u_{1L} Y^{23}_P  -\frac{g_s}{2}u_{2L} X^{31}_P  -\frac{g_s}{2}u_{2L} Y^{31}_P  -g u'_{3L}W^1_P, \\ 
\end{aligned}\hfill\\
\begin{aligned}
  d_{1R;P}&=\partial_P d_{1R} -g_s d_{1L}V^1_P  +g_s d_{2L}U^1_P  -g_s d_{3L}U^1_P \\
&+g_s u_{1L}Y^{23}_P  +\frac{g_s}{2}u_{2L}X^{31}_P  -\frac{g_s}{2}u_{2L}Y^{31}_P  -\frac{g_s}{2}u_{3L}X^{12}_P  -\frac{g_s}{2}u_{3L}Y^{12}_P, \\
  d_{2R;P}&=\partial_P d_{2R} -g_s d_{2L}V^2_P  +g_s d_{3L}U^2_P  -g_s d_{1L}U^2_P \\
&+g_s u_{2L}Y^{31}_P  +\frac{g_s}{2}u_{3L}X^{12}_P  -\frac{g_s}{2}u_{3L}Y^{12}_P  -\frac{g_s}{2}u_{1L}X^{23}_P  -\frac{g_s}{2}u_{1L}Y^{23}_P, \\
  d_{3R;P}&=\partial_P d_{3R} -g_s d_{3L}V^3_P  +g_s d_{1L}U^3_P  -g_s d_{2L}U^3_P \\
&+g_s u_{3L}Y^{12}_P  +\frac{g_s}{2}u_{1L}X^{23}_P  -\frac{g_s}{2}u_{1L}Y^{23}_P  -\frac{g_s}{2}u_{2L}X^{31}_P  -\frac{g_s}{2}u_{2L}Y^{31}_P, \\ 
\end{aligned}\hfill\\
\begin{aligned}
  u_{1L;P}&=\partial_P u_{1L} -g_s u_{1L}U^1_P  -\frac{g_s}{2}u_{2L}X^{12}_P  -\frac{g_s}{2}u_{2L}Y^{12}_P  -\frac{g_s}{2}u_{3L}X^{31}_P  +\frac{g_s}{2}u_{3L}Y^{31}_P \hfill\\
&-g_s d_{1L}X^{23}_P  +g_s d_{2L}Y^{23}_P  -g_s d_{3L}Y^{23}_P  -g d'_{1L}W^1_P, \hfill\\
\end{aligned}\hfill\\
\begin{aligned}
  u_{2L;P}&=\partial_P u_{2L} -g_s u_{2L}U^2_P  -\frac{g_s}{2}u_{3L}X^{23}_P  -\frac{g_s}{2}u_{3L}Y^{23}_P  -\frac{g_s}{2}u_{1L}X^{12}_P  +\frac{g_s}{2}u_{1L}Y^{12}_P \hfill\\
&-g_s d_{2L}X^{31}_P  +g_s d_{3L}Y^{31}_P  -g_s d_{1L}Y^{31}_P  -g d'_{2L}W^1_P, \hfill\\
\end{aligned}\hfill\\
\begin{aligned}
  u_{3L;P}&=\partial_P u_{3L} -g_s u_{3L}U^3_P  -\frac{g_s}{2}u_{1L}X^{31}_P  -\frac{g_s}{2}u_{1L}Y^{31}_P  -\frac{g_s}{2}u_{2L}X^{23}_P  +\frac{g_s}{2}u_{2L}Y^{23}_P \hfill\\
&-g_s d_{3L}X^{12}_P  +g_s d_{1L}Y^{12}_P  -g_s d_{2L}Y^{12}_P  -g d'_{3L}W^1_P, \hfill\\ 
\end{aligned}\hfill\\
\begin{aligned}
  u_{1R;P}&=\partial_P u_{1R} -g_s u_{1R}U^1_P  +\frac{g_s}{2}u_{2R}X^{12}_P  +\frac{g_s}{2}u_{2R}Y^{12}_P  +\frac{g_s}{2}u_{3R}X^{31}_P  -\frac{g_s}{2}u_{3R}Y^{31}_P, \\
\end{aligned}\hfill\\
\begin{aligned}
  u_{2R;P}&=\partial_P u_{2R} -g_s u_{2R}U^2_P  +\frac{g_s}{2}u_{3R}X^{23}_P  +\frac{g_s}{2}u_{3R}Y^{23}_P  +\frac{g_s}{2}u_{1R}X^{12}_P  -\frac{g_s}{2}u_{1R}Y^{12}_P, \\
\end{aligned}\hfill\\
\begin{aligned}
  u_{3R;P}&=\partial_P u_{3R} -g_s u_{3R}U^3_P  +\frac{g_s}{2}u_{1R}X^{31}_P  +\frac{g_s}{2}u_{1R}Y^{31}_P  +\frac{g_s}{2}u_{2R}X^{23}_P  -\frac{g_s}{2}u_{2R}Y^{23}_P. \\ 
\end{aligned}\hfill\\
\RightList$$

\sProof Substitute Definition \ref{DefinitionOfWeakElectroStrongPotential} into $\rho_{mn}$ and consider Definition \ref{DefinitionOfUnifiedSymmetryCondition}, then compute them, and then substitute $d'_{1L},\ d'_{2L},\ d'_{3L},\ u'_{1L},\ u'_{2L},\ u'_{3L}$ into them, we finally obtain the results.
\qed

\sRemark The above proposition shows a geometric origin of CKM mixing. We see that, in affine connection representation of gauge fields, $d'_{1L},\ d'_{2L},\ d'_{3L},\ u'_{1L},\ u'_{2L},\ u'_{3L}$ arise as geometric properties on manifold. Detailed equations of CKM mixing can be obtained on an additional condition such as
$$\begin{aligned}
\rho_{(\mathfrak{D}- 2)(\mathfrak{D}- 2)} 
&= a^{23}\rho_{(\mathfrak{D}- 2)(\mathfrak{D}- 3)} + a^{24}\rho_{(\mathfrak{D}- 2)(\mathfrak{D}- 4)} 
+ a^{32}\rho_{(\mathfrak{D}- 3)(\mathfrak{D}- 2)} + a^{42}\rho_{(\mathfrak{D}- 4)(\mathfrak{D}- 2)}, \\
\rho_{(\mathfrak{D}- 1)(\mathfrak{D}- 1)}
&= a^{13}\rho_{(\mathfrak{D}- 1)(\mathfrak{D}- 3)} + a^{14}\rho_{(\mathfrak{D}- 1)(\mathfrak{D}- 4)} 
+ a^{31}\rho_{(\mathfrak{D}- 3)(\mathfrak{D}- 1)} + a^{41}\rho_{(\mathfrak{D}- 4)(\mathfrak{D}- 1)}, \\
\rho_{\mathfrak{D}\mathfrak{D}}\ \ \ \ \ \ \ \ \ \ \ \ \ 
&= a^{03}\rho_{\mathfrak{D}(\mathfrak{D}- 3)} + a^{04}\rho_{\mathfrak{D}(\mathfrak{D}- 4)} 
+ a^{30}\rho_{(\mathfrak{D}- 3)\mathfrak{D}} + a^{40}\rho_{(\mathfrak{D}- 4)\mathfrak{D}}, \\
\rho_{(\mathfrak{D}- 2)(\mathfrak{D}- 1)}
&= b^{23}\rho_{(\mathfrak{D}- 2)(\mathfrak{D}- 3)} + b^{13}\rho_{(\mathfrak{D}- 1)(\mathfrak{D}- 3)} + b^{24}\rho_{(\mathfrak{D}- 2)(\mathfrak{D}- 4)} + b^{14}\rho_{(\mathfrak{D}- 1)(\mathfrak{D}- 4)}, \\
\rho_{(\mathfrak{D}- 1)(\mathfrak{D}- 2)}
&= b^{32}\rho_{(\mathfrak{D}- 3)(\mathfrak{D}- 2)} + b^{31}\rho_{(\mathfrak{D}- 3)(\mathfrak{D}- 1)} + b^{42}\rho_{(\mathfrak{D}- 4)(\mathfrak{D}- 2)} + b^{41}\rho_{(\mathfrak{D}- 4)(\mathfrak{D}- 1)}, \\
\rho_{(\mathfrak{D}- 2)\mathfrak{D}}\ \ \ \ \ \ 
&= b^{23}\rho_{(\mathfrak{D}- 2)(\mathfrak{D}- 3)} + b^{03}\rho_{\mathfrak{D}(\mathfrak{D}- 3)} + b^{24}\rho_{(\mathfrak{D}- 2)(\mathfrak{D}- 4)} + b^{04}\rho_{\mathfrak{D}(\mathfrak{D}- 4)}, \\
\rho_{\mathfrak{D}(\mathfrak{D}- 2)}\ \ \ \ \ \ 
&= b^{32}\rho_{(\mathfrak{D}- 3)(\mathfrak{D}- 2)} + b^{30}\rho_{(\mathfrak{D}- 3)\mathfrak{D}} + b^{42}\rho_{(\mathfrak{D}- 4)(\mathfrak{D}- 2)} + b^{40}\rho_{(\mathfrak{D}- 4)\mathfrak{D}}, \\
\rho_{(\mathfrak{D}- 1)\mathfrak{D}}\ \ \ \ \ \ 
&= b^{13}\rho_{(\mathfrak{D}- 1)(\mathfrak{D}- 3)} + b^{03}\rho_{\mathfrak{D}(\mathfrak{D}- 3)} + b^{14}\rho_{(\mathfrak{D}- 1)(\mathfrak{D}- 4)} + b^{04}\rho_{\mathfrak{D}(\mathfrak{D}- 4)}, \\
\rho_{\mathfrak{D}(\mathfrak{D}- 1)}\ \ \ \ \ \ 
&= b^{31}\rho_{(\mathfrak{D}- 3)(\mathfrak{D}- 1)} + b^{30}\rho_{(\mathfrak{D}- 3)\mathfrak{D}} + b^{41}\rho_{(\mathfrak{D}- 4)(\mathfrak{D}- 1)} + b^{40}\rho_{(\mathfrak{D}- 4)\mathfrak{D}} \\
\end{aligned}$$

\sDefinition\label{DefinitionOfLeptonAndHardron} A particle is not an existence at the place of an individual point, its concept is defined on the entire manifold. Concretely speaking, if the reference-system $\mathcal{F}$ satisfies 
$$\begin{aligned}
&\begin{aligned}
\rho_{(\mathfrak{D}- 2)(\mathfrak{D}- 2)}&=\rho_{(\mathfrak{D}-1)(\mathfrak{D}-1)}=\rho_{\mathfrak{D}\mathfrak{D}}=\rho_{(\mathfrak{D}- 2)(\mathfrak{D}-1)}=\rho_{(\mathfrak{D}-1)(\mathfrak{D}- 2)}=\rho_{(\mathfrak{D}-1)\mathfrak{D}}=\rho_{\mathfrak{D}(\mathfrak{D}-1)}=\rho_{\mathfrak{D}(\mathfrak{D}- 2)} \\
&=\rho_{(\mathfrak{D}- 2)\mathfrak{D}}=0,\\
\end{aligned}\\
&\begin{aligned}
\Gamma_{(\mathfrak{D}- 2)(\mathfrak{D}- 2)P}&=\Gamma_{(\mathfrak{D}-1)(\mathfrak{D}-1)P}=\Gamma_{\mathfrak{D}\mathfrak{D}P}=\Gamma_{(\mathfrak{D}- 2)(\mathfrak{D}-1)P}=\Gamma_{(\mathfrak{D}-1)(\mathfrak{D}- 2)P}=\Gamma_{(\mathfrak{D}-1)\mathfrak{D}P}=\Gamma_{\mathfrak{D}(\mathfrak{D}-1)P} \\
&=\Gamma_{\mathfrak{D}(\mathfrak{D}- 2)P}=\Gamma_{(\mathfrak{D}- 2)\mathfrak{D}P}=0,\\
\end{aligned}
\end{aligned}$$
we say $\mathcal{F}$ is a {\bf lepton}, otherwise $\mathcal{F}$ is a {\bf hadron}. 

Suppose $\mathcal{F}$ is a hadron. For $d_1$, $d_2$, $d_3$, $u_1$, $u_2$, $u_3$, if $\mathcal{F}$ satisfies that five of them are zero and the other one is non-zero, we say $\mathcal{F}$ is an {\bf individual quark}. 

\sProposition\label{PropositionOfColorConfinement} There does not exist an individual quark. In other words, if any five ones of $d_1, d_2, d_3, u_1, u_2, u_3$ are zero, then $d_1 =d_2 =d_3 =u_1 =u_2 =u_3 =0$.

For an individual down-type quark, the above proposition is evidently true. Without loss of generality, let $u_1 =u_2 =u_3 =0$ and $d_1 =d_2 =0$; thus, $\rho_{(\mathfrak{D}- 2)(\mathfrak{D}- 2)}=\rho_{(\mathfrak{D}-1)(\mathfrak{D}-1)}=\rho_{\mathfrak{D}\mathfrak{D}}=0$; hence, we must have $d_3 =0$.

For an individual up-type quark, this paper has not made progress on the proof yet. Nevertheless, Proposition \ref{PropositionOfColorConfinement} provides the color confinement with a new geometric interpretation, which is significant in itself. It involves a natural geometric constraint of the curvatures among different dimensions.

\SectionWithLabel{Conclusions}

1. An affine connection representation of gauge fields is established in this paper. It has the following main points of view. 

(i) The holonomic connection Eq.(\ref{AffineConnectionDefinitionA}) contains more geometric information than Levi-Civita connection. It can uniformly describe gauge field and gravitational field.

(ii) Time is the total spatial metric with respect to all dimensions of internal coordinate space and external coordinate space. 

(iii) Energy is the total momentum with respect to all dimensions of internal coordinate space and external coordinate space. 

(iv) On-shell evolution is described by gradient direction. 

(v) Quantum theory is a geometric theory of distribution of gradient directions. It has a geometric meaning discussed in section \ref{Quantum evolution as a distribution of gradient directions}. 

\noindent 2. In the affine connection representation of gauge fields, some physical objects are incorporated into the same geometric framework.

(i) Gauge field and gravitational field can both be represented by affine connection. They have a unified coordinate description. Some parts of $\Gamma_{NP}^{M}$ describe gauge fields such as electromagnetic, weak, and strong interaction fields. The other parts of $\Gamma_{NP}^{M}$ describe gravitational field. 

(ii) Gauge field and elementary particle field are both geometric entities constructed from semi-metric. The components $\rho_{mn}$ of $\rho_{MN}$ with $m,n\in \{4,5,\cdots,\mathfrak{D}\}$ describe leptons and quarks, the other components of $\rho_{MN}$ may describe particle fields of dark matters.

(iii) Physical evolutions of gauge field and elementary particle field have a unified geometric description. Their on-shell evolution and quantum evolution both present as geometric properties about gradient direction. 

(iv) CPT inversion can be geometrically interpreted as a joint transformation of full inversion of coordinates and full inversion of metrics.

(v) Rest-mass is the total momentum with respect to internal space. It originates from geometric property of internal space. Energy, momentum, and mass have no essential difference in geometric sense. 

(vi) Quantum theory and gravitational theory have a unified geometric interpretation and the same view of time and space. They both reflect intrinsic geometric properties of manifold.

(vii) The origination of coupling constants of interactions can be interpreted geometrically.

(viii) Chiral asymmetry, PMNS mixing, and CKM mixing arise as geometric properties on manifold. 

(ix) There exists a geometric interpretation to the color confinement of quarks.

\noindent In the affine connection representation, we can get better interpretations to these physical properties. Therefore, to represent gauge fields by affine connection is probably a necessary step towards the ultimate theory of physics.

\small
\bibliographystyle{spmpsci} 
\bibliography{ref.bib}




\end{document}